\documentclass[a4paper,twoside,11pt,english,fleqn,intlimits]{article}

\usepackage{amsmath,amsthm,amsfonts,amssymb}
\usepackage{graphicx}
\usepackage[dvips]{epsfig}
\usepackage{color} 
\usepackage{babel}
\usepackage[T1]{fontenc}
\usepackage{fancyhdr,times,euler,euscript,eufrak,color}
\usepackage{titlesec}
\usepackage{url,index,float}
\usepackage[all]{xy}

\DeclareMathAlphabet{\mathscr}{T1}{pzc}{m}{it}

\titleformat{\section}[block]{\scshape\filcenter\Large}{\thesection.}{.5em}{}
\titleformat{\subsection}[block]{\bfseries\filcenter\large}{\thesubsection.}{.5em}{}
\titleformat{\subsubsection}[runin]{\bfseries}{\thesubsubsection.}{.5em}{}[.]

\swapnumbers

\theoremstyle{definition}

\newcommand{\emptysectionmark}[1]{\markboth{\textbf{#1}}{\textbf{#1}}}

\fancypagestyle{plain}{
\fancyhf{}\fancyfoot[LC,RC]{}
\fancyfoot[LE,RO]{$\thepage$}

}

\UseTips
\SelectTips{eu}{11}

\newdir{ >}{{}*!/-10pt/@{>}}
\newdir{ -}{{}*!/-10pt/@{}}
\newdir{> }{{}*!/+10pt/@{>}}

\makeatletter

\xyletcsnamecsname@{dir4{}}{dir{}}
\xydefcsname@{dir4{-}}{\line@ \quadruple@\xydashh@}
\xydefcsname@{dir4{.}}{\point@ \quadruple@\xydashh@}
\xydefcsname@{dir4{~}}{\squiggle@ \quadruple@\xybsqlh@}
\xydefcsname@{dir4{>}}{\Tttip@}
\xydefcsname@{dir4{<}}{\reverseDirection@\Tttip@}

\xydef@\quadruple@#1{%
	\edef\Drop@@{%
		\dimen@=#1\relax
		\dimen@=.5\dimen@
		\A@=-\sinDirection\dimen@
		\B@=\cosDirection\dimen@
		\setboxz@h{%
			\setbox2=\hbox{\kern3\A@\raise3\B@\copy\z@}%
			\dp2=\z@ \ht2=\z@ \wd2=\z@ \box2
			\setbox2=\hbox{\kern\A@\raise\B@\copy\z@}%
			\dp2=\z@ \ht2=\z@ \wd2=\z@ \box2
			\setbox2=\hbox{\kern-\A@\raise-\B@\copy\z@}%
			\dp2=\z@ \ht2=\z@ \wd2=\z@ \box2
			\setbox2=\hbox{\kern-3\A@\raise-3\B@ \noexpand\boxz@}%
			\dp2=\z@ \ht2=\z@ \wd2=\z@ \box2
		}%
		\ht\z@=\z@ \dp\z@=\z@ \wd\z@=\z@ \noexpand\styledboxz@
	}%
}

\xydef@\Tttip@{\kern2pt \vrule height2pt depth2pt width\z@
	\Tttip@@ \kern2pt \egroup
	\U@c=0pt \D@c=0pt \L@c=0pt \R@c=0pt \Edge@c={\circleEdge}%
	\def\Leftness@{.5}\def\Upness@{.5}%
	\def\Drop@@{\styledboxz@}\def\Connect@@{\straight@{\dottedSpread@\jot}}}
	
\xydef@\Tttip@@{%
	\dimen@=.25\dimen@
 	\B@=\cosDirection\dimen@
	\setboxz@h\bgroup\reverseDirection@\line@ \wdz@=\z@ \ht\z@=\z@ \dp\z@=\z@
	{\vDirection@(1,-1)\xydashl@ \xyatipfont\char\DirectionChar}%
	{\vDirection@(1,+1)\xydashl@ \xybtipfont\char\DirectionChar}%
}

\xydef@\ar@form{
	\ifx \space@\next \expandafter\DN@\space{\xyFN@\ar@form}%
	\else\ifx ^\next \DN@ ^{\xyFN@\ar@style}\edef\arvariant@@{\string^}%
	\else\ifx _\next \DN@ _{\xyFN@\ar@style}\edef\arvariant@@{\string_}%
	\else\ifx 0\next \DN@ 0{\xyFN@\ar@style}\def\arvariant@@{0}%
	\else\ifx 1\next \DN@ 1{\xyFN@\ar@style}\def\arvariant@@{1}%
	\else\ifx 2\next \DN@ 2{\xyFN@\ar@style}\def\arvariant@@{2}%
	\else\ifx 3\next \DN@ 3{\xyFN@\ar@style}\def\arvariant@@{3}%
	\else\ifx 4\next \DN@ 4{\xyFN@\ar@style}\def\arvariant@@{4}%
	\else\ifx \bgroup\next \let\next@=\ar@style
	\else\ifx [\next \DN@[##1]{\ar@modifiers{[##1]}}
	\else\ifx *\next \DN@ *{\ar@modifiers}%
	\else\addLT@\ifx\next \let\next@=\ar@slide
	\else\ifx /\next \let\next@=\ar@curveslash
	\else\ifx (\next \let\next@=\ar@curveinout 
	\else\addRQ@\ifx\next \addRQ@\DN@{\ar@curve@}%
	\else\addLQ@\ifx\next \addLQ@\DN@{\xyFN@\ar@curve}%
	\else\addDASH@\ifx\next \addDASH@\DN@{\defarstem@-\xyFN@\ar@}%
	\else\addEQ@\ifx\next \addEQ@\DN@{\def\arvariant@@{2}\defarstem@-\xyFN@\ar@}%
	\else\addDOT@\ifx\next \addDOT@\DN@{\defarstem@.\xyFN@\ar@}%
	\else\ifx :\next \DN@:{\def\arvariant@@{2}\defarstem@.\xyFN@\ar@}%
	\else\ifx ~\next \DN@~{\defarstem@~\xyFN@\ar@}%
	\else\ifx !\next \DN@!{\dasharstem@\xyFN@\ar@}%
	\else\ifx ?\next \DN@?{\ar@upsidedown\xyFN@\ar@}%
	\else \let\next@=\ar@error
	\fi\fi\fi\fi\fi\fi\fi\fi\fi\fi\fi\fi\fi\fi\fi\fi\fi\fi\fi\fi\fi\fi\fi \next@}

\makeatother

\newcommand{\fl}{\to}

\newcommand{\dfl}{\Rightarrow}
\newcommand{\tfl}{\Rrightarrow}
\newcommand{\Llongrightarrow}{\equiv\joinrel\Rrightarrow}

\newcommand{\ens}[1]{\left\{#1\right\}}
\newcommand{\mon}[1]{\langle#1\rangle}
\newcommand{\abs}[1]{\left|#1\right|}
\newcommand{\norm}[1]{\abs{\abs{#1}}}

\DeclareMathOperator{\id}{Id}

\newcommand{\tens}{\otimes}
\newcommand{\cl}[1]{\overline{#1}}

\newcommand{\ie}{\emph{i.e.}}

\renewcommand{\phi}{\varphi}
\renewcommand{\epsilon}{\varepsilon}

\newcommand{\Nb}{\mathbb{N}}

\newcommand{\Zb}{\mathbb{Z}}

\newcommand{\dr}{\partial}

\newcommand{\Cr}{\EuScript{C}}
\newcommand{\Dr}{\EuScript{D}}
\newcommand{\Er}{\EuScript{E}}

\renewcommand{\Pr}{\EuScript{P}}

\newcommand{\Sr}{\EuScript{S}}

\newcommand{\Tr}{\EuScript{T}}

\newcommand{\Vr}{\EuScript{V}}

\newcommand{\commute}[2]{\ar@{}[#1]|-*+[o][F-]{\scriptscriptstyle{#2}}}

\newcommand{\Ab}{\mathbf{Ab}}

\newcommand{\Set}{\mathbf{Set}}

\newcommand{\Ord}{\mathbf{Ord}}
\newcommand{\Cat}[1]{\mathbf{Cat}_{#1}}
\newcommand{\Pol}[1]{\mathbf{Pol}_{#1}}
\newcommand{\Tck}[1]{\mathbf{Tck}_{#1}}
\newcommand{\Ct}[1]{\mathbf{C}{#1}}

\newcommand{\whisk}[1]{\mathbf{W}{#1}}

\newcommand{\figeps}[1]{\raisebox{-1.25mm}{\includegraphics{#1.eps}}}
\newcommand{\smallfigeps}[1]{\includegraphics[scale=0.66]{#1.eps}}
\newcommand{\scalefigeps}[2]{\raisebox{#1}{\includegraphics{#2.eps}}}

\newcommand{\onecell}[3]{%
\xymatrix@1{
{\scriptstyle #2}
	\ar@<-0.2ex> [r] ^-{#1}
& {\scriptstyle #3}
}%
}
\newcommand{\twocell}[5]{%
\xymatrix@1{
{\scriptstyle #4}
	\ar@/^2ex/ [r] ^{#2} _{}="source"
	\ar@/_2ex/ [r] _{#3} ^{}="target"
	\ar@2 "source" ; "target" ^-{#1}
& {\scriptstyle #5}
}%
}
\newcommand{\threecell}[7]{%
\xymatrix@1{
{\scriptstyle #6}
	\ar@/^2ex/ [r] ^{#4} _{}="source"
	\ar@/_2ex/ [r] _{#5} ^{}="target"
	\ar@2 "source" ; "target" ^-{#2} 
& {\scriptstyle #7}
}
\:\overset{#1}{\tfl}\:
\xymatrix@1{
{\scriptstyle #6}
	\ar@/^2ex/ [r] ^{#4} _{}="source"
	\ar@/_2ex/ [r] _{#5} ^{}="target"
	\ar@2 "source" ; "target" ^-{#3} 
& {\scriptstyle #7}
}%
}

\definecolor{orange}{rgb}{1,0.55,0}

\newcommand{\pdf}[1]{$#1$}

\begin{document}

\thispagestyle{empty}
\begin{center}
\textbf{\LARGE Higher-dimensional categories with finite derivation type} 

\vspace{4mm}
\begin{tabular}{c c c}
\textbf{\large Yves Guiraud} &$\qquad\qquad$& \textbf{\large Philippe Malbos} \\
INRIA Nancy && Université Lyon 1 \\
yves.guiraud@loria.fr && malbos@math.univ-lyon1.fr
\end{tabular}
\end{center} 

\vspace{2mm}
\begin{em}
\hrule height 1.5pt

\medskip
\noindent \textbf{Abstract --} 
We study convergent (terminating and confluent) presentations of n-categories. Using the notion of polygraph (or computad), we introduce the homotopical property of finite derivation type for n-categories, generalising the one introduced by Squier for word rewriting systems. We characterise this property by using the notion of critical branching. In particular, we define sufficient conditions for an n-category to have finite derivation type. Through examples, we present several techniques based on derivations of 2-categories to study convergent presentations by 3-polygraphs. 


\noindent \textbf{Support --} This work has been partially
supported by ANR INVAL project (ANR-05-BLAN-0267). 

\medskip
\hrule height 1.5pt
\end{em}

\section*{Introduction}
\emptysectionmark{Introduction}

\subsubsection*{Rewriting}

This is a combinatorial theory that studies presentations by generators and relations. For that, the latter are replaced by \emph{rewriting rules}, which are relations only  usable in one direction~\cite{Newman42}. There exist many flavours of rewriting, depending on the objects to be presented: word rewriting~\cite{BookOtto93}, for monoids; term rewriting~\cite{BaaderNipkow98,Klop92,Terese03}, for algebraic theories~\cite{Lawvere63}; rewriting on topological objects, such as Reidemeister moves, for braids and knots~\cite{Adams04}.

In this work, we study presentations by rewriting of higher-dimensional categories, which encompass the ones above~\cite{Burroni93,Lafont03,Guiraud04,Guiraud06jpaa}, plus many others, like Petri nets~\cite{Guiraud06tcs} or formal proofs of propositional calculus and linear logic~\cite{Guiraud06apal}.  

For example, the presentation of the monoid $\mon{ a \:|\: aa=a}$ by the word rewriting system $aa\fl a$ is interpreted as follows: the generator $a$ is a $1$-cell and the rewriting rule is a $2$-cell $aa\dfl a$ over the $1$-category freely generated by $a$. Similarly, the presentation of the associative theory by the term rewriting system $(x\cdot y)\cdot z \:\fl\: x\cdot(y\cdot z)$ becomes: the binary operation is treated as a $2$-cell~$\figeps{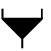}$, while the rewriting rule is seen as a $3$-cell over the $2$-category freely generated by~$\figeps{m}$, with shape
$$
\scalefigeps{-2.8mm}{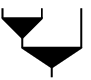} \:\tfl\: \scalefigeps{-2.8mm}{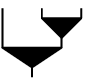} \:.
$$

\noindent Another example is the categorical presentation of the groups of permutations, used in particular for the explicit management of pointers in polygraphic programs~\cite{BonfanteGuiraud09}: it has one $2$-cell~$\figeps{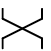}$, standing for a generating transposition, and the following two $3$-cells, respectively expressing that~$\figeps{tau}$ is an involution and that it satisfies the Yang-Baxter relation:
$$
\scalefigeps{-2.8mm}{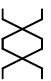}
	\:\tfl\: 
	\scalefigeps{-2.8mm}{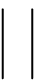}  
\qquad\text{and}\qquad\: 
\scalefigeps{-4.3mm}{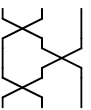} 
	\:\tfl\:
	\scalefigeps{-4.3mm}{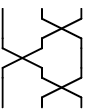} \:.
$$

\subsubsection*{Polygraphs}

The categorical rewriting systems presented in the previous paragraph are particular instances of objects called \emph{polygraphs} or \emph{computads}. Those objects are presentations by "generators" and "relations" of higher-dimensional categories~\cite{Street76,Burroni93,Street87,Street95} and they are defined by induction as follows. A $0$-polygraph is a set and a $1$-polygraph is a directed graph. An $(n+1)$-polygraph is given by an $n$-polygraph $\Sigma_n$, together with a family of $(n+1)$-cells between parallel $n$-cells of the $n$-category $\Sigma_n^*$ freely generated by $\Sigma_n$. The $n$-category presented by such an $n$-polygraph is the quotient of the free $n$-category $\Sigma_n^*$ by the congruence relation generated by the $(n+1)$-cells of $\Sigma_{n+1}$.

We recall the notions of polygraph and of presentation of $n$-categories in Section~\ref{Subsection:PolygraphsAndPresentations}, as originally described by Burroni~\cite{Burroni93,Metayer03}. Here we particularly focus on $n$-polygraphs for $n\leq 3$, because they contain well-known examples of rewriting systems: indeed, abstract rewriting systems, word rewriting systems and Petri nets are special instances of $1$-polygraphs, $2$-polygraphs and $3$-polygraphs, respectively, while term rewriting systems and formal proofs can be interpreted into $3$-polygraphs with similar computational properties.

Among those properties, we are mostly interested in \emph{convergence}: like other rewriting systems, a polygraph is \emph{convergent} when it is both \emph{terminating} and \emph{confluent}. The termination property ensures that no infinite reduction sequence exists, while the confluence property implies that all reduction sequences starting at the same point yield the same result. The aforegiven examples of $3$-polygraphs, for associativity and permutations, are convergent, as proved in Sections~\ref{Subsection:PolyMon} and~\ref{Subsection:PolyPerm}, respectively.

\subsubsection*{Homotopy type}

In order to study $n$-polygraphs from a homotopical point of view, we introduce the notion of \emph{higher-dimensional track category} in Section~\ref{Section:TDF}: a track $n$-category is an $(n-1)$-category enriched in groupoid (an $n$-category whose $n$-cells are invertible). This notion generalises track $2$-categories, introduced by Baues~\cite{Baues91} as an algebraic model of the homotopy type in dimension $2$.  

To an $n$-polygraph $\Sigma$, we associate the free track
$n$-category $\Sigma^{\top}$ it generates, used as a
combinatorial complex to describe the convergence property
of $\Sigma$.  Towards this goal, we define
in~\ref{SubSubsection:HomotopyBase} a \emph{homotopy
  relation} on $\Sigma^{\top}$ as a track $(n+1)$-category
with $\Sigma^{\top}$ as underlying $n$-category. Every
family of $(n+1)$-cells over $\Sigma^{\top}$ generates a
homotopy relation; a \emph{homotopy basis} of
$\Sigma^{\top}$ is such a family that generates a "full"
homotopy relation, \ie, a homotopy relation that identifies
any two parallel $n$-cells of $\Sigma^{\top}$.  

An $(n+1)$-polygraph $\Sigma$ has \emph{finite derivation type} when it is finite and when $\Sigma^{\top}$ admits a finite homotopy basis. This property is an invariant of the $n$-category being presented by $\Sigma$: when two $(n+1)$-polygraphs are \emph{Tietze-equivalent}, \ie, when they present the same $n$-category, then both or neither have finite derivation type (Proposition~\ref{TDFTietzeInvariant}). Hence, having finite derivation type is a finiteness property of $n$-categories in dimension $n+2$, in a way that is comparable to finite generation type (finiteness in dimension $n$) and finite presentation type (finiteness in dimension $n+1$). 

\subsubsection*{Critical branchings and homotopy bases} 

A \emph{critical branching} in a polygraph is a pair of reductions acting on overlapping "subcells" of the same cell (Definition~\ref{SubSubSection:DefinitionCriticalBranching}). The branching is  confluent when there exist two reduction sequences that close the diagram. For example, the $2$-polygraph $aa\dfl a$ has a unique, confluent critical branching:
$$
\xymatrix@ur{ 
aaa 
	\ar@2 [r] 
	\ar@2 [d]
& aa 
	\ar@2 [d]
\\
aa
	\ar@2 [r]
& a.
}
$$

\noindent The $3$-polygraph of associativity also has a unique, confluent critical branching, which is also known as the $2$-dimensional \emph{associahedron} or \emph{Stasheff polytope}:
$$
\xymatrix{
& {\figeps{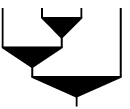}}
	\ar@3 [rr]
&& {\figeps{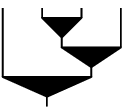}}
	\ar@3 [dr]
\\
{\figeps{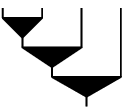}}
	\ar@3 [ur]
	\ar@3 [drr]
&&&& {\figeps{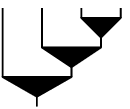}} \:.
\\
&& {\figeps{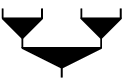}}
	\ar@3 [urr]
}
$$

\noindent Finally, the $3$-polygraph of permutations contains several critical branchings, given in~\ref{Subsubsection:ConfluencePerm}, all of which are confluent. Among them, one finds the $2$-dimensional \emph{permutohedron}, generated by an overlapping of the Yang-Baxter $3$-cell with itself: 
$$
\xymatrix{
& {\figeps{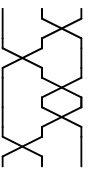}}
	\ar@3 [r] 
& {\figeps{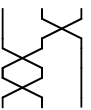}}
	\ar@3 [dr] 
\\
{\figeps{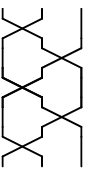} }
	\ar@3 [ur] 
	\ar@3 [dr] 
&&& {\figeps{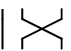}} \:.
\\
& {\figeps{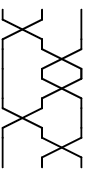}}
	\ar@3 [r] 
& {\figeps{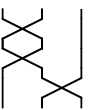}}
	\ar@3 [ur] 
}
$$

\noindent We prove that, when a polygraph is convergent, its critical branchings generate a homotopy basis (Proposition~\ref{Proposition3:BranchingHomotopyBases}). As a consequence, every finite and convergent polygraph with a finite number of critical branchings has finite derivation type (Proposition~\ref{PropfiniteCriticalBranchingFDT}).

This property is relevant when one considers higher-dimensional rewriting as a computational model, for example in the case of polygraphic programs~\cite{BonfanteGuiraud08,BonfanteGuiraud09}. Indeed, let us consider a convergent polygraph with finite derivation type: then, there exist finitely many elementary choices, corresponding to critical branchings, between parallel computation paths. Hence, Proposition~\ref{TDFTietzeInvariant} tells us that being of finite derivation type is a first step to ensure that an $n$-category admits a presentation by a rewriting system, together with a deterministic and finitely generated evaluation strategy. 

\subsubsection*{Convergence of \pdf{2}-polygraphs}

The notion of track $n$-category freely generated by an $n$-polygraph generalises the $2$-dimensional combinatorial complex associated to word rewriting systems~\cite{Squier94}. Squier introduced it to define finite derivation type for monoids and, then, linked this property with the possibility, for a finitely generated monoid, to have its word problem decided by the normal form algorithm. This procedure consists in finding a finite convergent presentation of the monoid $M$ by a word rewriting system $(X,R)$: given such a presentation, every element in the monoid $M$ has a canonical normal form in the free monoid $X^*$; hence, one can decide if $u$ and $v$ in $X^*$ represent the same element of $M$ by computing their unique normal forms for $R$ and, then, by checking if the results are equal or not in $X^*$. 

Squier has proved that, when a monoid admits a presentation by a finite and convergent word rewriting system, then it has finite derivation type. As a consequence, rewriting is not a universal way to decide the word problem of finitely generated monoids: to prove that, Squier has exhibited a finitely presented monoid whose word problem is decidable, yet lacking the property of finite derivation type.

Here, we recover Squier's convergence theorem as a consequence of Proposition~\ref{PropfiniteCriticalBranchingFDT}. Indeed, a $2$-polygraph has two kinds of critical branchings, 
namely \emph{inclusion} ones and \emph{overlapping} ones, respectively corresponding to the following situations:
$$
\xymatrix{
\strut
	\ar[r] 
	\ar@/^2pc/ [rrr] _{}="t1"
& \strut 
	\ar[r] ^{}="s1" _{}="s2"
	\ar@/_2pc/ [r] ^{}="t2"
& \strut
	\ar[r]
& \strut
\ar@2 "s1";"t1"
\ar@2 "s2";"t2"
}
\qquad\text{and}\qquad
\xymatrix{
\strut
	\ar [r]
	\ar@/^2pc/ [rr] _{}="t1"
	\ar@{} [rr] ^{}="s1"
& \strut
	\ar[r]
	\ar@/_2pc/ [rr] ^{}="t2"
	\ar@{} [rr] _{}="s2"
& \strut 
	\ar[r]
& \strut
\ar@2 "s1";"t1"
\ar@2 "s2";"t2"
} \:.
$$

\noindent Hence a finite $2$-polygraph can have only finitely many critical branchings, yielding a finite homotopy basis for its track $2$-category when it is also convergent.

\subsubsection*{Convergence of \pdf{3}-polygraphs}

This case is more complicated than the one of $2$-polygraphs, because of the nature of critical branchings generated by $3$-dimensional rewriting rules on $2$-cells. In Section~\ref{Section:CaseOf3Polygraphs}, we analyse the possible critical branchings a $3$-polygraph may have. We give a classification that unveil a new kind of these objects, that we call \emph{indexed critical branching} and that describes situations such as the following one:
$$
\scalebox{1.25}{\raisebox{-4mm}{\raisebox{-1.25mm}{\input{pc-ind-r-src.pstex_t}}}}
$$

\noindent where two $3$-cells respectively reduce the
$2$-cells
$$\scalebox{1.25}{\raisebox{-2.5mm}{\raisebox{-1.25mm}{\input{pc-ind-r-srcPart1.pstex_t}}}}
\qquad
\text{and}
\qquad
\scalebox{1.25}{\raisebox{-2.5mm}{\raisebox{-1.25mm}{\input{pc-ind-r-srcPart2.pstex_t}}}}\;.$$
There, the $2$-cell $k$ belongs to none of the considered $3$-cells. A \emph{normal instance} of the critical branching is such a situation where $k$ is a normal form (\ie, it cannot be reduced by any $3$-cell).

We prove that the existence of indexed critical branchings is an obstruction to get a generalisation of Squier's result on finiteness and convergence for higher dimensions. Indeed, for every natural number $n\geq 2$, there exists an $n$-category that lacks finite derivation type, even though it admits a presentation by a finite convergent $(n+1)$-polygraph (Theorem~\ref{Theorem:mainTheorem}). 

To get this result, we use the $3$-polygraph
$$
\scalefigeps{-2mm}{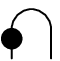} \tfl \scalefigeps{-2mm}{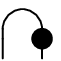} 
\:, \qquad
\scalefigeps{-2mm}{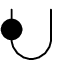} \tfl \scalefigeps{-2mm}{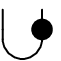}
\:, \qquad
\scalefigeps{-2mm}{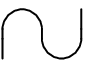} \tfl \scalefigeps{-2mm}{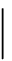} 
\:, \qquad
\scalefigeps{-2mm}{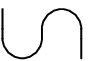} \tfl \scalefigeps{-2mm}{id1} 
\:,
$$

\noindent for which we prove, in Section~\ref{Subsection:main_counter_example}, that it is finite and convergent, but does not have finite derivation type. Let us note that this $3$-polygraph has a topological flavour: it presents a $2$-category whose $2$-cells are "planar necklaces with pearls" considered up to homotopy.

\subsubsection*{Finitely indexed \pdf{3}-polygraphs}

From our classification of critical branchings, we give a family of extra sufficient conditions that ensure that a finite convergent $3$-polygraph has finite derivation type. 

First, a finite convergent $3$-polygraph without indexed critical branching always has finite derivation type (Theorem~\ref{Thm:3PolyNonIndexe}): this is the case of the associativity one and of the monoid one. We illustrate the construction of a homotopy basis for this kind of $3$-polygraphs on this last example in Section~\ref{Subsection:PolyMon}: the basis corresponds to the coherence diagrams satisfied by a monoidal category. This yields a new formulation and proof of Mac Lane's coherence theorem asserting that, in a monoidal category, all the diagrams built from the monoidal structure are commutative~\cite{MacLane98}. 

More generally, we say that a $3$-polygraph is \emph{finitely indexed} when every indexed critical branching has finitely many normal instances (Definition~\ref{SubSubSection:FinitelyIndexed}). This is the case of the former class of non-indexed $3$-polygraphs, but also of many known ones such as the $3$-polygraph of permutations. We prove that a finite, convergent and finitely indexed $3$-polygraph has finite derivation type (Theorem~\ref{SubSubSection:MainTheorem1}).

In the case of finitely indexed $3$-polygraphs, building a homotopy basis requires a careful and comprehensive study of normal forms. We illustrate this construction in Section~\ref{Subsection:PolyPerm}, where we prove that the $3$-polygraph of permutations is finitely indexed. Such an observation was first made by Lafont~\cite{Lafont03} and we formalise it thanks to the notion of homotopy basis. 

\subsubsection*{Perspectives}

Our work gives methods to study, from a homotopical point of view, the convergence property of presentations of $2$-categories by $3$-polygraphs. We think that further research on these methods shall allow progress on questions such as the following ones. 

Our study of the $3$-polygraph of permutations adapts to polygraphic presentations of Lawvere algebraic theories~\cite{Lawvere63}. Indeed, there is a canonical translation of their presentations by term rewriting systems into $3$-polygraphs~\cite{Burroni93,Lafont03} and, when the original presentation is finite, left-linear and convergent, then the $3$-polygraph one gets is finite, convergent~\cite{Guiraud06jpaa} and finitely indexed~\cite{Lafont03}. Thus, if one proves that a given Lawvere algebraic theory does not have finite derivation type, one gets that it does not admit a presentation by a first-order functional program, which is a special kind of finite, left-linear and convergent term rewriting system.

We still do not know, for many special $2$-categories, if they admit a convergent presentation by a $3$-polygraph. Among these $2$-categories, we are particularly interested in the one of braids. It is known that it admits a presentation by a finite $3$-polygraph whose generators are, in dimension~$2$, the elementary crossings~$\figeps{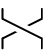}$ and~$\figeps{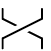}$ and, in dimension~$3$, the Reidemeister moves:
$$
\scalefigeps{-2.8mm}{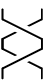} \:\tfl\:
\scalefigeps{-2.8mm}{id2} \:,
\qquad
\scalefigeps{-2.8mm}{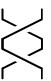} \:\tfl\:
\scalefigeps{-2.8mm}{id2} \:,
\qquad 
\scalefigeps{-4.3mm}{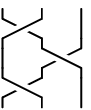} \:\tfl\:
\scalefigeps{-4.3mm}{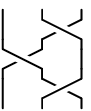} \:,
\qquad 
\scalefigeps{-4.3mm}{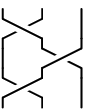} \:\tfl\:
\scalefigeps{-4.3mm}{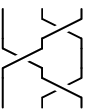} \:.
$$

\noindent As a consequence of this work, we know that the presence of indexed critical branchings in this $3$-polygraph, similar to the ones encountered for permutations, is one of the major obstructions to finding a convergent presentation of the $2$-category of braids. 

\setcounter{tocdepth}{1}
\tableofcontents

\bigskip\bigskip
\noindent In this work, we use known notions from the theories of categories, of $n$-categories and of rewriting that we do not necessarily explain in details. For more information on these subjects, we respectively recommend the books by Saunders Mac Lane~\cite{MacLane98}, by Eugenia Cheng and Aaron Lauda~\cite{ChengLauda04}, by Franz Baader and Tobias Nipkow~\cite{BaaderNipkow98}.

\section{Higher-dimensional categories presented by polygraphs}
\label{Section:nCategoriesPresentedByPolygraphs}

\subsection{Generalities on \pdf{n}-categories and \pdf{n}-functors}

In this document, we consider small, strict $n$-categories and strict $n$-functors between them. We denote by $\Cat{n}$ the (large) category they form. 

\subsubsection{Vocabulary and notations} 

If $\Cr$ is an $n$-category, we denote by $\Cr_k$ the set of $k$-cells of $\Cr$ and by $s_k$ and $t_k$ the $k$-source and $k$-target maps. If $f$ is a $k$-cell, $s_{k-1}(f)$ and $t_{k-1}(f)$ are respectively called its \emph{source} and \emph{target} and respectively denoted by $s(f)$ and $t(f)$. The source and target maps satisfy the \emph{globular relations}:
$$
s_k\circ s_{k+1} = s_k\circ t_{k+1} 
\qquad\text{and}\qquad
t_k\circ s_{k+1} = t_k\circ t_{k+1}.
$$

\noindent Two cells $f$ and $g$ are \emph{parallel} when they have same source and same target. A pair $(f,g)$ of parallel $k$-cells is called a \emph{$k$-sphere}. The \emph{boundary} of a $k$-cell is the $(k-1)$-sphere $\dr f=(s(f),t(f))$. The source and target maps are extended to a $k$-sphere $\gamma=(f,g)$ by $s(\gamma)=f$ and $t(\gamma)=g$.

A pair $(f,g)$ of $k$-cells of $\Cr$ is \emph{$i$-composable} when $t_i(f)=s_i(g)$ holds; when $i=k-1$, one simply says \emph{composable}. The $i$-composite of $(f,g)$ is denoted by $f\star_i g$, \ie, in the diagrammatic direction. The compositions satisfy the \emph{exchange relation} given, for every $j\neq k$ and every possible cells $f$, $f'$, $g$, $g'$, by:
$$
(f\star_j f')\star_k (g\star_j g') = (f\star_k g)\star_j(f'\star_k g').
$$

\noindent If $f$ is a $k$-cell, we denote by $1_f$ its identity $(k+1)$-cell and, by abuse, all the higher-dimensional identity cells it generates. When $1_f$ is composed with cells of dimension $k+1$ or higher, we abusively denote it by $f$ to make expressions easier to read. A cell is \emph{degenerate} when it is an identity cell. For $k\leq n$, a $k$-category $\Cr$ can be seen as an $n$-category, with only degenerate cells above dimension~$k$.

\subsubsection{Graphical representations}

Low-dimensional cells are written $u:p\fl q$, $f:u\dfl v$, $A:f\tfl g$ and pictured as usual (and so are $n$-categories, omitting the degenerate cells):
$$
\xymatrix{
{\scriptstyle p}
}
\qquad\qquad
\xymatrix{
{\scriptstyle p} \ar[rr] ^u && {\scriptstyle q}
}
\qquad\qquad
\xymatrix{
{\scriptstyle p}
	\ar@(ur,ul) [rr] ^u _{}="1"
	\ar@(dr,dl) [rr] _v ^{}="2"
	\ar@{=>} "1"; "2" ^ f
&& {\scriptstyle q}
}
\qquad\qquad
\xymatrix{
{\scriptstyle p}
	\ar@(ur,ul) [rr] ^u _{}="1"
	\ar@(dr,dl) [rr] _v ^{}="2"
	\ar@{=>}@/_2.5ex/ "1"; "2" _f ^{} ="3"
	\ar@{=>}@/^2.5ex/ "1"; "2" ^g _{} ="4"
	\ar@3{->} "3";"4" ^A
&& {\scriptstyle q}
}
$$

\noindent For readability, we also depict $3$-cells as "rewriting rules" on $2$-cells:
$$
\xymatrix@1{
{\scriptstyle p}
	\ar@(ur,ul) [rr] ^u _{}="1"
	\ar@(dr,dl) [rr] _v ^{}="2"
	\ar@{=>} "1"; "2" ^ f
&& {\scriptstyle q}
}
\quad\overset{A}{\Llongrightarrow}\quad
\xymatrix@1{
{\scriptstyle p}
	\ar@(ur,ul) [rr] ^u _{}="1"
	\ar@(dr,dl) [rr] _v ^{}="2"
	\ar@{=>} "1"; "2" ^ g
&& {\scriptstyle q}
}
$$

\noindent For $2$-cells, circuit-like diagrams are alternative representations, where $0$-cells are parts of the plane, $1$-cells are lines and $2$-cells are points, inflated for emphasis:
$$
\raisebox{-1.25mm}{\input{2-cellule.pstex_t}}
\qquad\qquad\qquad\qquad
\raisebox{-1.25mm}{\input{3-cellule.pstex_t}}
$$

\subsection{Standard cells and spheres}

\subsubsection{Suspension functors} 

For every natural number $n$, the \emph{suspension functor} 
$$
S_n:\Cat{n}\fl\Cat{n+1}
$$

\noindent lifts all the cells by one dimension, adding a formal $0$-source and a formal $0$-target for all of them; thus, in the $(n+1)$-category one gets, one has exactly the same compositions as in the original one. More formally, given an $n$-category $\Cr$, the $(n+1)$-category $S_n\Cr$ has the following cells:
$$
(S_n\Cr)_0 \:=\: \ens{-,+}
\qquad\text{and}\qquad
(S_n\Cr)_{k+1} \:=\: \Cr_k\amalg\ens{-,+}.
$$

\noindent Every cell has $0$-source $-$ and $0$-target $+$. The $(k+1)$-source and $(k+1)$-target of a non-degenerate cell are its $k$-source and $k$-target in $\Cr$. The $(k+1)$-composable pairs are the $k$-composable ones of $\Cr$, plus pairs where at least one of the cells is an identity of $-$ or $+$.

\subsubsection{Standard \pdf{n}-cells and \pdf{n}-spheres} 

By induction on $n$, we define the $n$-categories $\Er_n$
and $\Sr_n$, respectively called the \emph{standard
  $n$-cell} and the \emph{standard $n$-sphere}. We consider
them as the $n$-categorical equivalents of the 
standard topological $n$-ball and $n$-sphere, used to build
the $n$-categorical equivalents of (relative) CW-complexes. 

The standard $0$-cell $\Er_0$ is defined as any chosen single-element set and the standard $0$-sphere as any chosen set with two elements. Then, if $n\geq 1$, the $n$-categories $\Er_n$ and $\Sr_n$ are defined as the suspensions of $\Er_{n-1}$ and $\Sr_{n-1}$:
$$
\Er_n \:=\: S_{n-1}(\Er_{n-1})
\qquad\text{and}\qquad
\Sr_n \:=\: S_{n-1}(\Sr_{n-1}).
$$

\noindent For coherence, we define $\Sr_{-1}$ as the empty set. Thus, the standard $n$-cell~$\Er_n$ and $n$-sphere $\Sr_n$ have two non-degenerate $k$-cells $e_k^-$ and $e_k^+$ for every $k$ in $\ens{0,\dots,n-1}$, plus a non-degenerate $n$-cell $e_n$ in~$\Er_n$. Using the cellular representations, the standard cells $\Er_0$, $\Er_1$, $\Er_2$ and $\Er_3$ are respectively pictured as follows (for $\Sr_{-1}$, $\Sr_0$, $\Sr_1$ and $\Sr_2$, one removes the top-dimensional cell):
$$
\xymatrix{
{\scriptstyle e_0}
}
\qquad\qquad
\xymatrix{
{\scriptstyle e_0^-} \ar[rr]^-{e_1} && {\scriptstyle e_0^+}
}
\qquad\qquad
\xymatrix{
{\scriptstyle e_0^-}
	\ar@(ur,ul) [rr] ^-{e_1^-} _-{}="1"
	\ar@(dr,dl) [rr] _-{e_1^+} ^-{}="2"
	\ar@{=>} "1"; "2" ^-{e_2}
&& {\scriptstyle e_0^+}
}
\qquad\qquad
\xymatrix{
{\scriptstyle e_0^-}
	\ar@(ur,ul) [rr] ^-{e_1^-} _-{}="1"
	\ar@(dr,dl) [rr] _-{e_1^+} ^-{}="2"
	\ar@{=>}@/_2.5ex/ "1"; "2" _-{e_2^-} ^-{} ="3"
	\ar@{=>}@/^2.5ex/ "1"; "2" ^-{e_2^+} _-{} ="4"
	\ar@3{->} "3";"4" ^-{e_3}
&& {\scriptstyle e_0^+}
}
$$

\noindent If $\Cr$ is an $n$-category then, for every $k$ in $\ens{0,\dots,n}$, the $k$-cells and $k$-spheres of $\Cr$ are in bijective correspondence with the $n$-functors from $\Er_k$ to $\Cr$ and from $\Sr_k$ to $\Cr$, respectively. When the context is clear, we use the same notation for a $k$-cell or $k$-sphere and its corresponding $n$-functor.

As a consequence, if $I$ is a set, the $I$-indexed families of $k$-cells (resp. $k$-spheres) of~$\Cr$ are in bijective correspondence with the $n$-functors from $I\cdot\Er_k$ (resp. $I\cdot\Sr_k$) to~$\Cr$. We recall that, for a set $X$ and an $n$-category $\Dr$, the copower $X\cdot\Dr$ is the coproduct $n$-category $\coprod_{x\in X}\Dr$, whose set of $k$-cells is the product $X\times\Dr_k$.

\subsubsection{Inclusion and collapsing \pdf{n}-functors} 

For every $n$, the \emph{inclusion $n$-functor} $J_n$ and the \emph{collapsing $n$-functor} $P_n$
$$
J_n \::\: \Sr_n \:\fl\: \Er_{n+1}
\qquad\text{and}\qquad
P_n \::\: \Sr_n\:\fl\: \Er_n
$$

\noindent are respectively defined as the canonical inclusion of $\Sr_n$ into $\Er_{n+1}$ and as the $n$-functor sending both~$e_n^-$ and~$e_n^+$ to~$e_n$, leaving the other cells unchanged. 

\subsection{Adjoining and collapsing cells}

\subsubsection{Definition}

Let $\Cr$ be an $n$-category, let $k$ be in $\ens{0,\dots,n-1}$, let $I$ be a set and let $\Gamma:I\cdot\Sr_k\fl\Cr$ be an $n$-functor. The \emph{adjoining of~$\Gamma$ to $\Cr$} and the \emph{collapsing of $\Gamma$ in $\Cr$} are the $n$-categories respectively denoted by $\Cr[\Gamma]$ and $\Cr/\Gamma$ and defined by the following pushouts in $\Cat{n}$:   
$$
\xymatrix{
I\cdot \Sr_k
	\ar[r]^-{\Gamma}
	\ar[d]_-{I\cdot J_k}
	\ar@{}[dr] |-{\copyright}
& \Cr
	\ar[d]
\\
{I\cdot \Er_{k+1}}
	\ar[r] 
& \Cr[\Gamma]
}
\qquad\qquad
\xymatrix{
I\cdot\Sr_k
	\ar[r]^-{\Gamma}
	\ar[d]_-{I\cdot P_k}
	\ar@{}[dr] |-{\copyright}
& \Cr
	\ar[d] 
\\
{I\cdot \Er_k}
	\ar[r]
& \Cr/\Gamma
}
$$

\noindent When $k=n$, one defines $\Cr[\Gamma]$ by seeing $\Cr$ as an $(n+1)$-category with degenerate $(n+1)$-cells only. 

The $n$-category $\Cr[\Gamma]$ has the same cells as $\Cr$
up to dimension $k$; its $(k+1)$-cells are all the formal
composites made of the $(k+1)$-cells of $\Cr$, plus one
extra $(k+1)$-cell from $\Gamma(i,e_k^-)$ to
$\Gamma(i,e_k^+)$ for every $i$ in $I$; above dimension
$k+1$, its cells are the ones of $\Cr$, plus the identities
of each extra cell. 

The $n$-category $\Cr/\Gamma$ has the same cells as $\Cr$ up to dimension $k-1$; its $k$-cells are the equivalence classes of $k$-cells of $\Cr$, for the congruence relation generated by $\Gamma(i,e_k^-)\sim\Gamma(i,e_k^+)$, for every $i$ in $I$; above dimension $k$, its cells are the formal composites of the ones of $\Cr$, but with sources and targets considered modulo the previous congruence.

\subsubsection{Extensions of \pdf{n}-functors}

Let $\Cr$ and $\Dr$ be $n$-categories and let $\Gamma:I\cdot\Sr_k\fl\Cr$ be an $n$-functor. Then, by universal property of $\Cr[\Gamma]$, one extends an $n$-functor $F:\Cr\fl\Dr$ to a unique $n$-functor $F:\Cr[\Gamma]\fl\Dr$ by fixing, for every $\gamma$ in $\Gamma$, a $(k+1)$-cell $F(\gamma)$ in $\Dr$ such that the following two equalities hold:
$$
s(F(\gamma)) \:=\: F(s(\gamma))
\qquad \text{and} \qquad
t(F(\gamma)) \:=\: F(t(\gamma)).
$$ 

\subsubsection{Occurrences}

Here we see the group $\Zb$ of integers as an $n$-category: it has one cell in each dimension up to $n-1$ and $\Zb$ as set of $n$-cells; all the compositions of $n$-cells are given by the addition.

Let $\Cr$ be an $n$-category and let $\Gamma:I\cdot\Sr_k\fl\Cr$ be an $n$-functor. We denote by $\norm{\cdot}_{\Gamma}$ the $n$-functor from~$\Cr[\Gamma]$ to $\Zb$ defined by:  
$$
\norm{f}_{\Gamma} \:=\: 
\begin{cases}
1 &\text{if } f\in\Gamma, \\
0 &\text{otherwise.}
\end{cases}
$$

\noindent For every cell $f$, one calls $\norm{f}_{\Gamma}$
the \emph{number of occurrences of cells of $\Gamma$ in
  $f$}. 

\subsubsection{The \pdf{n}-category presented by an \pdf{(n+1)}-category}

Let $\Cr$ be an $(n+1)$-category. If $f$ is an $(n+1)$-cell of $\Cr$, then $\dr f$ is an $n$-sphere of $\Cr$. Thus, the set $\Cr_{n+1}$ of $(n+1)$-cells of $\Cr$ yields an $(n+1)$-functor from $\Cr_{n+1}\cdot\Sr_n$ to the underlying $n$-category of $\Cr$: the \emph{$n$-category presented by $\Cr$} is the $n$-category denoted by $\cl{\Cr}$ one gets by collapsing the $(n+1)$-cells of $\Cr$ in its underlying $n$-category.

\subsection{Polygraphs and presentations of \pdf{n}-categories}
\label{Subsection:PolygraphsAndPresentations}

\emph{Polygraphs} (or \emph{computads}) are presentations by "generators" and "relations" of some higher-dimensional categories~\cite{Street76,Burroni93}, see also~\cite{Street87,Street95}. We define $n$-polygraphs by induction on the natural number $n$. 

The category $\Pol{0}$ of $0$-polygraphs and morphisms between them is the one of sets and maps.  A $0$-polygraph is \emph{finite} when it is finite as a set. A \emph{$0$-cell} of a $0$-polygraph is one of its elements. The \emph{free $0$-category functor} is the identity functor $\Pol{0}\fl\Cat{0}$. 

Now, let us fix a non-zero natural number $n$ and let us assume that we have defined the category $\Pol{n-1}$ of $(n-1)$-polygraphs and morphisms between them, finite $(n-1)$-polygraphs, $k$-cells of an $(n-1)$-polygraph and the free $(n-1)$-category functor $\Pol{n-1}\fl\Cat{n-1}$, sending an $(n-1)$-polygraph~$\Sigma$ to the $(n-1)$-category $\Sigma^*$.

\subsubsection{\pdf{n}-polygraphs}

An \emph{$n$-polygraph} is a pair $\Sigma=(\Sigma_{n-1},\Sigma_n)$ made of an $(n-1)$-polygraph $\Sigma_{n-1}$ and a family $\Sigma_n$ of $(n-1)$-spheres of the $(n-1)$-category $\Sigma_{n-1}^*$. 

An \emph{$n$-cell of $\Sigma$} is an element of $\Sigma_n$ and, if $k<n$, a \emph{$k$-cell of $\Sigma$} is a $k$-cell of the $(n-1)$-polygraph~$\Sigma_{n-1}$. The set of $k$-cells of $\Sigma$ is abusively denoted by $\Sigma_k$, thus identifying it to the $k$-polygraph underlying~$\Sigma$. An $n$-polygraph is \emph{finite} when it has a finite number of cells in every dimension. The \emph{size} of a $k$-cell $f$ in $\Sigma^\ast$, denoted by $\norm{f}$, is the natural number $\norm{f}_{\Sigma_k}$, giving the number of $k$-cells of $\Sigma$ that $f$ is made of. For $1$-cells, we also use $\abs{\cdot}$ instead of $\norm{\cdot}$.

The original paper~\cite{Burroni93} contains an equivalent description of $n$-polygraphs, where they are defined as diagrams
$$
\xymatrix{
\Sigma_0
	\ar@{=}[d] 
&& \Sigma_1
	\ar[dll] |-{s_0,t_0} 
	\ar@{ >->}[d]
&& (\cdots)
	\ar[dll] |-{s_1,t_1}
&& \Sigma_{n-1}
	\ar[dll] |-{s_{n-2},t_{n-2}}
	\ar@{ >->}[d]
&& \Sigma_n
	\ar[dll] |-{s_{n-1},t_{n-1}}
\\
\Sigma_0
&& \Sigma_1^\ast 
	\ar[ll] ^{\cl{s}_0,\cl{t}_0}
&&	(\cdots)
	\ar[ll]^{\cl{s}_1,\cl{t}_1}
&& \Sigma_{n-1}^{\ast}
	\ar[ll] ^-{\cl{s}_{n-2},\cl{t}_{n-2}}
}
$$

\noindent of sets and maps such that, for any $k$ in $\ens{0,\dots, n-1}$, the following two conditions hold:

\begin{itemize}

\item The diagram $
\xymatrix@1{
\Sigma_0^\ast 
& \Sigma_1^\ast
	\ar@<+0.2pc>[l]^-{\cl{s}_0}
	\ar@<-0.2pc>[l]_-{\cl{t}_0}
& (\cdots)
	\ar@<+0.2pc>[l]^-{\cl{s}_1}
	\ar@<-0.2pc>[l]_-{\cl{t}_1}
& \Sigma_{k}^\ast
	\ar@<+0.2pc>[l]^-{\cl{s}_{k-1}}
	\ar@<-0.2pc>[l]_-{\cl{t}_{k-1}}
}
$ is a $k$-category.

\item The diagram $
\xymatrix@1{
\Sigma_0^\ast 
& \Sigma_1^\ast
	\ar@<+0.2pc>[l]^-{\cl{s}_0}
	\ar@<-0.2pc>[l]_-{\cl{t}_0}
& (\cdots)
	\ar@<+0.2pc>[l]^-{\cl{s}_1}
	\ar@<-0.2pc>[l]_-{\cl{t}_1}
& \Sigma_{k}^\ast
	\ar@<+0.2pc>[l]^-{\cl{s}_{k-1}}
	\ar@<-0.2pc>[l]_-{\cl{t}_{k-1}}
& \Sigma_{k+1}
	\ar@<+0.2pc>[l]^-{s_k}
	\ar@<-0.2pc>[l]_-{t_k}	
}
$ is a $(k+1)$-graph. 
\end{itemize}

\subsubsection{Morphisms of \pdf{n}-polygraphs}

Let $\Sigma$ and $\Xi$ be two $n$-polygraphs. A \emph{morphism of $n$-polygraphs from $\Sigma$ to $\Xi$} is a pair $F=(F_{n-1},F_n)$ where $F_{n-1}$ is a morphism of $(n-1)$-polygraphs from $\Sigma_{n-1}$ to~$\Xi_{n-1}$ and where $F_n$ is a map from $\Sigma_n$ to $\Xi_n$ such that the following two diagrams commute:
$$
\xymatrix{
\Sigma_n 
	\ar [r] ^-{F_n} 
	\ar [d] _-{s_{n-1}}
	\ar@{} [dr] |-{\copyright}
& \Xi_n
	\ar [d] ^-{s_{n-1}}
\\
\Sigma_{n-1}^* 
	\ar [r] _-{F_{n-1}^*}
& \Xi_{n-1}^*
}
\qquad\qquad\qquad
\xymatrix{
\Sigma_n 
	\ar [r] ^-{F_n} 
	\ar [d] _-{t_{n-1}}
	\ar@{} [dr] |-{\copyright}
& \Xi_n
	\ar [d] ^-{t_{n-1}}
\\
\Sigma_{n-1}^* 
	\ar [r] _-{F_{n-1}^*}
& \Xi_{n-1}^*
}
$$

\noindent Alternatively, if $\Sigma_n:I\cdot\Sr_{n-1}\fl\Sigma_{n-1}^*$ and $\Xi_n:J\cdot\Sr_{n-1}\fl\Sigma_{n-1}^*$ are seen as $(n-1)$-functors, then $F_n$ is a map from $I$ to $J$ such that the following diagram commutes in $\Cat{n-1}$:
$$
\xymatrix{
I\cdot\Sr_{n-1}
	\ar[r] ^-{\Sigma_n}
	\ar[d] _-{F_n\cdot 1_{\Sr_{n-1}}}
	\ar@{}[dr] |-{\copyright}
& \Sigma_{n-1}^*
	\ar[d] ^-{F_{n-1}^*}
\\
J\cdot\Sr_{n-1}
	\ar[r] _-{\Xi_n}
& \Xi_{n-1}^*
}
$$

\noindent We denote by $\Pol{n}$ the category of polygraphs and morphisms between them.  

\subsubsection{The free \pdf{n}-category functor}

Let $\Sigma$ be an $n$-polygraph. The \emph{$n$-category freely generated by $\Sigma$} is the $n$-category $\Sigma^*$ defined as follows:
$$
\Sigma^* \:=\: \Sigma_{n-1}^*[\Sigma_n].
$$

\noindent This construction extends to an $n$-functor $(\cdot)^\ast : \Pol{n}\fl\Cat{n}$ called the \emph{free $n$-category functor}.  

\subsubsection{The \pdf{n}-category presented by an \pdf{(n+1)}-polygraph}

Let $\Sigma$ be a $(n+1)$-polygraph. The \emph{$n$-category presented by $\Sigma$} is the $n$-category denoted by $\overline{\Sigma}$ and defined as follows: 
$$
\cl{\Sigma} \: = \: \Sigma_n^*/\Sigma_{n+1}.
$$

\noindent Two $n$-polygraphs are \emph{Tietze-equivalent} when the $(n-1)$-categories they present are isomorphic. If $\Cr$ is an $n$-category, a \emph{presentation of $\Cr$} is an $(n+1)$-polygraph $\Sigma$ such that $\Cr$ is isomorphic to the $n$-category~$\cl{\Sigma}$ presented by $\Sigma$. One says that an $n$-category $\Cr$ is \emph{finitely generated} when it admits a presentation by an $(n+1)$-polygraph $\Sigma$ whose underlying $n$-polygraph $\Sigma_n$ is finite. One says that $\Cr$ is \emph{finitely presented} when it admits a finite presentation.

\subsubsection{Example: a presentation of the \pdf{2}-category of permutations} 

The $2$-category $\Pr erm$ of permutations has one $0$-cell, one $1$-cell for each natural number and, for each pair $(m,n)$ of natural number, its set of $2$-cells from $m$ to $n$ is the group $S_n$ of permutations if $m=n$ and the empty set otherwise. The $0$-composition of $1$-cells is the addition of natural numbers. The $0$-composition of two $2$-cells $\sigma\in S_m$ and $\tau\in S_n$ is the permutation $\sigma\star_0\tau$ defined by: 
$$
\sigma\star_0\tau(i) \:=\:
\begin{cases}
\sigma(i) &\text{if } 1\leq i\leq n, \\
\tau(i-n) &\text{otherwise.}
\end{cases}
$$

\noindent Finally the $1$-composition of $2$-cells is the composition of permutations. The $2$-category $\Pr erm$ is presented by the $3$-polygraph with one $0$-cell, one $1$-cell, one $2$-cell, pictured by~$\figeps{tau}$, and the following two $3$-cells:
$$
\scalefigeps{-2.8mm}{tautau} \tfl \scalefigeps{-2.8mm}{id2} 
\qquad\text{and}\qquad 
\scalefigeps{-4.3mm}{yb1} \tfl \scalefigeps{-4.3mm}{yb2}\:\:.
$$

\section{Contexts, modules and derivations of \pdf{n}-categories}
\label{Section:ContexteModuleDerivation}

\subsection{The category of contexts of an \pdf{n}-category}
\label{Subsection:Contexte}

Throughout this section, $n$ is a fixed natural number and $\Cr$ is a fixed $n$-category. 

\subsubsection{Contexts of an \pdf{n}-category}

A \emph{context of $\Cr$} is a pair $(x,C)$ made of an $(n-1)$-sphere $x$ of $\Cr$ and an $n$-cell~$C$ in $\Cr[x]$ such that $\norm{C}_x=1$. We often denote by $C[x]$, or simply by $C$, such a context.  

Let $x$ and $y$ be $(n-1)$-spheres of $\Cr$ and let $f$ be an $n$-cell in $\Cr[x]$ such that $\dr f=y$ holds. We denote by~$C[f]$ the image of a context $C[y]$ of $\Cr$ by the functor $\Cr[y] \fl \Cr[x]$ that extends $\id_{\Cr}$ with $y\mapsto f$.

\subsubsection{The category of contexts}

The \emph{category of contexts of $\Cr$} is the category
denoted by $\Ct{\Cr}$, whose objects are the $n$-cells of
$\Cr$ and whose morphisms from $f$ to $g$ are the contexts
$C[\dr f]$ of $\Cr$ such that $C[f]=g$ holds. If $C:f\fl g$
and $D:g\fl h$ are morphisms of $\Ct{\Cr}$ then $D\circ
C:f\fl h$ is $D[C]$. The identity context on an
$n$-cell $f$ of $\Cr$ is the context $\dr f$. When
$\Sigma$ is an $n$-polygraph, one writes $\Ct{\Sigma}$
instead of $\Ct{\Sigma^*}$. 

\subsubsection{Proposition}
\label{eqnDecompositionContexts}

\begin{em}
Every context of $\Cr$ has a decomposition
$$
f_n \star_{n-1} ( f_{n-1} \star_{n-2} \cdots (
f_1 \star_0 x\star_0 g_1 ) \cdots \star_{n-2} g_{n-1} )
\star_{n-1} g_n, 
$$

\noindent where $x$ is an $(n-1)$-sphere and, for every $k$ in $\ens{1,\dots,n}$, $f_k$ and $g_k$ are $n$-cells of $\Cr$. Moreover, one can choose these cells so that~$f_k$ and~$g_k$ are (the identities of) $k$-cells.  
\end{em}

\begin{proof}
The set of $n$-cells $f$ of $\Cr[x]$ such that $\norm{f}_x=1$ is a quotient of the following inductively defined set~$X$: the $n$-cell $x$ is in $X$; if $C$ is in $X$ and $f$ is an $n$-cell of $\Cr$ such that $t_i(f)=s_i(C)$ (resp. $t_i(C)=s_i(f)$) holds for some $i$, then $f\star_i C$ (resp. $C\star_i f$) is in $X$.

Using the associativity and exchange relations satisfied by the compositions of $\Cr$, one can order these successive compositions to reach the required shape, or to reach the same shape with $f_k$ and $g_k$ being identities of $k$-cells.
\end{proof}

\subsubsection{Whiskers}

A \emph{whisker of $\Cr$} is a context with a decomposition
$$
f_{n-1} \star_{n-2}  \cdots ( f_1 \star_0 x\star_0 g_1 ) \cdots \star_{n-2} g_{n-1} 
$$

\noindent such that, for every $k$ in $\ens{1,\dots,n-1}$, $f_k$ and $g_k$ are $k$-cells. We denote by $\whisk{\Cr}$ the subcategory of~$\Ct{\Cr}$ with the same objects and with whiskers as morphisms. When $\Sigma$ is an $n$-polygraph, we write $\whisk{\Sigma}$ instead of $\whisk{\Sigma^*}$.

\subsubsection{Proposition}
\label{Eqn:DecompoCellPoly}

\begin{em}
Let $\Sigma$ be an $n$-polygraph. Every $n$-cell $f$ in $\Sigma^{\ast}$ with size $k\geq 1$ has a decomposition 
$$
f \:=\: C_1[\gamma_1] \star_{n-1} \cdots \star_{n-1} C_k[\gamma_k].
$$

\noindent where $\gamma_1$, $\dots$, $\gamma_k$ are $n$-cells in $\Sigma$ and $C_1$, $\dots$, $C_k$ are whiskers of $\Sigma^*$.
\end{em}

\begin{proof}
We proceed by induction on the size of the $n$-cell $f$. 
If it has size $1$, then it contains exactly one $n$-cell $\gamma$ of $\Sigma$, possibly composed with other ones of lower dimension. Using the relations satisfied by compositions in an $n$-category, one can write $f$ as $C[\gamma]$, with~$C$ a context of $\Sigma^*$. Moreover, this context must be a whisker, since $f$ has size $1$.

Now, let us assume that we have proved that every $n$-cell with size at most $k$, for a fixed non-zero natural number $k$, admits a decomposition as in Proposition~\ref{Eqn:DecompoCellPoly}. Then let us consider an $n$-cell $f$ with size $k+1$. Since $\norm{f}\geq 2$ and by construction of $\Sigma^*=\Sigma^*_{n-1}[\Sigma_n]$, one gets that $f$ can be written $g\star_i h$, where $(g,h)$ is a pair of $i$-composable $n$-cells of $\Sigma^*$, for some $i$ in $\ens{0,\dots,n-1}$, with $\norm{g}$ and $\norm{h}$ at least $1$. One can assume that $i=n-1$ since, otherwise, one considers the following alternative decomposition of $f$, thanks to the exchange relation between $\star_i$ and $\star_{n-1}$:
$$
f \:=\: \left( g \star_i s(h) \right) \star_{n-1} \left( t(g) \star_i h \right) .
$$

\noindent Since $\norm{f}=\norm{g}+\norm{h}$, one must have $\norm{g}\leq k$ and $\norm{h}\leq k$. We use the induction hypothesis to decompose $g$ and $h$ as in~\ref{Eqn:DecompoCellPoly}, where $j$ denotes $\norm{g}$:
$$
g \:=\: C_1[\gamma_1] \star_{n-1} \cdots \star_{n-1} C_j [\gamma_j]
\qquad\text{and}\qquad
h \:=\: C_{j+1}[\gamma_{j+1}] \star_{n-1} \cdots \star_{n-1} C_k [\gamma_k].
$$

\noindent We compose the right members and use the associativity of $\star_{n-1}$ to conclude.
\end{proof}

\subsection{Contexts in low dimensions}

\subsubsection{Contexts of a \pdf{1}-category as factorizations}

From Proposition \ref{eqnDecompositionContexts}, we know that the contexts of a $1$-category~$\Cr$ have the following shape: 
$$
u\star_0 x \star_0 v,
$$

\noindent where $x$ is a $0$-sphere and $u$, $v$ are $1$-cells of $\Cr$. The morphisms in $\Ct{\Cr}$ from $w:p\fl q$ to $w':p'\fl q'$ are the pairs $(u:p'\fl p,v:q\fl q')$ of $1$-cells of $\Cr$ such that $u\star_0 w\star_0 v=w'$ holds in $\Cr$:
$$
\xymatrix{
p
	\ar[d] _-{w} 
	\ar@{} [dr] |-{\copyright}
& p'
	\ar[l] _-{u}
	\ar[d] ^-{w'} 
\\
q 
	\ar[r] _-{v}
& q'
}
$$

\noindent When $\Cr$ is freely generated by a $1$-polygraph, the $1$-cells $u$ and $v$ are uniquely defined by the context. Moreover, the contexts from $w$ to $w'$ are in bijective correspondence with the occurrences of the word~$w$ in the word $w'$.
The category $\Ct{\Cr}$ has been introduced by Quillen under the name \emph{category of factorizations of $\Cr$} \cite{Quillen72}. It has been used by Leech to introduce cohomological properties of congruences on monoids \cite{Leech85} and by Baues and Wirsching for the cohomology of small categories \cite{BauesWirsching85}.  

\subsubsection{Contexts of \pdf{2}-categories}

Let $\Cr$ be a $2$-category. From Proposition
\ref{eqnDecompositionContexts}, a context of $\Cr$
has the following shape: 
$$
h\star_1(g_1\star_0 x \star_0 g_2)\star_1 k
$$

\noindent where $x$ is a $1$-sphere and $g_1$, $g_2$, $h$, $k$ are
$2$-cells. Morphisms in $\Ct{\Cr}$ from a $2$-cell $f$ to
a $2$-cell $f'$ are the contexts $h\star_1(g_1\star_0
x \star_0 g_2)\star_1 k$ of $\Cr$ such that  
$$
h\star_1(g_1\star_0 f \star_0 g_2)\star_1 k \:=\: f'
$$ 

\noindent holds in $\Cr$. This last relation is graphically represented as follows:
$$ 
\xymatrix{
\bullet
	\ar @/^9ex/ [rrr]  _-{} ="sourceA"
	\ar @/_9ex/ [rrr]  ^-{} ="targetB"
	\ar @/^3ex/ [r]  _-{} ="source1" 
	\ar @/_3ex/ [r]  ^-{} ="target1" 
	\ar @2 "source1" ; "target1" ^- {g_1}
& \bullet 
	\ar @/^1.2pc/ [r]  _-{} ="source2" ^-{}="targetA"
	\ar @/_1.2pc/ [r]  ^-{} ="target2" _-{}="sourceB"
	\ar @2 "source2" ; "target2" ^- f 
	\ar @2 "sourceA" ; "targetA" ^- {h} 
	\ar @2 "sourceB" ; "targetB" ^- {k}
& \bullet
	\ar @/^3ex/ [r]  _-{} ="source3"
	\ar @/_3ex/ [r]  ^-{} ="target3"
	\ar @2 "source3" ; "target3" ^- {g_2}
& \bullet }
\qquad = \qquad
\xymatrix{
\bullet 
	\ar @/^3ex/ [r]  _-{} ="source1"
	\ar @/_3ex/ [r]  ^-{} ="target1"
	\ar @2 "source1" ; "target1" ^- {f'}
& \bullet
}
$$

\noindent However, the exchange relation between the two
compositions~$\star_0$ and~$\star_1$ implies that this
decomposition is not unique. Two decompositions 
$$
h\star_1(g_1\star_0 x \star_0 g_2)\star_1 k
\qquad\text{and}\qquad
h'\star_1(g_1'\star_0 x' \star_0 g_2')\star_1 k'
$$

\noindent represent the same context if and only if $x=x'$ and there exist $2$-cells $l_1$, $l_2$, $m_1$, $m_2$ such that the following four relations are defined and statisfied in $\Cr$: 
$$
\xymatrix@R=15mm{
&& \bullet
	\ar@/^4ex/ [rrr] _-{} ="src h"
	\ar[r] 
& \bullet
	\ar[r] ^-{} ="tgt h" _-{} ="src eq h"
& \bullet
	\ar[r]
& \bullet 
	\ar @2 "src h";"tgt h" ^-{h}
\\
&& \bullet
	\ar @/^4ex/ [rrr] _-{} ="src h p" ^-{} ="tgt eq h"
	\ar @/^1ex/ [r] _-{} ="src l1"
	\ar @/_3ex/ [r] ^-{} ="tgt l1"
& \bullet
	\ar[r] ^-{} ="tgt h p"
& \bullet
	\ar @/^1ex/ [r] _-{} ="src l2"
	\ar @/_3ex/ [r] ^-{} ="tgt l2"
& \bullet 
	\ar @2 "src h p" ; "tgt h p" ^-{h'}
	\ar @2 "src l1" ; "tgt l1" ^-{l_1}
	\ar @2 "src l2" ; "tgt l2" ^-{l_2}
	\ar @{} "src eq h" ; "tgt eq h" |-{=}
\\ 
\bullet
	\ar@/^6ex/ [r] _-{}="1"
	\ar@/^2ex/ [r] ^-{}="2" _-{}="3"
	\ar@/_2ex/ [r] ^-{}="4" _-{}="5"
	\ar@/_6ex/ [r] ^-{}="6"
	\ar@2 "1";"2" ^-{l_1}
	\ar@2 "3";"4" ^-{g_1}
	\ar@2 "5";"6" ^-{m_1}
& \bullet
	\ar @{} [r] |-{=}
& \bullet
	\ar @/^3ex/ [r] _-{} ="1" 
	\ar @/_3ex/ [r] ^-{} ="2" 
	\ar @2 "1" ; "2" ^- {g'_1}
& \bullet
& \bullet
	\ar @/^3ex/ [r] _-{} ="1" 
	\ar @/_3ex/ [r] ^-{} ="2" 
	\ar @2 "1" ; "2" ^- {g'_2}
& \bullet
	\ar @{} [r] |-{=}
& \bullet
	\ar@/^6ex/ [r] _-{}="1"
	\ar@/^2ex/ [r] ^-{}="2" _-{}="3"
	\ar@/_2ex/ [r] ^-{}="4" _-{}="5"
	\ar@/_6ex/ [r] ^-{}="6"
	\ar@2 "1";"2" ^-{l_2}
	\ar@2 "3";"4" ^-{g_2}
	\ar@2 "5";"6" ^-{m_2}
& \bullet
\\
&& \bullet
	\ar @/_4ex/ [rrr] ^-{} ="tgt k p" _-{} ="src eq k"
	\ar @/^3ex/ [r] _-{} ="src m1"
	\ar @/_1ex/ [r] ^-{} ="tgt m1"
& \bullet
	\ar [r] ^-{} ="src k p"
& \bullet
	\ar @/^3ex/ [r] _-{} ="src m2"
	\ar @/_1ex/ [r] ^-{} ="tgt m2"
& \bullet 
	\ar @2 "src k p";"tgt k p" ^-{k'}
	\ar @2 "src m1";"tgt m1"^-{m_1}
	\ar @2 "src m2";"tgt m2" ^-{m_2}
\\ 
&& \bullet
	\ar@/_4ex/ [rrr] ^-{} ="tgt k" 
	\ar[r] 
& \bullet
	\ar[r] ^-{} ="src k" ^-{} ="tgt eq k"
& \bullet
	\ar[r]
& \bullet 
	\ar @2 "src k";"tgt k" ^-{k}
	\ar @{} "src eq k"; "tgt eq k" |-{=}
}
$$

\noindent 

\subsection{Modules over \pdf{n}-categories}
\label{Subsection:Module}

\subsubsection{Definition}

Let $\Cr$ be an $n$-category. A \emph{$\Cr$-module} is a functor from the category of contexts $\Ct{\Cr}$ to the category $\Ab$ of abelian groups. Hence, a $\Cr$-module $M$ is specified by an abelian group $M(f)$, for every $n$-cell $f$ in $\Cr$, and a morphism $M(C):M(f)\fl M(g)$ of groups, for every context $C :f\fl g$ of~$\Cr$. When no confusion may occur, one writes $C[m]$ instead of $M(C)(m)$ and, when $C$ has shape $h\star_i x$ (resp. $x\star_i h$), one writes $h\star_i m$ (resp. $m\star_i h$) instead of $M(C)(m)$. 

\subsubsection{Proposition}
\label{Proposition:ModuleDecomposition}

\begin{em}
Let $\Cr$ be an $n$-category. A $\Cr$-module $M$ is entirely and uniquely defined by its values on the following contexts of $\Cr$:
$$
f\star_i x \qquad\text{and}\qquad x\star_i f
$$

\noindent for every $i$ in $\ens{0,\dots,n-1}$ and every non-degenerate $(i+1)$-cell $f$ in $\Cr$. 

Moreover, when $\Sigma$ is an $n$-polygraph, then a $\Sigma^*$-module $M$ is entirely and uniquely defined by its values on the following contexts of $\Sigma^*$:
$$
C[\phi]\star_i x \qquad\text{and}\qquad x\star_i C[\phi]
$$

\noindent for every $i$ in $\ens{0,\dots,n-1}$, every generating $(i+1)$-cell $\phi$ in $\Sigma_{i+1}$ and every whisker $C[\dr\phi]$ of $\Sigma_{i+1}^*$. 
\end{em}

\begin{proof}
Let $h$, $h'$ be two $n$-cells of $\Cr$ and let $C[x]:h\fl h'$ be a morphism of $\Ct{\Cr}$. We use Proposition~\ref{eqnDecompositionContexts} to decompose $C[x]$ as follows:
$$
C[x] \:=\: f_n \star_{n-1} \cdots \star_1\left( f_1 \star_0
  x\star_0 g_1 \right)\star_1 \cdots \star_{n-1}g_n, 
$$

\noindent in such a way that, for every $k$ in $\ens{1,\dots,n}$, $f_k$ and $g_k$ are $k$-cells. Thus, in the category $\Ct{\Cr}$, the context~$C[x]$ decomposes into
$$
C[x] \:=\: C_n[x_n] \circ \cdots \circ C_1[x_1],
$$

\noindent where $x_1=x$ and, for every $i$ in $\ens{1,\dots,n}$, one has $C_i[x_i]=f_i\star_{i-1} x_i\star_{i-1} g_i$ and $x_{i+1}=\dr C_i[x_i]$. Moreover, each $C_i[x_i]$ splits into:
$$
C_i[x_i] \:=\: \left(y_i\star_{i-1} g_i\right) \circ \left(f_i\star_{i-1} x_i\right),
$$

\noindent where $y_i = \dr(f_i\star_{i-1} x_i)$. Thus, since $M$ is a functor, it is entirely defined by its values on the contexts with shape $f\star_i x$ or $x\star_i f$, with $i$ in $\ens{0,\dots,n-1}$ and $f$ a non-degenerate $(i+1)$-cell (indeed, when $f$ is degenerate as a $i$-cell, one has $x\star_i f=x$ and $M(x)$ is always an identity). This proves the first part of the result.

Now, let us continue, assuming that $\Cr$ is freely generated by an $n$-polygraph $\Sigma$. Let us consider the $n$-context $f\star_i x$, where $f$ is an $(i+1)$-cell of size $k\geq 1$. We decompose it as in Proposition~\ref{Eqn:DecompoCellPoly}:
$$
f \:=\: C_1[\phi_1] \star_i \cdots \star_i C_k[\phi_k],
$$

\noindent where $\phi_1$, $\dots$, $\phi_k$ are generating
$(i+1)$-cells and $C_1$, $\dots$, $C_k$ are
$i$-contexts. Thus, a context $f\star_i x$ decomposes into
$\Ct{\Sigma}$ as follows: 
$$
f\star_i x \:=\: \left( C_1[\phi_1] \star_i x_1 \right) \circ \cdots \circ \left( C_k[\phi_k] \star_i x_k \right),
$$

\noindent where $x_k=x$ and $x_j=\dr(C_{j+1}[\phi_{j+1}] \star_i x_{j+1})$. Proceeding similarly with contexts of the shape $x\star_i f$, one gets the result.
\end{proof}

\subsubsection{Example: the trivial module}

Let $\Cr$ be an $n$-category. The \emph{trivial $\Cr$-module} sends each $n$-cell of $\Cr$ to $\Zb$ and each context of $\Cr$ to the identity of $\Zb$.

\subsubsection{Example of modules over $2$-categories}
\label{Example:modulesOver2categories}

Let $\Vr$ be a concrete category. We view it as a
$2$-category with one $0$-cell, objects as $1$-cells and
morphisms as $2$-cells. The $0$-composition in given by the
cartesian product and the $1$-composition by the composition
of morphisms. 

Let us fix an internal abelian group $G$ in $\Vr$, a $2$-category $\Cr$ and $2$-functors $X:\Cr\fl\Vr$ and $Y:\Cr^{\text{co}}\fl\Vr$, where $\Cr^{\text{co}}$ is $\Cr$ where one has exchanged the source and target of every $2$-cell. Then, using Proposition~\ref{Proposition:ModuleDecomposition}, the following assignments yield a $\Cr$-module $M_{X,Y,G}$: 
\begin{itemize}
\item Every $2$-cell $f: u \dfl v$ is sent to the abelian group of morphisms: 
$$
M_{X,Y,G}(f) \:=\: \Vr\big(\: X(u)\times Y(v), \: G \big).
$$

\item If $w$ and $w'$ are $1$-cells of $\Cr$ and $C=w\star_0 x \star_0 w'$ is a context from $f: u \dfl v$ to $w\star_0 f \star_0 w'$, then $M_{X,Y,G}(C)$ sends a morphism $a:X(u)\times Y(v)\fl G$ in $\Vr$ to:
$$
\begin{array}{r c l}
X(w)\times X(u) \times X(w') \times Y(w)\times Y(v) \times Y(w') &\:\longrightarrow\:& G \\
(x',x,x'',y',y,y'') &\:\longmapsto\:& a(x,y).
\end{array}
$$

\item If $g:u'\dfl u$ and $h:v\dfl v'$ are $2$-cells of $\Cr$ and $C=g\star_1 x \star_1 h$ is a context from $f: u\dfl v$ to $g\star_1 f \star_1 h$, then $M_{X,Y,G}(C)$ sends a morphism $a:X(u) \times Y(v)\fl G$ in $\Vr$ to $a\circ (X\times Y)$, that is:
$$
\begin{array}{r c l}
X(u')\times Y(v') &\:\longrightarrow\:& G \\
(x,y) &\:\longmapsto\:& a \left( \; X(g)(x), \: Y(h)(y) \; \right).
\end{array}
$$
\end{itemize}

\noindent When $X$ or $Y$ is trivial, \ie, sends all the cells of $\Cr$ to the terminal object of $\Vr$, one denotes the corresponding $\Cr$-module by $M_{\ast,Y,G}$ or $M_{X,\ast,G}$. In particular, $M_{\ast,\ast,\Zb}$ is the trivial $\Cr$-module.

By construction, a $\Cr$-module $M_{X,Y,G}$ is uniquely and entirely defined by the values $X(u)$ and $Y(u)$, for every $1$-cell $u$, and by the morphisms $X(f)$ and $Y(f)$ for every $2$-cell $f$. As a consequence, when $\Cr$ is freely generated by a $2$-polygraph $\Sigma$, the $\Cr$-module $M_{X,Y,G}$ is uniquely and entirely determined by: 
\begin{itemize}
\item The objects $X(a)$ and $Y(a)$ of $\Vr$, for every generating $1$-cell $a$ in $\Sigma_1$. 
\item The morphisms $X(\gamma):X(u)\fl X(v)$ and $Y(\gamma):Y(v)\fl Y(u)$ of $\Vr$, for every generating $2$-cell $\phi:u\dfl v$ in~$\Sigma_2$. 
\end{itemize}

\noindent In the sequel, we consider this kind of $\Cr$-module with $\Vr$ being the category $\Set$ of sets and maps or the category $\Ord$ of partially ordered sets and monotone maps. For this last situation, we recall that an internal abelian group in $\Ord$ is a partially ordered set equipped with a structure of abelian group whose addition is monotone in both arguments.

\subsection{Derivations of \pdf{n}-categories}
\label{Subsection:Derivation}

\subsubsection{Definition}

Let $\Cr$ be an $n$-category and let $M$ be a $\Cr$-module. A \emph{derivation of $\Cr$ into $M$} is a map sending every $n$-cell $f$ of $\Cr$ to an element $d(f)$ of $M(f)$ such that the following relation holds, for every $i$-composable pair $(f,g)$ of $n$-cells of $\Cr$: 
$$
d(f\star_i g) \:=\: f\star_i d(g) + d(f) \star_i g.
$$

\noindent Given a derivation $d$ on $\Cr$, we define its values on contexts by
$$
d(C)= \sum_{i=-n}^n\;f_n \star_{n-1}(f_{n-1} \star_{n-2} \cdots 
(d(f_i)\star_{i-1}\cdots ( f_1 \star_0 x\star_0 f_{-1} )
\cdots \star_{n-1} f_{-n},
$$

\noindent for any context $C[x]=f_n \star_{n-1} \cdots ( f_1
\star_0 x\star_0 f_{-1} ) \cdots \star_{n-1} f_{-n}$ of
$\Cr$. This gives a mapping $d(C)$ taking an $n$-cell $f$ of
$\Cr$ with boundary $x$ to an element $d(C)[f]$ of the
abelian group $M(C[f])$. In this way a derivation from $\Cr$
into $M$ satisfies: 
$$
d(C[f]) = d(C)[f] + C[d(f)].
$$


\subsubsection{Proposition}

\begin{em}
Let $\Cr$ be an $n$-category, let $M$ be a $\Cr$-module and
let $d$ be a derivation of $\Cr$ into~$M$. Then, for every
degenerate $n$-cell $f$ of $\Cr$, we have
$d(f)=0$. Moreover, when $\Cr$ is the $n$-category freely
generated by an $n$-polygraph $\Sigma$, then $d$ is entirely
and uniquely determined by its values on the generating
cells of $\Sigma$. 
\end{em}

\begin{proof}
Let $f$ be a degenerate $n$-cell of $\Cr$. We have:
$$
d(f) \:=\: d(f\star_{n-1} f) \:=\: f\star_{n-1} d(f) + d(f) \star_{n-1} f \:=\: 2\cdot d(f).
$$

\noindent Since $d(f)$ is an element of the abelian group $M(f)$, then we have $d(f)=0$. 

As a consequence of its definition, a derivation is compatible with the associativity, unit and exchange relations. This implies that the values of $d$ on an $n$-cell $f$ of $\Sigma^*$ can be uniquely computed from its values on the generating $n$-cells $f$ is made of.
\end{proof}

\subsubsection{Example: occurrences}

If $\Cr$ is an $n$-category and $\Gamma:I\cdot\Sr_{n-1}\fl\Cr$ is an $n$-functor, we have defined the $n$-functor $\norm{\cdot}_\Gamma$ counting the number of occurrences of $n$-cells of $\Gamma$ in an $n$-cell of $\Cr[\Gamma]$. This construction is a derivation of $\Cr$ into the trivial $\Cr$-module, sending each $n$-cell of $\Cr$ to $0$ and each $n$-cell of $\Gamma$ to $1$.

\subsubsection{Example: derivations of free \pdf{2}-categories}

Let us consider a $2$-polygraph $\Sigma$, a concrete category~$\Vr$ and a module of the shape $M_{X,Y,G}$, as defined in \ref{Example:modulesOver2categories}. Then, by construction of $\Sigma^*$, a derivation~$d$ of $\Sigma^*$ into $M_{X,Y,G}$ is entirely and uniquely determined by a family $(d\phi)_{\phi\in\Sigma_2}$ made of a morphism
$$
d\phi \::\: X(u)\times Y(v) \:\fl\: G
$$ 

\noindent of $\Vr$ for each $2$-cell $\phi:u\dfl v$ of $\Sigma$. 

\section{Higher-dimensional categories with finite derivation type} 
\label{Section:TDF}
  
\subsection{Track \pdf{n}-categories}
\label{Subsection:TracknCat}

\subsubsection{Definitions}

In an $n$-category $\Cr$, a $k$-cell $f$ is
\emph{invertible} when there exists a $k$-cell $g$ from
$t(f)$ to~$s(f)$ in $\Cr$ such that both
$f\star_{k-1}g=s(f)$ and $g\star_{k-1}f=t(f)$ hold. In that
case, $g$ is unique and denoted by $f^{-1}$. The following
relations are satisfied: 
$$
(1_x)^{-1} \:=\: 1_x 
\qquad \text{and} \qquad
(f\star_i g)^{-1} \:=\: 
	\begin{cases}
		f^{-1} \star_i g^{-1} &\text{if } i<k-1 \\
		g^{-1} \star_{k-1} f^{-1} &\text{otherwise.}
	\end{cases}
$$

\noindent Moreover, if $F:\Cr\fl\Dr$ is an $n$-functor, one has:
$$
F(f^{-1}) \:=\: F(f)^{-1}.
$$

\noindent A \emph{track $n$-category} is an $n$-category whose $n$-cells are invertible, \ie, an $(n-1)$-category enriched in groupoid. We denote by $\Tck{n}$ the category of track $n$-categories and $n$-functors between them. 

\subsubsection{Example}

Let $\Cr$ be an $n$-category. Given two $n$-cells $f$ from $u$ to $v$ and $g$ from $v$ to $u$ in $\Cr$, we denote by $I_{f,g}$ the following $n$-sphere of $\Cr$:
$$
I_{f,g} = \left( f\star_{n-1}g ,\: 1_u \right).
$$

\noindent If $\gamma=(f,g)$ is an $n$-sphere of $\Cr$, we denote by $\gamma^{-1}$ the $n$-sphere $(g,f)$ of $\Cr$. Then we define the track $(n+1)$-category $\Cr(\gamma)$ by:
$$
\Cr(\gamma) = \Cr \left[ \gamma, \; \gamma^{-1} \right] \: \big/ \: \ens{I_{\gamma,\gamma^{-1}}, \; I_{\gamma^{-1},\gamma}}.
$$ 

\noindent This construction is extended to a set $\Gamma$ of $n$-spheres, yielding a track $(n+1)$-category $\Cr(\Gamma)$. 

\subsubsection{The free track \pdf{n}-category functor}

Given an $n$-polygraph $\Sigma$, the \emph{track $n$-category freely generated by $\Sigma$} is the $n$-category denoted by~$\Sigma^\top$ and defined by:
$$
\Sigma^\top \:=\: \Sigma_{n-1}^{\ast}(\Sigma_n),
$$

\noindent This construction extends to a functor $(\cdot)^\top : \Pol{n} \fl \Tck{n}$ called the \emph{free track $n$-category functor}.

\subsection{Homotopy bases}
\label{SubSubsection:HomotopyBase}

\subsubsection{Homotopy relation}

Let $\Cr$ be an $n$-category. A \emph{homotopy relation} on $\Cr$ is a track $(n+1)$-category~$\Tr$ with $\Cr$ as underlying $n$-category. Given an $n$-sphere $(f,g)$ in $\Cr$, one denotes by $f\approx_{\Tr}g$ the fact that there exists an $(n+1)$-cell from $f$ to $g$ in $\Tr$. If $\Gamma$ is a set of $n$-spheres of $\Cr$, one simply writes~$\approx_\Gamma$ instead of $\approx_{\Cr(\Gamma)}$ and calls it the \emph{homotopy relation on $\Cr$ generated by $\Gamma$}. 

One has $f\approx_{\Tr} g$ if and only if $\pi(f)=\pi(g)$ holds, where $\pi$ is the canonical projection from $\Tr$ to the $n$-category $\cl{\Tr}$ presented by $\Tr$, \ie, $\Cr/\Tr_{n+1}$. As a consequence, the relation $\approx_{\Tr}$ is a congruence relation on the parallel $n$-cells of $\Cr$, \ie, it is an equivalence relation compatible with every composition of $\Cr$. 

\subsubsection{Homotopy basis}

A set $\Gamma$ of $n$-spheres of $\Cr$ is a \emph{homotopy basis of $\Cr$} when, for every $n$-sphere~$(f,g)$ of $\Cr$, one has $f\approx_{\Gamma} g$. In other words, $\Gamma$ is a homotopy basis if and only if, for every $n$-sphere $\gamma$ of $\Cr$, there exists an $(n+1)$-cell $\cl{\gamma}$ such that $\dr\cl{\gamma}=\gamma$ holds, \ie, such that the following diagram commutes in $\Cat{n+1}$:
$$
\xymatrix{
\Sr_n 
	\ar[r] ^-{\gamma}
	\ar[d] _-{J_n}
	\ar@{} [dr] |-{\copyright}
& \Cr 
	\ar[d]
\\
\Er_{n+1}
	\ar[r] _-{\cl{\gamma}}
& \Cr(\Gamma)
}
$$

\subsubsection{Proposition}
\label{PropExtractFiniteHomotopyBase}

\begin{em}
Let $\Cr$ be an $n$-category and let $\Gamma$ be a homotopy
basis of $\Cr$. If $\Cr$ admits a finite homotopy basis,
then there exists a finite subset of $\Gamma$ that is a
homotopy basis of $\Cr$. 
\end{em}

\begin{proof}
Let $\Gamma'$ be a finite homotopy basis of $\Cr$. Let $\gamma$ be an $n$-sphere of $\Cr$ in $\Gamma'$. Since $\Gamma$ is a homotopy basis of $\Cr$, there exists an $(n+1)$-cell $\phi_{\gamma}$ in $\Cr(\Gamma)$ with boundary $\gamma$. This defines an $(n+1)$-functor $F$ from $\Cr(\Gamma')$ to $\Cr(\Gamma)$ which is the identity on cells of $\Cr$ and which sends each $\gamma$ in $\Gamma'$ to $\phi_{\gamma}$. For each $\phi_{\gamma}$, we fix a representative in $\Cr[\Gamma,\Gamma^{-1}]$ and denote by $\ens{\phi_{\gamma}}_{\Gamma}$ the set of cells of $\Gamma$ occurring in this representative. Let us denote by $\Gamma_0$ the following subset of $\Gamma$
$$
\Gamma_0 \:=\: \bigcup_{\gamma\in\Gamma'} \ens{\phi_{\gamma}}_{\Gamma} \:,
$$

\noindent consisting of all the cells of $\Gamma$ contained in the cells $\phi_{\gamma}$. The subset $\Gamma_0$ is finite since $\Gamma'$ and each $\ens{\phi_{\gamma}}_{\Gamma}$ are. Now let us see that it is an homotopy basis of $\Cr$. Let us fix an $n$-sphere $(f,g)$ of $\Cr$. By hypothesis, there exists an $(n+1)$-cell $A$ in $\Cr(\Gamma')$ with boundary $(f,g)$. By application of $F$, one gets an $(n+1)$-cell $F(A)$ in $\Cr(\Gamma)$ with boundary $(f,g)$. Moreover, the $(n+1)$-cell $F(A)$ is a composite of cells of the shape~$\phi_{\gamma}$: hence, it lives in $\Cr(\Gamma_0)$. As a consequence, one gets $f\approx_{\Gamma_0}g$, which concludes the proof.
\end{proof}

\subsection{Polygraphs with finite derivation type}

\subsubsection{Definitions}

One says that an $n$-polygraph $\Sigma$ has \emph{finite
  derivation type} when it is finite and when the track
$n$-category $\Sigma^\top$ it generates admits a finite
homotopy basis.
An $n$-category has \emph{finite derivation type} when it
admits a presentation by an $(n+1)$-polygraph with finite
derivation type.

\subsubsection{Lemma}
\label{LemmaTech1:FDTTietzeEquivalent}

\begin{em}
Let $\Sigma$ and $\Sigma'$ be $n$-polygraphs. We denote by $\pi: \Sigma^\ast_{n-1} \fl \cl{\Sigma}$ and by $\pi': {\Sigma'}^\ast_{n-1}\fl\cl{\Sigma}'$ the canonical $(n-1)$-functors. Then every $(n-1)$-functor $F$ from $\cl{\Sigma}$ to $\cl{\Sigma}'$ can be lifted to an $n$-functor $\widetilde{F} : \Sigma^\top \fl {\Sigma'}^\top$ such that the following diagram commutes in $\Cat{n-1}$: 
$$
\xymatrix{
\Sigma^\ast_{n-1} 
	\ar[r] ^-{\pi}
	\ar[d] _-{\widetilde{F}}
	\ar@{}[dr] |- {\copyright}
& \cl{\Sigma}
	\ar[d] ^-{F}
\\
{\Sigma'}^\ast_{n-1}
	\ar[r] _-{\pi'} 
& \cl{\Sigma}'
}
$$
\end{em}

\begin{proof} 
For every $k$-cell $u$ in $\Sigma^\ast$, with $k$ in $\ens{0,\dots,n-2}$, we take $\widetilde{F}(u)=F(u)$. Since $\pi$ and $\pi'$ are identities on cells up to dimension $n-2$, we have the relation $F\circ\pi(u)=\pi'\circ\widetilde{F}(u)$. 

Now, let us consider an $(n-1)$-cell $u$ in $\Sigma$. One arbitrarily chooses an $(n-1)$-cell of ${\Sigma'}^\ast$, hence of ${\Sigma'}^\top$, that is sent on $F\circ\pi(u)$ by $\pi'$, and one fixes $\widetilde{F}(u)$ to that $(n-1)$-cell. One extends $\widetilde{F}$ to every $(n-1)$-cell of $\Sigma^\ast$ thanks to the universal property of $\Sigma^\ast$. 

Then, let $f$ be an $n$-cell from $u$ to $v$ in $\Sigma$. Then $\pi(u)=\pi(v)$ holds by definition of $\pi$. Applying~$F$ on both members and using the property satisfied by $\widetilde{F}$, one gets $\pi'\circ\widetilde{F}(u)=\pi'\circ\widetilde{F}(v)$. By definition of~$\pi'$ and of ${\Sigma'}^\top$, this means that there exists an $n$-cell from $\widetilde{F}(u)$ to $\widetilde{F}(v)$ in ${\Sigma'}^\top$. One takes one such $n$-cell for $\widetilde{F}(f)$. Finally, one extends $\widetilde{F}$ to every $n$-cell of $\Sigma^\top$. 
\end{proof} 

\subsubsection{Lemma}
\label{lemmaImageofGamma}

\begin{em}
Let $\Sigma$ and $\Sigma'$ be $n$-polygraphs and let $F:\Sigma^\top \fl {\Sigma'}^\top$ be an $n$-functor. Given a set $\Gamma$ of $n$-spheres of $\Sigma^\top$, we define $F(\Gamma)$ as the following set of $n$-spheres of ${\Sigma'}^{\top}$:
$$
F(\Gamma) \:=\: \ens{ \: (F(g),F(g')) \:\big|\: (g,g')\in \Gamma \:} \:.
$$

\noindent Then, for every $n$-sphere $(f,f')$ of $\Sigma^\top$ such that $f\approx_\Gamma f'$ holds, we have $F(f) \approx_{F(\Gamma)} F(f')$. 
\end{em}

\begin{proof}
We use the functoriality of $F$.
\end{proof}

\subsubsection{Proposition}
\label{TDFTietzeInvariant}

\begin{em}
Let $\Sigma$ and $\Sigma'$ be Tietze-equivalent finite $n$-polygraphs. Then $\Sigma$ has finite derivation type if and only if $\Sigma'$ has.
\end{em}

\begin{proof}
Let us assume that $\Sigma$ and $\Sigma'$ are $n$-polygraphs which present the same $(n-1)$-category, say~$\Cr$.  Let us assume that $\Sigma$ has finite derivation type, so that we can fix a finite homotopy basis $\Gamma$ of $\Sigma^\top$. Using Lemma~\ref{LemmaTech1:FDTTietzeEquivalent} twice on the $(n-1)$-functor $\id_{\Cr}$, we get two $n$-functors $F:\Sigma^\top\fl {\Sigma'}^\top$ and $G: {\Sigma'}^\top \fl \Sigma^\top$ such that the following diagrams commute in $\Cat{n-1}$: 
$$
\xymatrix{
\Sigma^\ast_{n-1}
	\ar[r]^-{\pi} 
	\ar[d]_-F
	\ar@{}[dr] |-{\copyright}
& \Cr 
	\ar[d]^{\id_{\Cr}} 
\\
{\Sigma'}^\ast_{n-1}
	\ar[r]_-{\pi'} 
& \Cr
}
\qquad\qquad
\xymatrix{
\Sigma^\ast_{n-1}
	\ar[r]^-{\pi}  
	\ar@{}[dr] |-{\copyright}
& \Cr \\
{\Sigma'}^\ast_{n-1} 
	\ar[r]_-{\pi'} 
	\ar[u]^-G
& \Cr \ar[u]_{\id_{\Cr}}
}
$$

\noindent In particular, both $\pi$ and $\pi'$ are the identity on $k$-cells, for every $k<n-1$, hence so are $F$ and $G$. 

Let us consider an $(n-1)$-cell $a$ in $\Sigma'$. Then $\pi'\circ FG(a)=\pi\circ G(a)=\pi'(a)$. Thus, there exists an $n$-cell denoted by $f_a$ from $a$ to $FG(a)$ in ${\Sigma'}^\top$. From these cells, we define $f_u$ for every $(n-1)$-cell~$u$ in~${\Sigma'}^\ast$, hence of ${\Sigma'}^\top$, using the following relations: 
\begin{itemize}
\item for every degenerate $(n-1)$-cell $u$, $f_u$ is defined as $u$,
\item for every $i$-composable pair $(u,v)$ of $(n-1)$-cells, $f_{u\star_i v}$ is defined as $f_u\star_i f_v$. 
\end{itemize}

\noindent We have that, for every $(n-1)$-cell $u$, the $n$-cell $f_u$ goes from $u$ to $FG(u)$: to check this, we argue that~$FG$ is an $n$-functor which is the identity on degenerate $(n-1)$-cells. 

Now, let us consider an $n$-cell $g$ from $u$ to $v$ in ${\Sigma'}^\top$. We denote by $f_g$ the following $n$-cell from~$u$ to~$u$ in ${\Sigma'}^\top$, with a cellular representation giving the intuition for the case $n=2$: 
$$
f_g=
g\star_{n-1}f_v\star_{n-1}
FG(g)^{-1}\star_{n-1}f_u^{-1}
\qquad\qquad
\xymatrix{
\bullet 
\ar @/^4pc/ [rrr] ^{FG(u)} _{} ="tGammaU" ="tFGalpha" 
\ar @/^1.5pc/ [rrr]  ^u _(.25){} ="sGammaU" ^(.25){} ="salpha"
\ar @/_1.5pc/ [rrr] _(.75){} ="sFGalpha" ="tGammaV" ^{FG(v)} 
\ar @/_4pc/ [rrr] _{v} ^{} ="sGammaV" ="talpha"
\ar @{=>}@/^0.7pc/ "sGammaU" ; "tGammaU" _{f_u }
\ar @{=>}@/_1.1pc/ "salpha" ; "talpha" ^(.75){g}
\ar @{=>}@/_0.7pc/ "sGammaV" ; "tGammaV" ^{f_v}
\ar @{=>}@/_1.3pc/ "sFGalpha" ; "tFGalpha" _{\;\;FG(g)^{-1}}
&&&\bullet
}
$$

\noindent Let us prove that, for any composable pair $(g,h)$ of $n$-cells in ${\Sigma'}^\top$, we have: 
$$
f_{g\star_{n-1} h} \:=\: 
g \star_{n-1} f_h \star_{n-1} g^{-1}\star_{n-1} f_g.
$$

\noindent For that, we assume that $g$ has source $u$ and target $v$, while $h$ has source $v$ and target $w$. Then we compute: 
\begin{align*}
&\:\: g 
	\star_{n-1} f_h 
	\star_{n-1} g^{-1}
	\star_{n-1} f_{g} 
\\ 
=&\:\: g 
	\star_{n-1} \left( 
		h 
		\star_{n-1} f_w
		\star_{n-1} FG(h)^{-1} 
		\star_{n-1} f_v^{-1}
	\right)
\\
&\:\:	\star_{n-1} g^{-1} 
	\star_{n-1} \left(
		g 
		\star_{n-1} f_v 
		\star_{n-1} FG(g)^{-1} 
		\star_{n-1} f_u^{-1}
	\right) 
\\ 
=&\:\: g 
	\star_{n-1} h 
	\star_{n-1} f_w
	\star_{n-1} FG(h)^{-1} 
	\star_{n-1} FG(g)^{-1}
	\star_{n-1} f_u^{-1}
\\
=&\:\: (g \star_{n-1} h)
	\star_{n-1} f_w
	\star_{n-1} FG(g\star_{n-1} h)^{-1} 
	\star_{n-1} f_u^{-1}
\\
=&\:\: f_{g\star_{n-1}h}.
\end{align*}

\noindent Now, let us consider an $n$-cell $g$ and a whisker $C[x]$ in $\Sigma^\top$ such that $x=\dr(g_{n-1}^-)$. We note that, by definition of $f_g$, it has the same $(n-1)$-source and $(n-1)$-target as $g$, so that $C[f_g]$ is defined. Let us prove that the following relation holds: 
$$
f_{C[g]}= C[f_g].
$$

\noindent From the decomposition of contexts, it is sufficient to prove that the following relation holds 
$$
f_{u\star_i g\star_i v} \:=\: u\star_i f_g \star_i v
$$

\noindent for every $n$-cell $g$, every possible $k$-cells $u$ and $v$, with $k<n-1$, and every $i<k$ such that $u\star_i g \star_i v$ is defined. Let us assume that $g$ has source $w$ and target $w'$ and compute, from the left-hand side of this relation: 
\begin{align*}
f_{u\star_i g \star_i v}
	\:&=\: (u\star_i g \star_i v) 
	\star_{n-1} f_{u\star_i w'\star_i v} 
	\star_{n-1} FG(u\star_i g\star_i v)^{-1}
	\star_{n-1} f_{u\star_i w\star_i v}^{-1} 
\\ 
	\:&=\: (u\star_i g \star_i v) 
	\star_{n-1} (u\star_i f_{w'}\star_i v) 
	\star_{n-1} (u\star_i FG(g)^{-1} \star_i v)
	\star_{n-1} (u\star_i f_{w}^{-1}\star_i v)
\\
	\:&=\: u\star_i f_g \star_i v.
\end{align*}

\noindent Now, we denote by $\Gamma'$ the set of $n$-spheres $(f_g,1_{s(g)})$, for every $n$-cell $g$ in $\Sigma'$. Then, it follows from the previous relations that, for every $n$-cell $g$ in ${\Sigma'}^\top$, one has: 
$$
f_g \approx_{\Gamma} 1_{s(g)}.
$$

\noindent Let us consider an $n$-sphere $(g,g')$ in ${\Sigma'}^\top$. Then $(G(g),G(g'))$ is an $n$-sphere in $\Sigma^\top$. Since $\Gamma$ is a homotopy basis for $\Sigma^\top$, we have $G(g)\approx_{\Gamma}G(g')$, so that, by Lemma \ref{lemmaImageofGamma}, one gets $FG(g)~\approx_{F(\Gamma)}~FG(g')$.  

Finally, let us denote $\Gamma''$ the set of $n$-spheres of ${\Sigma'}^\top$ defined by $\Gamma''=\Gamma'\cup F(\Gamma)$ and let us prove that~$\Gamma''$ is a finite homotopy basis of ${\Sigma'}^\top$. Since both $\Sigma'_n$ and $\Gamma$ are finite, so is $\Gamma''$. Let us consider an $n$-sphere $(g,g')$ in ${\Sigma'}^\top$, with source $w$ and target $w'$, and let us prove that $g\approx_{\Gamma''} g'$ holds. We start by using the definition of $f_g$ to get: 
$$
g \:=\: f_g \star_{n-1} f_{w'}^{-1} \star_{n-1} FG(g) \star_{n-1} f_w^{-1}.
$$

\noindent Using the definition of $f_{g'}$, one gets a similar formula for $g'$. We have seen that $f_g\approx_{\Gamma'} w$, $f_{g'}\approx_{\Gamma'}w$ and $FG(g)\approx_{\Gamma''} FG(g')$ hold. Thus one gets $g\approx_{\Gamma''} g'$. 
\end{proof}

\subsubsection{Remark}

Proposition \ref{TDFTietzeInvariant} shows that the property
of having finite derivation type is invariant by
Tietze-equivalence for finite polygraphs. We will illustrate
in Example \ref{Example:infinitePolTDF} that this is not the
case for infinite ones.   

\section{Critical branchings and finite derivation type} 
\label{Section:RewritingTDF}

\subsection{Rewriting properties of polygraphs}
\label{Subsection:Rewriting}

We fix an $(n+1)$-polygraph $\Sigma$ and an $n$-cell $f$ in 
$\Sigma^\ast$. 

\subsubsection{Reductions and normal forms}

One says that $f$ \emph{reduces} into some $n$-cell $g$ when there exists a non-degenerate  $(n+1)$-cell $A$ from $f$ to $g$ in $\Sigma^*$. A \emph{reduction sequence} is a family $(f_k)_k$ of $n$-cells such that each $f_k$ reduces into $f_{k+1}$. One says that $f$ is \emph{a normal form (for $\Sigma_{n+1}$)} when every $(n+1)$-cell with source $f$ is degenerate, \ie, it does not reduce into any $n$-cell. A \emph{normal form for $f$} is a normal form~$g$ such that~$f$ reduces into $g$. The polygraph $\Sigma$ is \emph{normalizing at $f$} when $f$ admits a normal form. It is \emph{normalizing} when it is at every $n$-cell of $\Sigma^\ast$. 

\subsubsection{Termination}

One says that $\Sigma$ \emph{terminates at $f$} when there
exists no infinite reduction sequence starting at $f$. One
says that $\Sigma$ \emph{terminates} when it does at every
$n$-cell of $\Sigma^\ast$.  
If $\Sigma$ terminates at $f$, then it is normalizing at
$f$, \ie, every $n$-cell has at least one normal
form. Moreover, in case of termination, one can prove
properties using \emph{Noetherian induction}. For that, one
proves the property on normal forms; then one fixes an
$n$-cell $f$, one assumes that the result holds for every
$g$ such that $f$ reduces into $g$ and one proves that,
under those hypotheses, the $n$-cell $f$ satisfies the
property.  

\subsubsection{Confluence}

A \emph{branching of $\Sigma$} is a pair $(A,B)$ of $(n+1)$-cells of $\Sigma^\ast$ with same source; this $n$-cell is called the \emph{source} of the branching $(A,B)$. A branching $(A,B)$ is \emph{local} when $\norm{A}=\norm{B}=1$. A \emph{confluence of $\Sigma$} is a pair $(A,B)$ of $(n+1)$-cells of $\Sigma^\ast$ with same target. A branching $(A,B)$ is \emph{confluent} when there exists a confluence $(A',B')$ such that both $t_n(A)=s_n(A')$ and $t_n(B)=s_n(B')$ hold, as in the following diagram:
$$
\xymatrix{
& \ar[dl] _A \ar[dr]^B &\\
\ar[dr]_{A'}&&\ar[dl]^{B'} \\
&&}
$$

\noindent Such a pair $(A',B')$ is called a \emph{confluence for $(A,B)$}. Branchings and confluences are only considered up to symmetry, so that $(A,B)$ and $(B,A)$ are considered equal. The polygraph $\Sigma$ is \emph{(locally) confluent at $f$} when every (local) branching with source $f$ is confluent. It is \emph{(locally) confluent} when it is at every $n$-cell.

If $\Sigma$ is confluent then every $n$-cell of $\Sigma^*$ has at most one normal form. Thus, normalization and confluence imply that the $n$-cell $f$ has exactly one normal form, written $\widehat{f}$. In a terminating polygraph, local confluence and confluence are equivalent: this was proved in the case of word rewriting systems (a subcase of $2$-polygraphs) by Newman~\cite{Newman42} and, since then, the result is called Newman's lemma.

\subsubsection{Convergence}

The polygraph $\Sigma$ is \emph{convergent at $f$} when it terminates and it is confluent at $f$. It is \emph{convergent} when it is at every $n$-cell.
If $\Sigma$ is convergent at $f$, then $f$ has exactly one normal form. Thanks to Newman's lemma, one gets convergence from termination and local confluence. If $\Sigma$ is convergent, we have $f\approx_{\Sigma_{n+1}}g$ if and only if the equality $\widehat{f}=\widehat{g}$ holds. As a consequence, a finite and convergent $(n+1)$-polygraph provides a decision procedure to the equivalence of $n$-cells in the $n$-category it presents. 

\subsubsection{Critical branchings in polygraphs}
\label{SubSubSection:DefinitionCriticalBranching}

Given a branching $b=(A,B)$ of $\Sigma$ with source $f$ and
a whisker $C[\dr f]$ of $\Sigma^\ast$, the pair
$C[b]=(C[A],C[B])$ is a branching of $\Sigma$, with source
$C[f]$. Furthermore, if~$b$ is local,
then $C[b]$ is also local. We define by $\preccurlyeq$ the
order relation on branchings of $\Sigma$ given by
$b\preccurlyeq b'$ when there exists a whisker $C$ such that
$C[b]=b'$ holds.  

A branching is \emph{minimal} when it is minimal for the order relation $\preccurlyeq$. A branching is \emph{trivial} when it can be written either as $(A,A)$, for a $(n+1)$-cell $A$, or as $(A\star_is_n(B),s_n(A)\star_i B)$, for $(n+1)$-cells $A$ and $B$ and a $i$ in $\ens{0,\dots,n-1}$. A branching is \emph{critical} when it is minimal and not trivial.  

In order to prove that $\Sigma$ is locally confluent, it is sufficient to prove that all its critical branchings are confluent. Indeed, trivial branchings are always confluent and a non-minimal branching is confluent if and only if the corresponding minimal branching is (to prove that the latter exists, we proceed by induction on the size of the source of the local branching, which is an $n$-cell in the free $n$-category $\Sigma^*_n$).

\subsection{Using derivations for proving termination of a \pdf{3}-polygraph}

A method to prove termination of a $3$-polygraph has been introduced in~\cite{Guiraud04}, see also~\cite{Guiraud06jpaa,Guiraud06apal}; in special cases, it can also provide complexity bounds~\cite{BonfanteGuiraud09}. It turns out that the method uses interpretations that are a special case of derivations, as described here. Here we only give an outline of the proof.

\subsubsection{Theorem}
\label{Theorem:TerminationDerivation}

\begin{em}
Let $\Sigma$ be a $3$-polygraph such that there exist:
\begin{itemize}
\item Two $2$-functors $X:\Sigma^*_2\fl\Ord$ and $Y:(\Sigma^*_2)^{\text{co}}\fl\Ord$ such that, for every $1$-cell $a$ in $\Sigma_1$, the sets~$X(a)$ and $Y(a)$ are non-empty and, for every $3$-cell $\alpha$ in $\Sigma_3$, the inequalities $X(s\alpha)\geq X(t\alpha)$ and $Y(s\alpha)\geq Y(t\alpha)$ hold.
\item An abelian group $G$ in $\Ord$ whose addition is strictly monotone in both arguments and such that every decreasing sequence of non-negative elements of $G$ is stationary.
\item A derivation $d$ of $\Sigma_2^*$ into the module $M_{X,Y,G}$ such that, for every $2$-cell $f$ in $\Sigma_2^*$, we have $d(f)\geq 0$ and, for every $3$-cell $\alpha$ in $\Sigma_3$, the strict inequality $d(s\alpha)>d(t\alpha)$ holds.
\end{itemize}

\noindent Then the $3$-polygraph $\Sigma$ terminates.
\end{em}

\begin{proof}
Let us assume that $A:f\tfl g$ is a $3$-cell of $\Sigma^*$ with size~$1$. Then there exists a $3$-cell $\alpha:\phi\tfl\psi$ of $\Sigma$ and a context $C$ of $\Sigma_2^*$ such that $A=C[\alpha]$ holds, \ie, such that $f = C[\phi]$ and $g=C[\psi]$ hold. Thus, one gets:
$$
d(f) \:=\: d(C)[\phi] + C[d(\phi)] \qquad\text{and}\qquad d(g) \:=\: d(C)[\psi] + C[d(\psi)].
$$

\noindent We use the fact $d(\phi)>d(\psi)$ holds by hypothesis to get $C[d(\phi)]>C[d(\psi)]$. Moreover, since $X$ and $Y$ are $2$-functors into $\Ord$ and since $d$ sends every $2$-cell to a monotone map, one gets $d(C)[\phi]\geq d(C)[\psi]$. Finally, one uses the hypothesis on the strict monotony of addition in $G$ to get $d(f)>d(g)$. Then one deduces that, for every non-degenerate $3$-cell $A:f\tfl g$, one has $d(f)>d(g)$. Thus, every infinite reduction sequence $(f_k)_k$ would produce an infinite, strictly decreasing sequence $(d(f_k)_k)$ of non-negative elements in $G$, the existence of which is prohibited by hypothesis.
\end{proof}

\subsubsection{Special cases}

The sequel contains several examples where derivations are used to prove termination. Other examples can be found in~\cite{Guiraud06jpaa} or~\cite{BonfanteGuiraud09}. Often, we take the trivial $2$-functor for at least one of the $2$-functors~$X$ and~$Y$ and $\Zb$ for $G$. One can check that those situations match the hypotheses of Theorem~\ref{Theorem:TerminationDerivation}.

\subsection{Branchings and homotopy bases}
\label{Subsection:BranchingHomotopyBases}

In the case of convergent word rewriting systems, \emph{i.e.} convergent $2$-polygraphs with exactly one $0$-cell, the critical branchings generate a homotopy basis~\cite{Squier94}. In this section, we generalise this result to any polygraph. 
In particular, we recover Squier's theorem as Corollary~\ref{Corollary2:BranchingHomotopyBases}, stating that a finite and convergent $2$-polygraph has finite derivation type. However, this result fails to generalise to higher-dimensional polygraphs, as stated in Theorem~\ref{Theorem:mainTheorem}. Indeed, for every $n\geq 3$, there exists at least a finite and convergent $n$-polygraph with an infinite number of critical branchings. The detailed proof can be found in~\ref{Subsection:main_counter_example}. 

\subsubsection{Notation}

When $\Sigma$ is a locally confluent $(n+1)$-polygraph, we assume that, for every critical branching $b=(A,B)$, a confluence $(A',B')$ has been chosen. We denote by $\Gamma_{\Sigma}$ the set of all the $(n+1)$-spheres $(A\star_n A', B\star_n B')$ of $\Sigma$, for each critical branching $b=(A,B)$.

\subsubsection{Lemma}
\label{Lemma1:BranchingHomotopyBases}

\begin{em}
Let $\Sigma$ be a locally confluent $(n+1)$-polygraph. Then every local branching $b=(A,B)$ admits a confluence $(A',B')$ such that $A\star_n A' \approx_{\Gamma_{\Sigma}} B\star_n B'$ holds. 
\end{em}

\begin{proof}
First, let us examine the case where $b$ is a trivial branching. If $A=B$, then $(t_n(A),t_n(B))$ is a confluence that satisfies the required property. Otherwise, let us assume that there exist $(n+1)$-cells $A_1$ and $B_1$ in $\Sigma^\ast$ and an $i$ in $\ens{0,\dots,n-2}$ such that $A=A_1\star_i s_n(B_1)$ and $B=s_n(A_1)\star_i B_1$ hold: then $(t_n(A_1)\star_i B_1, A_1\star_i t_n(B_1))$ is a confluence that satisfies the required property. 

Now, let us assume that $b$ is not trivial. Let $b_1=(A_1,B_1)$ be a minimal branching such that $b_1\preccurlyeq b$, with a whisker $C$ such that $b=C[b_1]$ holds. Since $(A,B)$ is not trivial, then $b_1$ cannot be trivial, so that it is critical. Then we consider its fixed confluence $(A',B')$. Then $(C[A'],C[B'])$ is a confluence for $(A,B)$. Furthermore, one has: 
$$
A\star_n C[A']
	\:=\: C[A_1]\star_n C[A']
	\:=\: C[A_1\star_n A'].
$$

\noindent Similarly, one gets $B\star_n C[B'] = C[B_1\star_n B']$. Since $C$ is a whisker and since, by definition of~$C$, one has $A_1\star_n A'\approx_{\Gamma_{\Sigma}} B_1\star_n B'$, one gets that $(C[A'],C[B'])$ satisfies the required property.  
\end{proof}

\subsubsection{Lemma}
\label{Lemma2:BranchingHomotopyBases}

\begin{em}
Let $\Sigma$ be a convergent $(n+1)$-polygraph and let $(A,B)$ be a branching of $\Sigma$ such that both $t_n(A)$ and $t_n(B)$ are normal forms. Then one has $t_n(A)=t_n(B)$ and $A\approx_{\Gamma_{\Sigma}} B$. 
\end{em}

\begin{proof}
Since $\Sigma$ is terminating, we can prove the result by induction on the source of the branching.

First, if this source $f$ is a normal form, then by definition of normal form, both $A$ and $B$ must be identities. Hence $t_n(A)$ and $t_n(B)$ are equal, and so are $A$ and $B$. Thus $A\approx_{\Gamma_{\Sigma}} B$ holds. 

Now, we fix an $n$-cell $f$, which is not a normal form. We assume that the result holds for every branching $(A,B)$ such that the targets of $A$ and $B$ are normal forms and such that there exists a non-trivial $(n+1)$-cell from $f$ to their source. Let $(A,B)$ be a branching with source $f$ and such that the targets of $A$ and $B$ are normal forms. Since $f$ is not a normal form, $A$ and $B$ cannot be identities, hence one can decompose them into $A=A_1\star_n A_2$ and $B=B_1\star_n B_2$ with $A_1$ and $B_1$ being $(n+1)$-cells of size $1$. 

The pair $(A_1,B_1)$ is a local branching. Thus, using Lemma~\ref{Lemma1:BranchingHomotopyBases}, one gets a confluence $(A'_1,B'_1)$ for $(A_1,B_1)$ such that  $A_1\star_n A'_1\approx_{\Gamma_{\Sigma}} B_1\star_n B'_1$ holds. Let us denote by $g$ the common target of $A'_1$ and $B'_1$, by $e$ its normal form and by $A_3$ an $n$-cell from $g$ to $e$.   

Then we consider the branching $(A_2,A_1'\star_n A_3)$, whose source is denoted by $h$. The targets of $A_2$ and $A_1'\star_n A_3$ are normal forms and $A_1$ is a non-trivial $(n+1)$-cell from $f$ to $h$: thus, the induction hypothesis can be applied to this branching, yielding that $A_2$ has target $e$ and that $A_2\approx_{\Gamma_{\Sigma}} A_1'\star_n A_3$ holds.  

We proceed similarly to prove that $B_2$ satisfies the same properties, so that one gets that $A$ and $B$ have the same target and that $A\approx_{\Gamma_{\Sigma}} B$ holds. The constructions we have done are summarized in the following diagram:  
$$
\xymatrix{
&& h
	\ar[dr] |-{A'_1}
	\ar@/^2ex/[rrrd] ^-{A_2} _-{}="5"
\\
f 
	\ar@/^14ex/[rrrrr] ^-A _-{}="1"
	\ar@/_14ex/[rrrrr] _-B ^-{}="3"
	\ar[urr] ^-{A_1}
	\ar[drr] _-{B_1}
	\ar@{} [rrr] |-{\approx_{\scriptscriptstyle \Gamma_{\Sigma}}}
&&& g  
	\ar[rr] |-{A_3}
&& e
\\
&& k 
	\ar[ru] |-{B'_1}
	\ar@/_2ex/[rrru] _-{B_2} ^-{}="8"
\ar@{} "1" ; "1,3" |-{=}
\ar@{} "3" ; "3,3" |-{=}
\ar@{} "5" ; "2,4" |-{\approx_{\Gamma_{\Sigma}}}
\ar@{} "8" ; "2,4" |-{\approx_{\Gamma_{\Sigma}}}
}
$$
\end{proof}

\subsubsection{Proposition}
\label{Proposition3:BranchingHomotopyBases}

\begin{em}
Let $\Sigma$ be a convergent $(n+1)$-polygraph. Then $\Gamma_{\Sigma}$ is a homotopy basis for $\Sigma^\top$.  
\end{em}

\begin{proof}
Let $(A_1,A_2)$ be an $(n+1)$-sphere in $\Sigma^\top$, with target $n$-cell $f$. Since $\Sigma$ is convergent, we can choose an $(n+1)$-cell $B$ from $f$ to its normal form. Then $(A_1\star_nB,A_2\star_n B)$ satisfies the hypotheses of Lemma~\ref{Lemma2:BranchingHomotopyBases}, yielding $A_1\star_n B \approx_{\Gamma_{\Sigma}} A_2\star_n B$, hence $A_1\approx_{\Gamma_{\Sigma}} A_2$.
\end{proof}

\subsubsection{Proposition}
\label{PropfiniteCriticalBranchingFDT}

\begin{em}
A finite convergent polygraph with a finite set of critical branchings has finite derivation type.  
\end{em}

\begin{proof}
If $\Sigma$ has a finite set of critical branchings, then the set $\Gamma_\Sigma$ is finite. 
\end{proof}

\subsubsection{Corollary}
\label{Corollary1:BranchingHomotopyBases}

\begin{em}
A terminating polygraph with no critical branching has finite derivation type.
\end{em}

\subsubsection{Corollary (\cite{Squier94})} 
\label{Corollary2:BranchingHomotopyBases}

\begin{em}
A finite convergent $2$-polygraph has finite derivation type. 
\end{em}

\begin{proof}
If $\Sigma$ is a finite convergent $2$-polygraph with one $0$-cell, \ie, a word rewriting system, then its set of critical branchings is finite. Indeed, it is equal to the number of possible overlaps between the words corresponding to the sources of $2$-cells: there are finitely many $2$-cells and finitely many letters in each word. If $\Sigma$ has more than one $0$-cell, then the number of possible overlaps is bounded by the number of overlaps in $\Sigma'$, built from $\Sigma$ by identification of all its $0$-cells.
\end{proof}

\noindent From this result Squier has proved that, if a
monoid admits a presentation by a finite convergent word
rewriting system, then it has finite derivation
type~\cite{Squier94}. Now we prove that this result is false
for $n$-categories when $n\geq 2$. 

\subsubsection{Proposition}
\label{Proposition:ofmainTheorem}

\begin{em}
For every natural number $n\geq 3$, there exists a finite convergent $n$-polygraph without finite derivation type. 
\end{em}

\begin{proof}
We consider the $3$-polygraph $\Sigma$ with one $0$-cell, one $1$-cell, three $2$-cells~$\figeps{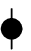}$, $\figeps{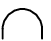}$, $\figeps{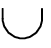}$ and the following four $3$-cells: 
$$
\scalefigeps{-2mm}{on} \tfl \scalefigeps{-2mm}{no} 
\:, \qquad
\scalefigeps{-2mm}{ou} \tfl \scalefigeps{-2mm}{uo}
\:, \qquad
\scalefigeps{-2mm}{nu} \tfl \scalefigeps{-2mm}{id1} 
\:, \qquad
\scalefigeps{-2mm}{un} \tfl \scalefigeps{-2mm}{id1} 
\:.
$$

\noindent The $3$-polygraph $\Sigma$ is finite and convergent. However, the first and second $3$-cells create an infinite number of critical branchings whose confluence diagrams cannot be presented by a finite homotopy basis. These facts are proved in~\ref{Subsection:main_counter_example}.

Then we apply suspension functors on $\Sigma$ to get an $n$-polygraph, for any $n\geq 3$. It has exactly the same cells and compositions in dimensions $n-3$, $n-2$, $n-1$ and $n$ as $\Sigma$ has in dimensions $0$, $1$, $2$ and $3$; on top of that, it has two cells in each dimension up to $n-4$ and no other possible compositions, except with degenerate cells. Thus, we conclude that the $n$-polygraph we have built is finite and convergent, yet it still fails to have finite derivation type.
\end{proof}

\subsubsection{Theorem}
\label{Theorem:mainTheorem}

\begin{em}
For every natural number $n\geq 2$, there exists an
$n$-category which does not have finite derivation type and
admits a presentation by a finite convergent
$(n+1)$-polygraph. 
\end{em} 

\begin{proof}
For any $n\geq 2$, Proposition~\ref{Proposition:ofmainTheorem} implies that there exists a finite convergent $(n+1)$-polygraph~$\Sigma$ without finite derivation type. By Proposition~\ref{TDFTietzeInvariant}, no finite $(n+1)$-polygraph presenting the $n$-category~$\cl{\Sigma}$ can  have finite derivation type. Thus, $\cl{\Sigma}$ does not have finite derivation type.
\end{proof}

\subsubsection{Example}
\label{Example:infinitePolTDF}

We end this section with an example proving that the property of finite derivation type is not Tietze-invariant for \emph{infinite} polygraphs. Let $\Cr$ be the $2$-category presented by the $3$-polygraph~$\Sigma$ with one $0$-cell, one $1$-cell, three $2$-cells $\figeps{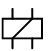}$, $\figeps{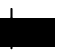}$, $\figeps{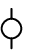}$ and the following two $3$-cells: 
$$
\scalefigeps{-2.8mm}{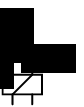} \overset{\alpha}{\tfl}
\scalefigeps{-2.8mm}{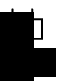}
\qquad\text{and}\qquad 
\scalefigeps{-2.8mm}{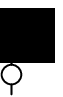} \overset{\beta}{\tfl}
\scalefigeps{-2.8mm}{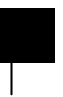}\:. 
$$

\noindent The polygraph $\Sigma$ terminates and does not have critical branching. By Corollary~\ref{Corollary1:BranchingHomotopyBases} it follows that $\Sigma$ has finite derivation type and, thus, so does~$\Cr$. 

Now let us consider another presentation of the $2$-category $\Cr$, namely the $3$-polygraph $\Xi$ defined the same way as~$\Sigma$ except for the $3$-cells:
$$
\scalefigeps{-2.8mm}{exAlphaR} \overset{\alpha}{\tfl}
\scalefigeps{-2.8mm}{exAlphaL}  
\qquad\text{and}\qquad 
\scalefigeps{-2.8mm}{exBetaL} \overset{\beta}{\tfl}
\scalefigeps{-2.8mm}{exBetaR}\:. 
$$

\noindent The $3$-polygraph $\Xi$ still terminates, but it has the following non-confluent critical branching: 
$$
\xymatrix{
{\scalefigeps{-2mm}{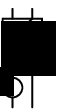}}
	\ar@3 [rr] ^-{\beta}
	\ar@3 [dr] _-{\alpha} ^{}="target"
&& {\scalefigeps{-2mm}{exAlphaR}} 
	\ar@3 [dr] ^-{\alpha}_{}="source" 
\\
& {\scalefigeps{-2mm}{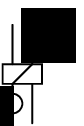}}
&& {\figeps{exAlphaL}} 
}
$$

\noindent We define, by induction on the natural number $k\geq 1$, the $2$-cell $\figeps{exPor1}_k$ as follows:  
$$
\figeps{exPor1}_1 
	\:=\:	\figeps{exPor1}  
\qquad \text{and} \qquad
\figeps{exPor1}_{k+1} 
	\:=\: \raisebox{-1.25mm}{\begin{picture}(0,0)%
\includegraphics{exPor1K.pstex}%
\end{picture}%
\setlength{\unitlength}{4144sp}%
\begingroup\makeatletter\ifx\SetFigFont\undefined%
\gdef\SetFigFont#1#2#3#4#5{%
  \reset@font\fontsize{#1}{#2pt}%
  \fontfamily{#3}\fontseries{#4}\fontshape{#5}%
  \selectfont}%
\fi\endgroup%
\begin{picture}(511,339)(79,557)
\put(406,704){\makebox(0,0)[lb]{\smash{{\SetFigFont{8}{9.6}{\rmdefault}{\mddefault}{\updefault}{\color[rgb]{0,0,0}$k$}%
}}}}
\end{picture}%
} .
$$

\noindent Then, we complete the $3$-polygraph $\Xi$ into an infinite convergent polygraph $\Xi_\infty = \Xi \amalg \ens{\beta_k,\: k\geq 1}$, where~$\beta_0$ is $\beta$ and $\beta_k$ is the following $3$-cell:
$$
\raisebox{-1.25mm}{\begin{picture}(0,0)%
\includegraphics{exBetaKL.pstex}%
\end{picture}%
\setlength{\unitlength}{4144sp}%
\begingroup\makeatletter\ifx\SetFigFont\undefined%
\gdef\SetFigFont#1#2#3#4#5{%
  \reset@font\fontsize{#1}{#2pt}%
  \fontfamily{#3}\fontseries{#4}\fontshape{#5}%
  \selectfont}%
\fi\endgroup%
\begin{picture}(421,474)(169,557)
\put(406,704){\makebox(0,0)[lb]{\smash{{\SetFigFont{8}{9.6}{\rmdefault}{\mddefault}{\updefault}{\color[rgb]{0,0,0}$k$}%
}}}}
\end{picture}%
}
\quad \overset{\beta_k}{\tfl} \quad
\raisebox{-1.25mm}{\begin{picture}(0,0)%
\includegraphics{exBetaKR.pstex}%
\end{picture}%
\setlength{\unitlength}{4144sp}%
\begingroup\makeatletter\ifx\SetFigFont\undefined%
\gdef\SetFigFont#1#2#3#4#5{%
  \reset@font\fontsize{#1}{#2pt}%
  \fontfamily{#3}\fontseries{#4}\fontshape{#5}%
  \selectfont}%
\fi\endgroup%
\begin{picture}(421,358)(169,673)
\put(406,704){\makebox(0,0)[lb]{\smash{{\SetFigFont{8}{9.6}{\rmdefault}{\mddefault}{\updefault}{\color[rgb]{0,0,0}$k$}%
}}}}
\end{picture}%
} .
$$

\noindent The $3$-polygraph $\Xi_{\infty}$ has one confluent critical branching for every natural number $k$:
$$
\xymatrix{
{\raisebox{-1.25mm}{\begin{picture}(0,0)%
\includegraphics{exCritBranchH.pstex}%
\end{picture}%
\setlength{\unitlength}{4144sp}%
\begingroup\makeatletter\ifx\SetFigFont\undefined%
\gdef\SetFigFont#1#2#3#4#5{%
  \reset@font\fontsize{#1}{#2pt}%
  \fontfamily{#3}\fontseries{#4}\fontshape{#5}%
  \selectfont}%
\fi\endgroup%
\begin{picture}(376,609)(169,557)
\put(361,659){\makebox(0,0)[lb]{\smash{{\SetFigFont{8}{9.6}{\rmdefault}{\mddefault}{\updefault}{\color[rgb]{0,0,0}$k$}%
}}}}
\end{picture}%
}}
	\ar@3 [rr] ^-{\raisebox{-1.25mm}{\input{exCritBranchreduc1.pstex_t}}}
	\ar@3 [dr] _-{\raisebox{-1.25mm}{\begin{picture}(0,0)%
\includegraphics{exCritBranchreduc2.pstex}%
\end{picture}%
\setlength{\unitlength}{4144sp}%
\begingroup\makeatletter\ifx\SetFigFontNFSS\undefined%
\gdef\SetFigFontNFSS#1#2#3#4#5{%
  \reset@font\fontsize{#1}{#2pt}%
  \fontfamily{#3}\fontseries{#4}\fontshape{#5}%
  \selectfont}%
\fi\endgroup%
\begin{picture}(339,384)(124,782)
\put(181,884){\makebox(0,0)[lb]{\smash{{\SetFigFontNFSS{8}{9.6}{\rmdefault}{\mddefault}{\updefault}{\color[rgb]{0,0,0}$\beta_k$}%
}}}}
\end{picture}%
}} 
&& {\raisebox{-1.25mm}{\begin{picture}(0,0)%
\includegraphics{exCritBranchR1.pstex}%
\end{picture}%
\setlength{\unitlength}{4144sp}%
\begingroup\makeatletter\ifx\SetFigFont\undefined%
\gdef\SetFigFont#1#2#3#4#5{%
  \reset@font\fontsize{#1}{#2pt}%
  \fontfamily{#3}\fontseries{#4}\fontshape{#5}%
  \selectfont}%
\fi\endgroup%
\begin{picture}(560,474)(169,557)
\put(406,704){\makebox(0,0)[lb]{\smash{{\SetFigFont{8}{9.6}{\rmdefault}{\mddefault}{\updefault}{\color[rgb]{0,0,0}$k+1$}%
}}}}
\end{picture}%
}}
	\ar@3 [dr] ^-{\beta_{k+1}} 
\\
& {\raisebox{-1.25mm}{\begin{picture}(0,0)%
\includegraphics{exCritBranchR2.pstex}%
\end{picture}%
\setlength{\unitlength}{4144sp}%
\begingroup\makeatletter\ifx\SetFigFont\undefined%
\gdef\SetFigFont#1#2#3#4#5{%
  \reset@font\fontsize{#1}{#2pt}%
  \fontfamily{#3}\fontseries{#4}\fontshape{#5}%
  \selectfont}%
\fi\endgroup%
\begin{picture}(376,538)(169,628)
\put(361,659){\makebox(0,0)[lb]{\smash{{\SetFigFont{8}{9.6}{\rmdefault}{\mddefault}{\updefault}{\color[rgb]{0,0,0}$k$}%
}}}}
\end{picture}%
}}
	\ar@3 [rr] _-{\raisebox{-1.25mm}{\input{exCritBranchR2reduc.pstex_t}}}
&& {\raisebox{-1.25mm}{\begin{picture}(0,0)%
\includegraphics{exCritBranchN.pstex}%
\end{picture}%
\setlength{\unitlength}{4144sp}%
\begingroup\makeatletter\ifx\SetFigFont\undefined%
\gdef\SetFigFont#1#2#3#4#5{%
  \reset@font\fontsize{#1}{#2pt}%
  \fontfamily{#3}\fontseries{#4}\fontshape{#5}%
  \selectfont}%
\fi\endgroup%
\begin{picture}(560,358)(169,673)
\put(406,704){\makebox(0,0)[lb]{\smash{{\SetFigFont{8}{9.6}{\rmdefault}{\mddefault}{\updefault}{\color[rgb]{0,0,0}$k+1$}%
}}}}
\end{picture}%
}}
\ar@4{ ->} "1,3" ; "2,2" |-*+[o]{\scriptstyle\alpha\beta_k}
}
$$
By Proposition \ref{Proposition3:BranchingHomotopyBases}, the set $\Gamma = \ens{ \alpha\beta_k \;|\; k\in\Nb }$ is a homotopy basis of the $3$-category $\Xi_\infty^\top$. 

Let us prove that the $3$-polygraph  $\Xi_\infty$ does not have finite derivation type. On the contrary, let us assume that $\Xi_\infty$ has finite derivation type. Then, following Proposition~\ref{PropExtractFiniteHomotopyBase}, there exists a finite subset~$\Gamma_0$ of~$\Gamma$ which is a homotopy base of $\Xi_\infty^\top$. Thus, there exists a natural number~$l$ such that, for every $k\geq l$, the $4$-cell $\alpha\beta_k$ is not in $\Gamma_0$. However, since $\Gamma_0$ is a homotopy base we still have: 
$$
s\left(\alpha\beta_l\right) \approx_{\Gamma_0}
t\left(\alpha\beta_l\right).
$$
Hence, there exists a $4$-cell $\Phi$ in $\Xi^\top_\infty(\Gamma_0)$ such that $s\Phi = s\left(\alpha\beta_l\right)$ and $t\Phi = t\left(\alpha\beta_l\right)$ hold. Let us prove that this is not possible, thanks to the derivation $d$ of $\Xi^{\top}_{\infty}$ into the trivial module given by:
$$
d(\alpha) \:=\: 0 
\qquad\text{and}\qquad 
d(\beta_k) \:=\: 
\begin{cases} 
	0 &\text{if }k\leq l, \\
	1 &\text{if }k\geq l+1.
\end{cases}
$$
Then, for every $k\leq l$, we have $d(s(\alpha\beta_k))=d(t(\alpha\beta_k))=0$. As a consequence, for every $4$-cell $\Psi$ in~$\Xi^{\top}_{\infty}(\Gamma_0)$, we have $d(s\Psi)=d(t\Psi)$. In particular, when $\Psi=\Phi$, we get $d(s(\alpha\beta_l))=d(t(\alpha\beta_l))$. This is not possible since, by definition of $d$, we have $d(s(\alpha\beta_l))=1$ and $d(t(\alpha\beta_l))=0$. This proves that~$\Xi_\infty$ does not have finite derivation type.



\section{The case of \pdf{3}-polygraphs}
\label{Section:CaseOf3Polygraphs}

\subsection{Classification of critical branchings}

\subsubsection{Types of critical branchings}

Let $\Sigma$ be a $3$-polygraph and let $(A,B)$ be a critical branching of $\Sigma$. Let us denote by~$\alpha$ and $\beta$ the $3$-cells of $\Sigma$ that generate $A$ and $B$. Then $(A,B)$ falls in one of three cases.

The first possibility is that there exists a context $C$ of $\Sigma_2^*$ such that $s\alpha=C[s\beta]$ holds. Then, the source of the branching $(A,B)$ is: 
$$
\scalebox{1.3}{\raisebox{-4mm}{\raisebox{-1.25mm}{\begin{picture}(0,0)%
\includegraphics{pc-inclusion-a.pstex}%
\end{picture}%
\setlength{\unitlength}{4144sp}%
\begingroup\makeatletter\ifx\SetFigFont\undefined%
\gdef\SetFigFont#1#2#3#4#5{%
  \reset@font\fontsize{#1}{#2pt}%
  \fontfamily{#3}\fontseries{#4}\fontshape{#5}%
  \selectfont}%
\fi\endgroup%
\begin{picture}(474,564)(529,-568)
\put(766,-331){\makebox(0,0)[b]{\smash{{\SetFigFont{8}{9.6}{\familydefault}{\mddefault}{\updefault}$s\alpha$}}}}
\end{picture}%
}}} 
\quad=\quad 
\scalebox{1.3}{\raisebox{-4mm}{\raisebox{-1.25mm}{\begin{picture}(0,0)%
\includegraphics{pc-inclusion-b.pstex}%
\end{picture}%
\setlength{\unitlength}{4144sp}%
\begingroup\makeatletter\ifx\SetFigFont\undefined%
\gdef\SetFigFont#1#2#3#4#5{%
  \reset@font\fontsize{#1}{#2pt}%
  \fontfamily{#3}\fontseries{#4}\fontshape{#5}%
  \selectfont}%
\fi\endgroup%
\begin{picture}(501,564)(529,-568)
\put(766,-331){\makebox(0,0)[b]{\smash{{\SetFigFont{8}{9.6}{\familydefault}{\mddefault}{\updefault}$s\beta$}}}}
\put(923,-489){\makebox(0,0)[b]{\smash{{\SetFigFont{8}{9.6}{\familydefault}{\mddefault}{\updefault}$C$}}}}
\end{picture}%
}}} \:.
$$

\noindent In that case, $(A,B)$ is an \emph{inclusion} critical branching.

If the branching $(A,B)$ is not an inclusion one, the second possibility is that there exist $1$-cells $u$, $v$ and $2$-cells $f$, $g$, $h$ such that $s\alpha$ and $s\beta$ decompose in one of the following ways.

\begin{itemize}

\item One has $s\alpha=f\star_1 (u\star_0 h)$ and $s\beta=(h\star_0 v)\star_1 g$, so that the source of $(A,B)$ is:
$$
\scalebox{1.3}{\raisebox{-4mm}{\raisebox{-1.25mm}{\input{pc-reg-1-a.pstex_t}}}}
\quad=\quad
\scalebox{1.3}{\raisebox{-4mm}{\raisebox{-1.25mm}{\input{pc-reg-1-src.pstex_t}}}}
\quad=\quad
\scalebox{1.3}{\raisebox{-4mm}{\raisebox{-1.25mm}{\input{pc-reg-1-b.pstex_t}}}} \:.
$$

\item One has $s\alpha=f\star_1 (h\star_0 u)$ and $s\beta=(v\star_0 h)\star_1 g$:
$$
\scalebox{1.3}{\raisebox{-4mm}{\raisebox{-1.25mm}{\input{pc-reg-3-a.pstex_t}}}}
\quad=\quad
\scalebox{1.3}{\raisebox{-4mm}{\raisebox{-1.25mm}{\input{pc-reg-3-src.pstex_t}}}}
\quad=\quad
\scalebox{1.3}{\raisebox{-4mm}{\raisebox{-1.25mm}{\input{pc-reg-3-b.pstex_t}}}} \:.
$$

\item One has $s\alpha=f\star_1 (u\star_0 h\star_0 v)$ and $s\beta=h\star_1 g$:
$$
\scalebox{1.3}{\raisebox{-4mm}{\raisebox{-1.25mm}{\input{pc-reg-2-a.pstex_t}}}}
\quad=\quad
\scalebox{1.3}{\raisebox{-4mm}{\raisebox{-1.25mm}{\input{pc-reg-2-src.pstex_t}}}}
\quad=\quad
\scalebox{1.3}{\raisebox{-4mm}{\raisebox{-1.25mm}{\input{pc-reg-2-b.pstex_t}}}} \:.
$$

\item One has $s\alpha=f\star_1 h$ and $s\beta=(u\star_0 h\star_0 v)\star_1 g$:
$$
\scalebox{1.3}{\raisebox{-4mm}{\raisebox{-1.25mm}{\input{pc-reg-4-a.pstex_t}}}}
\quad=\quad
\scalebox{1.3}{\raisebox{-4mm}{\raisebox{-1.25mm}{\input{pc-reg-4-src.pstex_t}}}}
\quad=\quad
\scalebox{1.3}{\raisebox{-4mm}{\raisebox{-1.25mm}{\input{pc-reg-4-b.pstex_t}}}} \:.
$$

\end{itemize}

\noindent If $(A,B)$ matches one of these cases, then it is called a \emph{regular} critical branching.

Finally, when the branching $(A,B)$ is not an inclusion or regular one, there exist $1$-cells $u$, $v$ and $2$-cells $f$, $g$, $h$ such $s\alpha$ and $s\beta$ decompose in one of the following ways.

\begin{itemize}

\item One has $s\alpha=f\star_1 (h\star_0 u)$ and $s\beta=(h\star_0 v)\star_1 g$, so that there exists a $2$-cell $k$ such that the source of $(A,B)$ is:
$$
\scalebox{1.3}{\raisebox{-4mm}{\raisebox{-1.25mm}{\input{pc-ind-r-a.pstex_t}}}}
\quad=\quad
\scalebox{1.3}{\raisebox{-4mm}{\raisebox{-1.25mm}{\input{pc-ind-r-src.pstex_t}}}}
\quad=\quad
\scalebox{1.3}{\raisebox{-4mm}{\raisebox{-1.25mm}{\input{pc-ind-r-b.pstex_t}}}} \:.
$$

\noindent In that case, one can write $(A,B)=(C[k],D[k])$ for appropriate contexts $C$ and $D$ of $\Sigma^*$. The family $(C[k],D[k])_k$, where $k$ ranges over the $2$-cells with appropriate boundary and such that $(C[k],D[k])$ is a minimal branching, is called a \emph{right-indexed} critical branching.

\item One has $s\alpha=f\star_1 (u\star_0 h)$ and $s\beta=(v\star_0 h)\star_1 g$, so that there exists a $2$-cell $k$ such that the source of $(A,B)$ is:
$$
\scalebox{1.3}{\raisebox{-4mm}{\raisebox{-1.25mm}{\input{pc-ind-l-a.pstex_t}}}}
\quad=\quad
\scalebox{1.3}{\raisebox{-4mm}{\raisebox{-1.25mm}{\input{pc-ind-l-src.pstex_t}}}}
\quad=\quad
\scalebox{1.3}{\raisebox{-4mm}{\raisebox{-1.25mm}{\input{pc-ind-l-b.pstex_t}}}} \:.
$$

\noindent In that case, one can write $(A,B)=(C[k],D[k])$ for appropriate contexts $C$ and $D$ of $\Sigma^*$. The family $(C[k],D[k])_k$, where $k$ ranges over the $2$-cells with appropriate boundary and such that $(C[k],D[k])$ is a minimal branching, is called a \emph{left-indexed} critical branching.

\item One is not in the right-indexed or left-indexed cases and one has
$$
s\alpha \:=\: f\star_1(u_0 \star_0 h_1\star_0 u_1\star_0 h_2\star_0 \cdots \star_0 u_{n-1} \star_0 h_n \star_0 u_n) 
$$

\noindent and
$$
s\beta \:=\: (v_0\star_0 h_1\star_0 v_1\star_0 h_2\star_0 \cdots \star_0 v_{n-1} \star_0 h_n \star_0 v_n)\star_1 g \:,
$$

\noindent so that there exist $2$-cells $k_0$, $\dots$, $k_n$ such that the source of $(A,B)$ is as follows, where we write~$p$ instead of~$n-1$ for size reasons:
\begin{align*}
& \quad \scalebox{1.3}{\raisebox{-4mm}{\raisebox{-1.25mm}{\input{pc-ind-m-a.pstex_t}}}} \\
=& \quad \scalebox{1.3}{\raisebox{-4mm}{\raisebox{-1.25mm}{\input{pc-ind-m-src.pstex_t}}}} \\
=& \quad \scalebox{1.3}{\raisebox{-4mm}{\raisebox{-1.25mm}{\input{pc-ind-m-b.pstex_t}}}} \:. 
\end{align*}

\noindent In that case, one can write $(A,B)=(C[k_0,\dots,k_n],D[k_0,\dots,k_n])$ for appropriate $3$-cells~$C$ and $D$ in some $\Sigma^*[x_0,\dots,x_n]$. The family $(C[k_0,\dots,k_n],D[k_0,\dots,k_n]) _{k_0,\dots,k_n}$, where the $k_i$'s range over the $2$-cells with appropriate boundary and such that $(C[k_0,\dots,k_n],D[k_0,\dots,k_n])$ is a minimal branching, is called a \emph{multi-indexed} critical branching.

\end{itemize}

\noindent In all those indexed cases, the branching $(A,B)$ is said to be an \emph{instance} of the corresponding right-indexed or left-indexed or multi-indexed one. It is a \emph{normal} instance when the indexing $2$-cell $k$ (resp. $2$-cells $k_0$, $\dots$, $k_n$) is a normal form (resp. are normal forms).

\subsubsection{Definitions}
\label{SubSubSection:FinitelyIndexed}

A $3$-polygraph is \emph{non-indexed} when each of its critical branchings is an inclusion one or a regular one. It is \emph{right-indexed} (resp. \emph{left-indexed}) when each of its critical branchings is either an inclusion one, a regular one or an instance of a right-indexed (resp. left-indexed) one.
A $3$-polygraph is \emph{finitely indexed} when each of its indexed critical branchings has a finite number of normal instances.

\subsubsection{Proposition}
\label{Proposition:InclusionRegular}

\begin{em}
A $3$-polygraph with a finite set of $3$-cells has a finite number of inclusion and regular critical branchings.
\end{em}

\begin{proof}
Let $\Sigma$ be a $3$-polygraph with $\Sigma_3=\ens{\alpha_1,\dots,\alpha_p}$ finite. As a consequence, for any  $i,j\in\ens{1,\dots,p}$, the set of morphisms from $s\alpha_i$ to $s\alpha_j$ in $\whisk{\Sigma}$ is finite. Thus $\Sigma$ has a finite number of inclusion branchings.

Now, let us fix $i,j\in\ens{1,\dots,p}$ and let us assume that there exist two whiskers $C$ and $D$ of $\Sigma^*$ such that the pair $(C[\alpha_i], D[\alpha_j])$ is a regular branching, with source $f$. Then there exist a $2$-cell $h$ and whiskers~$C'$ and $D'$ of $\Sigma^*$ that satisfy $C[s\alpha_i]=C'[h]=D'[h]=D[s\alpha_j]$. Since the sets $\whisk{\Sigma}(s\alpha_i,f)$, $\whisk{\Sigma}(s\alpha_j,f)$, $\whisk{\Sigma}(h,C[s\alpha_i])$ and $\whisk{\Sigma}(h,C[s\alpha_j])$ are finite, there exist finitely many regular branchings of this form, with $i$, $j$ fixed. Since $\Sigma_3$ is finite, the $3$-polygraph $\Sigma$ has finitely many regular branchings.
\end{proof}

\subsubsection{Theorem}
\label{Thm:3PolyNonIndexe}

\begin{em}
A finite, convergent, non-indexed $3$-polygraph has finite derivation type. 
\end{em}

\begin{proof}
We use Proposition~\ref{Proposition:InclusionRegular} and, then, we apply Proposition~\ref{PropfiniteCriticalBranchingFDT}.
\end{proof}

\subsection{Mac Lane's coherence theorem revisited}
\label{Subsection:PolyMon}

\subsubsection{Monoidal categories}

A \emph{monoidal category} is a data $(\Cr,\tens,e,a,l,r)$ made of a category $\Cr$,  a bifunctor $\tens:\Cr\times\Cr\fl\Cr$, an object $e$ of $\Cr$ and three natural isomorphisms
$$
a_{x,y,z} \: : \: (x\tens y)\tens z \:\fl\: x\tens(y\tens z) \:,
\qquad
l_x \::\: e \tens x \:\fl\: x,
\qquad
r_x \::\: x\tens e \:\fl\: x,
$$

\noindent such that the following two diagrams commute in $\Cr$:
$$
\xymatrix@C=-1em{
& {\scriptstyle (x\tens (y\tens z))\tens t}
	\ar[rr] ^-{a}
&& {\scriptstyle x\tens ((y\tens z) \tens t)}
	\ar[dr] ^-{a}
\\
{\scriptstyle ((x\tens y)\tens z) \tens t}
	\ar[ur] ^-{a}
	\ar[drr] _-{a}
&& {\scriptstyle \copyright}
&& {\scriptstyle x\tens(y\tens (z\tens t))}
\\
&& {\scriptstyle (x\tens y)\tens (z\tens t)}
	\ar[urr] _-{a}
}
\qquad
\raisebox{-1.5em}{
\xymatrix@C=1em{
& {\scriptstyle x\tens (e\tens y)}
	\ar[dr] ^-{l}
\\
{\scriptstyle (x\tens e)\tens y}
	\ar[ur] ^-{a}
	\ar[rr] _-{r} ^-{}="1"
&& {\scriptstyle x\tens y}
\ar@{} "1,2" ; "1" |-{\copyright}
}
}
$$

\noindent Mac Lane's coherence theorem~\cite{MacLane98} states that, in a such monoidal category, all the diagrams whose arrows are built from $\tens$, $e$, $l$ and $r$ commute. Thereafter, we give a proof of this fact by building a homotopy basis of a $3$-polygraph.

\subsubsection{The \pdf{3}-polygraph of monoids}

We denote by $\Sigma$ the $3$-polygraph with one $0$-cell, one $1$-cell, two $2$-cells~$\figeps{m}$ and~$\figeps{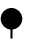}$ and the following three $3$-cells:
$$
\figeps{ass-src} \:\tfl_{\alpha}\: \figeps{ass-tgt} \:,
\qquad\qquad
\figeps{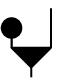} \:\tfl_{\lambda}\: \figeps{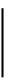} \:,
\qquad\qquad
\figeps{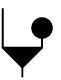} \:\tfl_{\rho}\: \figeps{un-tgt} \:.
$$

\noindent We denote by $\Gamma$ the set made of the following $4$-cells $\alpha\alpha$ and $\alpha\rho$, where we commit the abuse of denoting a $3$-cell of $\Sigma^*$ with size $1$ like its generating $3$-cell:
$$
\xymatrix@C=1em@R=0.9em{
& {\figeps{ass-ass-2}}
	\ar@3 [rr] ^-{\alpha} _-{}="1" 
&& {\figeps{ass-ass-3}}
	\ar@3 [dr] ^-{\alpha} 
\\
{\figeps{ass-ass-1}}
	\ar@3 [ur] ^-{\alpha}
	\ar@3 [drr] _-{\alpha}
&&&& {\figeps{ass-ass-5}}
\\
&& {\figeps{ass-ass-4}}
	\ar@3 [urr] _-{\alpha}
	\ar@4{ ->} "1"; |-*+[o]{\scriptstyle\alpha\alpha}
}
\qquad
\raisebox{-1em}{
\xymatrix{
& {\figeps{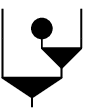}}
	\ar@3 [dr] ^-{\lambda}
\\
{\figeps{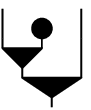}}
	\ar@3 [ur] ^-{\alpha}
	\ar@3 [rr] _-{\rho} ^-{}="1"
&& {\figeps{m}}
\ar@4{ ->} "1,2" ; "1" |-*+[o]{\scriptstyle\alpha\rho}
}
}
$$

\subsubsection{Theorem}

\begin{em}
The set $\Gamma$ of $4$-cells forms a homotopy basis of the track $3$-category $\Sigma^{\top}$.
\end{em}

\begin{proof}
Let us prove that $\Sigma$ terminates. We consider the
$\Sigma_2^*$-module $M_{X,\ast,\Zb}$ and the derivation $d$
of $\Sigma_2^*$ into $M_{X,\ast,\Zb}$ generated by the
following values:
$$
X(\:\figeps{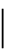}\:) \:=\: \Nb\setminus\ens{0} , 
\quad
X(\figeps{m})(i,j) \:=\: i+j,
\quad
X(\figeps{i}) \:=\: 1,
$$
$$
d\left(\figeps{m}\right) (i,j) \:=\: i,
\quad
d\left(\figeps{i}\right) \:=\: 0.
$$

\noindent We check that the $2$-functor $X$ satisfies the (in)equalities
$$
X\left(\scalefigeps{-2.8mm}{ass-src}\right) (i,j,k) 
	\:=\: i+j+k 
	\:=\: X\left(\scalefigeps{-2.8mm}{ass-tgt}\right) (i,j,k) \;,
$$
$$
X\left(\scalefigeps{-2.8mm}{ung-src}\right) (i) 
	\:=\: i  
	\:=\: X\left(\:\scalefigeps{-2.8mm}{un-tgt}\:\right)(i) \;,
\qquad\qquad
X\left(\scalefigeps{-2.8mm}{und-src}\right) (i) 
	\:=\: i  
	\:=\: X\left(\:\scalefigeps{-2.8mm}{un-tgt}\:\right)(i)
$$

\noindent and that the derivation $d$ satisfies the strict inequalities
$$
d\left(\scalefigeps{-2.8mm}{ass-src}\right) (i,j,k) 
	\:=\: 2i+j 
	\:>\: i+j 
	\:=\: d\left(\scalefigeps{-2.8mm}{ass-tgt}\right)(i,j,k) \;,
$$
$$
d\left(\scalefigeps{-2.8mm}{ung-src}\right) (i) 
	\:=\: 1 
	\:>\: 0 
	\:=\: d\left(\:\scalefigeps{-2.8mm}{un-tgt}\:\right)(i) \;,
\qquad\qquad
d\left(\scalefigeps{-2.8mm}{und-src}\right) (i) 
	\:=\: i 
	\:>\: 0 
	\:=\: d\left(\:\scalefigeps{-2.8mm}{un-tgt}\:\right)(i) \;.
$$
We apply Theorem~\ref{Theorem:TerminationDerivation} to get termination.

The $3$-polygraph $\Sigma$ has five critical branchings. All of them are regular ones and confluent. Their confluence diagrams are given by the boundaries of the two $4$-cells of $\Gamma$ and of the following three ones:
$$
\xymatrix{
& {\figeps{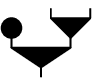}}
	\ar@3 [dr] ^-{\lambda}
\\
{\figeps{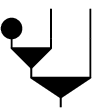}}
	\ar@3 [ur] ^-{\alpha}
	\ar@3 [rr] _-{\lambda} ^-{}="1"
&& {\figeps{m}}
\ar@4{ ->} "1,2" ; "1" |-*+[o]{\scriptstyle\lambda\alpha}
}
\qquad
\xymatrix{
& {\figeps{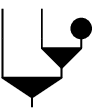}}
	\ar@3 [dr] ^-{\rho}
\\
{\figeps{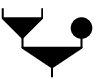}}
	\ar@3 [ur] ^-{\alpha}
	\ar@3 [rr] _-{\rho} ^-{}="1"
&& {\figeps{m}}
\ar@4{ ->} "1,2" ; "1" |-*+[o]{\scriptstyle\rho\alpha}
}
\qquad
\raisebox{-1.5em}{
\xymatrix{
{\figeps{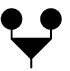}}
	\ar@3@/^6ex/ [rr] ^-{\lambda} _-{\hole}="1"
	\ar@3@/_6ex/ [rr] _-{\rho} ^-{\hole}="2"
&& {\figeps{i}}
\ar@4{ ->} "1"; "2" |-*+[o]{\scriptstyle\lambda\rho}
}
}
$$

\noindent Since $\Sigma$ terminates and has all its critical branchings confluent, it is convergent as a consequence of Newman's lemma. Thus we know that the set $\ens{ \alpha\alpha, \alpha\rho, \lambda\rho, \lambda\alpha, \rho\alpha }$ of $4$-cells is a (finite) homotopy basis of $\Sigma^{\top}$. To get the result, we check that $\lambda\alpha$, $\rho\alpha$ and $\lambda\rho$ are superfluous in this homotopy basis, \ie, that their boundaries are also the ones of $4$-cells of~$\Sigma^{\top}(\Gamma)$. 

For $\lambda\alpha$, we consider the $4$-cell $\left(\figeps{i}\:\figeps{i}\:\:\figeps{1-cell}\:\:\figeps{1-cell}\:\right) \star_1\alpha\alpha$ which is in $\Sigma^{\top}(\Gamma)$. We partially fill its boundary with other~$4$-cells of $\Sigma^{\top}(\Gamma)$ and equalities, yielding a $3$-sphere of $\Sigma^{\top}$ denoted by $\gamma$:
$$
\xymatrix{
&& {\figeps{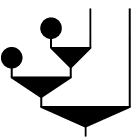}}
	\ar@3 [rr] ^-{\alpha} 
	\ar@3 [d] |-*+[o]{\scriptstyle\lambda}
	\ar@{} [drr] |-{=}
&& {\figeps{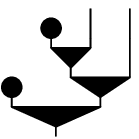}}
	\ar@3 [drr] ^-{\alpha} _-{}="2"
	\ar@3 [d] |-*+[o]{\scriptstyle\lambda} 
\\
{\figeps{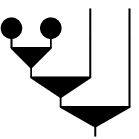}}
	\ar@3 [urr] ^-{\alpha} _-{}="1"
	\ar@3 [rr] |-*+[o]{\scriptstyle\rho} 
	\ar@3@/_3ex/ [drrr] _-{\alpha}
&&	{\figeps{ass-ung-1}}
	\ar@3 [rr] |-*+[o]{\scriptstyle\alpha}
	\ar@{} [dr] |-{=}
&& {\figeps{ass-ung-2}}
&& {\figeps{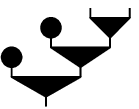}}
	\ar@3 [ll] |-*+[o]{\scriptstyle\lambda}
\\
&&& {\figeps{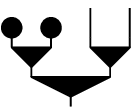}}
	\ar@3 [ur] |-*+[o]{\scriptstyle\rho} 
	\ar@3@/_3ex/ [urrr] _-{\alpha} ^-{}="3"
\ar@{} "1" ; "2,3" |-{\alpha\rho}
\ar@{} "2" ; "2,5" |-{\gamma}
\ar@{} "3" ; "2,5" |-{\alpha\rho}
}
$$

\noindent As a consequence of this construction, we have $s\gamma\approx_{\Gamma}t\gamma$. Then we build the following diagram, proving that $s(\lambda\alpha) \approx_{\Gamma} t(\lambda\alpha)$ also holds:
$$
\xymatrix@R=1.5em{
&&& {\figeps{pent-ee-3}} 
	\ar@3@/_3ex/ [dddlll] _-{\lambda} ^-{}="s1"
	\ar@3 [d] ^-{\lambda}
	\ar@3@/^3ex/ [dddrrr] ^-{\alpha} _-{}="s2"
\\ 
&&& {\figeps{ass-ung-1}}
	\ar@3 [dl] |-*+[o]{\scriptstyle\lambda} _-{}="t1"
	\ar@3 [dr] |-*+[o]{\scriptstyle\alpha} ^-{}="t2"
\\ 
&& {\figeps{m}}
&& {\figeps{ass-ung-2}}
	\ar@3 [ll] |-*+[o]{\scriptstyle\lambda} ^-{}="s3"
\\
{\figeps{ass-ung-2}}
	\ar@3 [urr] |-*+[o]{\scriptstyle\lambda}
&&&&&& {\figeps{pent-ee-5}}
	\ar@3 [ull] |-*+[o]{\scriptstyle\lambda}
	\ar@3 [llllll] ^-{\lambda} _-{}="t3"
\ar@{} "s1" ; "t1" |-{=}
\ar@{} "s2" ; "t2" |-{=}
\ar@{} "s3" ; "t3" |-{=}
}
$$

\noindent For the $4$-cell $\rho\alpha$, one proceeds in a similar way, starting with the $4$-cell $\left(\:\figeps{1-cell}\:\:\figeps{1-cell}\:\:\figeps{i}\:\figeps{i}\: \right) \star_1\alpha\alpha$. 

Finally, let us consider the case of the $4$-cell $\lambda\rho$. First, we consider the $3$-cell $\figeps{i}\star_1\rho\star_1\rho$. Thanks to the exchange relation between $\star_1$ and $\star_2$, we decompose this $3$-cell in two ways. This yields a (trivial) $3$-sphere that we partially fill, using $4$-cells of $\Sigma^{\top}(\Gamma)$, as follows, producing another $3$-sphere $\delta$ of $\Sigma^{\top}(\Gamma)$:
$$
\xymatrix{
{\figeps{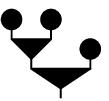}}
	\ar@3 @/^9ex/ [rrrr] ^-{\rho} _-{}="s1"
	\ar@3 [rr] |-*+[o]{\scriptstyle\alpha} ^-{}="t1" _-{}="s2"
	\ar@3 @/_9ex/ [rrrr] _-{\rho} ^-{}="t2"
&& {\figeps{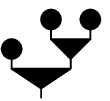}}
	\ar@3 @/^3ex/ [rr] ^(0.4){\scriptstyle\lambda} _-{}="s3"
	\ar@3 @/_3ex/ [rr] _(0.4){\scriptstyle\rho} ^-{}="t3"
&& {\figeps{ung-und}}
	\ar@3 [rr] ^-{\rho}
&& {\figeps{i}}
\ar@{} "s1";"t1" |-{\scriptstyle\alpha\rho} 
\ar@{} "s2";"t2" |-{\scriptstyle\rho\alpha}
\ar@{} "s3";"t3" |-{\scriptstyle\delta}
}
$$

\noindent As a consequence, we have $s\delta\approx_{\Gamma}t\delta$, hence $s\delta\star_2(\figeps{i}\star_1\lambda) \:\approx_{\Gamma}\: t\delta\star_2(\figeps{i}\star_1\lambda)$. The following diagram yields $s(\lambda\rho)\approx_{\Gamma}t(\lambda\rho)$, thus concluding the proof:
$$
\xymatrix{
&& {\figeps{ung-und}} 
	\ar@3 @/^3ex/ [drr] ^-{\lambda}
\\
{\figeps{1-11}}
	\ar@3 @/^3ex/ [urr] ^-{\lambda} _-{}="s1"
	\ar@3 [rr] |-*+[o]{\scriptstyle\lambda} 
	\ar@3 @/_3ex/ [drr] _-{\rho} ^{}="t2"
&& {\figeps{ung-und}}
	\ar@3 @/^3ex/ [rr] ^(0.4){\scriptstyle\lambda}
	\ar@3 @/_3ex/ [rr] _(0.4){\scriptstyle\rho}
&& {\figeps{i}}
\\
&& {\figeps{ung-und}}
	\ar@3 @/_3ex/ [urr] _-{\lambda}
\ar@{} "s1";"2,3" |-{=}
\ar@{} "2,3";"t2" |-{=}
}
$$
\end{proof}

\subsubsection{Corollary (Mac Lane's coherence theorem~\cite{MacLane98})}

\begin{em}
In a monoidal category $(\Cr,\tens,e,a,l,r)$, all the diagrams whose arrows are built from $\tens$, $e$, $a$, $l$ and $r$ are commutative.
\end{em}

\begin{proof}
We see $\Cat{1}$ as a (large) $3$-category with one $0$-cell, categories as $1$-cells, functors as $2$-cells and natural transformations as $3$-cells. The $0$-composition is the cartesian product of categories, the $1$-composition is the composition of functors and the $2$-composition is the "vertical" composition of natural transformations.

Then monoidal categories are exactly the $3$-functors from $\Sigma^{\top}/\Gamma$ to $\Cat{1}$. The correspondence between a monoidal category $(\Cr,\tens,e,a,l,r)$ and such a $3$-functor $M$ is given by:
$$
M(\:\figeps{1-cell}\:) \:=\: \Cr,
\quad
M(\figeps{m}) \:=\: \tens,
\quad
M(\figeps{i}) \:=\: e,
\quad
M(\alpha) \:=\: a,
\quad
M(\lambda) \:=\: l,
\quad
M(\rho) \:=\: r.
$$

\noindent As a consequence, a diagram $D$ in $\Cr$ whose arrows are built from $\tens$, $e$, $a$, $l$ and $r$ is the image by $M$ of a $3$-sphere $\gamma$ of $\Sigma^{\top}$. Since $\Gamma$ is a homotopy basis of $\Sigma^{\top}$, we have $s\gamma\approx_{\Gamma} t\gamma$. Since $M$ is a $3$-functor from $\Sigma^{\top}/\Gamma$ to $\Cat{1}$, we have $M(s\gamma) = M(t\gamma)$, which means that the diagram $D=M(\gamma)$ commutes.
\end{proof}

\subsubsection{Remark}

The definition of monoidal category we have given is minimal, in the sense that both coherence diagrams are required in order to get Mac Lane's coherence theorem. Otherwise, this would mean that either $\alpha\alpha$ or $\alpha\rho$ is superfluous in the homotopy basis $\Gamma$ of $\Sigma^{\top}$. Let us prove that this is not the case. Let $d_1$ be the derivation of $\Sigma^{\top}$ into the trivial module given by:
$$
d_1(\alpha) \:=\: 0,
\qquad
 d_1(\lambda) \:=\: 1,
\qquad 
d_1(\rho) \:=\: 0.
$$

\noindent Then we have $d_1(s\alpha\alpha)=d_1(t\alpha\alpha)=0$. As a consequence, for every $4$-cell $A$ in $\Sigma^{\top}(\alpha\alpha)$, we have $d_1(sA)=d_1(tA)$. Thus, if $\ens{\alpha\alpha}$ was a homotopy basis of $\Sigma^{\top}$, we would have $d_1(F)=d_1(G)$ for every $3$-sphere $(F,G)$ of $\Sigma^{\top}$. In particular, $d_1(s\alpha\rho)=d_1(t\alpha\rho)$ would be true. This is impossible since we have $d_1(s\alpha\rho)=1$ and $d_1(t\alpha\rho)=0$. 

In order to prove that $\ens{\alpha\rho}$ is not a homotopy basis either, we proceed similarly with the derivation~$d_2$ of $\Sigma^{\top}$ into the trivial module given by:
$$
d_2(\alpha) \:=\: 1,
\qquad
 d_2(\lambda) \:=\: -1,
\qquad 
d_2(\rho) \:=\: 0.
$$

\noindent We check that $d_2(s\alpha\rho)=d_2(t\alpha\rho)=0$ holds. Thus, if $\ens{\alpha\rho}$ was a homotopy basis of $\Sigma^{\top}$, the equality $d_2(s\alpha\alpha)=d_2(t\alpha\alpha)$ would be satisfied. However, we have $d_2(s\alpha\alpha)=3$ and $d_2(t\alpha\alpha)=2$.

\subsection{Right-indexed and left-indexed \pdf{3}-polygraphs}

\subsubsection{Proposition}
\label{Theorem:RightIndexedConfluence}

\begin{em}
Let $\Sigma$ be a terminating right-indexed (resp. left-indexed) $3$-polygraph. Then $\Sigma$ is confluent if and only if every inclusion critical branching, every regular critical branching and every instance of every right-indexed (resp. left-indexed) critical branching is confluent.
\end{em}

\begin{proof}
If $\Sigma$ is confluent then, by definition, all of its branchings are confluent: in particular, its inclusion and regular critical branchings and the normal instances of its right-indexed or left-indexed ones.

Conversely, let us assume that $\Sigma$ is a terminating right-indexed $3$-polygraph (the left-indexed case is similar) such that all of its inclusion and regular critical branchings and all of the normal instances of its right-indexed critical branchings are confluent. It is sufficient to prove  that every non-normal instance of its right-indexed critical branchings is confluent. 

Let us consider a right-indexed critical branching $(A[k],B[k])_k$, which has the following shape, by definition:
$$
\xymatrix@R=-1em{
& {\raisebox{-1.25mm}{\input{pcind-tgt1-k.pstex_t}}}
\\
{\raisebox{-1.25mm}{\input{pcind-src-k.pstex_t}}}
	\ar@3 [ur] ^-{A[k]}
	\ar@3[dr] _-{B[k]}
\\
& {\raisebox{-1.25mm}{\input{pcind-tgt2-k.pstex_t}}}
}
$$

\noindent Let $f$ be a $2$-cell such that $(A[f],B[f])$ is a non-normal instance of $(A[k],B[k])_k$. Since $\Sigma$ terminates,~$f$ admits a normal form, say $g$. We denote by $F$ a $3$-cell from $f$ to $g$. Since $g$ is a normal form, the branching $(A[g],B[g])$ is a normal instance of $(A[k],B[k])_k$ so that, by hypothesis, it is confluent: let us denote by $(G,H)$ a confluence for this branching, with target $h$. With all those ingredients, one builds the following confluence diagram for the critical branching $(A[f],B[f])$, thus concluding the proof:
$$
\xymatrix{
& {\raisebox{-1.25mm}{\input{pcind-tgt1-f.pstex_t}}}
	\ar@3 [rr] ^-{D[F]}
&& {\raisebox{-1.25mm}{\input{pcind-tgt1-g.pstex_t}}}
	\ar@3 [dr] ^-{G}
\\
{\raisebox{-1.25mm}{\input{pcind-src-f.pstex_t}}}
	\ar@3 [ur] ^-{A[f]}
	\ar@3 [dr] _-{B[f]}
	\ar@3 [rr] |-{C[F]}
&& {\raisebox{-1.25mm}{\input{pcind-src-g.pstex_t}}}
	\ar@3 [ur] ^-{A[g]}
	\ar@3 [dr] _-{B[g]}
&& {\raisebox{-1.25mm}{\begin{picture}(0,0)%
\includegraphics{nf-h.pstex}%
\end{picture}%
\setlength{\unitlength}{4144sp}%
\begingroup\makeatletter\ifx\SetFigFont\undefined%
\gdef\SetFigFont#1#2#3#4#5{%
  \reset@font\fontsize{#1}{#2pt}%
  \fontfamily{#3}\fontseries{#4}\fontshape{#5}%
  \selectfont}%
\fi\endgroup%
\begin{picture}(294,294)(-911,-1108)
\put(-764,-1006){\makebox(0,0)[b]{\smash{{\SetFigFont{8}{9.6}{\familydefault}{\mddefault}{\updefault}{\color[rgb]{0,0,0}$h$}%
}}}}
\end{picture}%
}}
\\
& {\raisebox{-1.25mm}{\input{pcind-tgt2-f.pstex_t}}}
	\ar@3 [rr] _-{E[F]}
&& {\raisebox{-1.25mm}{\input{pcind-tgt2-g.pstex_t}}}
	\ar@3 [ur] _-{H}
}
$$

\end{proof}

\subsubsection{Homotopy bases of indexed \pdf{3}-polygraphs}

Let $\Sigma$ be a locally confluent and right-indexed (resp. left-indexed) $3$-polygraph. We assume that a confluence has been chosen for each inclusion and regular critical branching and each normal instance of each right-indexed (resp. left-indexed) critical branching. We denote by $\Gamma_{\Sigma}$ the collection of the $2$-spheres of $\Sigma^*$ corresponding to these confluence diagrams. 

\subsubsection{Proposition}
\label{Theorem:RightIndexedTDF}

\begin{em}
Let $\Sigma$ be a convergent right-indexed (resp. left-indexed) $3$-polygraph. Then $\Gamma_{\Sigma}$ is a homotopy basis of $\Sigma^\top$. 
\end{em}

\begin{proof}
The proof follows the same scheme as the results of~\ref{Subsection:BranchingHomotopyBases}, where it was proved that the family of $3$-spheres associated to the confluence diagrams of all the critical branchings was a homotopy basis.

First, we prove that every local branching of $(A,B)$ of $\Sigma$ admits a confluence $(A',B')$ such that $A\star_2 A' \approx_{\Gamma_{\Sigma}} B\star_2 B'$ holds. The proof is the same as in~\ref{Subsection:BranchingHomotopyBases} when $(A,B)$ is a trivial or when it is generated by an inclusion or a regular critical branching.

There remains to check the cases of local branchings of the shape $C(A[f],B[f])$, where $(A[k],B[k])_k$ is a right-indexed (resp. left-indexed) critical branching and where $C$ is a context. For that, we proceed by Noetherian induction on the indexing $2$-cell $f$, thanks to the termination of $\Sigma$. 

When $f$ is a normal form, then $(A[f],B[f])$ is a normal instance of the branching $(A[k],B[k])_k$. To build $\Gamma_{\Sigma}$ we have fixed a confluence for this branching, say $(A',B')$. Then we have:
$$
C[A[f]] \star_2 A' \:\approx_{\Gamma_{\Sigma}}\: C[B[f]] \star_2 B'.
$$

\noindent Let us assume that $f$ is a $2$-cell which is not a normal form and such that $(A[f],B[f])$ is an instance of the branching $(A[k],B[k])_k$. Moreover, we assume that, for every $2$-cell $g$ such that $f$ reduces into $g$ and $(A[g],B[g])$ is an instance of $(A[k],B[k])_k$, there exists a confluence $(A',B')$ for $(A[g],B[g])$ such that $A[g]\star_2 A' \approx_{\Gamma_{\Sigma}} B[g] \star_2 B'$ holds. 

Since $f$ in not a normal form, we can choose a $2$-cell $g$ such that $f$ reduces into $g$, through a $3$-cell $F$. Since $f$ and $g$ have the same boundary, we have an instance $(A[g],B[g])$ of the branching $(A[k],B[k])_k$. We apply the induction hypothesis to $g$ to get a confluence $(A',B')$, with target denoted by $h$, such that $A[g]\star_2 A' \approx_{\Gamma_{\Sigma}} B[g] \star_2 B'$ holds. Moreover, the branchings $(C[A[f]],C[sA[F]])$ and $(C[B[f]],C[sB[F]])$ are trivial branchings, yielding:
$$
C[A[f]] \star_2 C[tA[F]] \:\approx_{\Gamma_{\Sigma}}\: C[sA[F]] \star_2 C[A[g]]
$$

\noindent and
$$
C[B[f]] \star_2 C[tB[F]] \:\approx_{\Gamma_{\Sigma}}\: C[sB[F]] \star_2 C[B[g]].
$$

\noindent With these constructions, we build the following diagram, where we have assumed that the considered branching was right-indexed -- the case of a left-indexed critical branching is similar:
$$
\xymatrix{
& {\raisebox{-1.25mm}{\input{hb-tgt1-f.pstex_t}}}
	\ar@3 [rr] ^-{ C[tA[F]] }
	\ar@{} [dr] |-{\approx_{\Gamma_{\Sigma}}}
&& {\raisebox{-1.25mm}{\input{hb-tgt1-g.pstex_t}}}
	\ar@3 [dr] ^-{ C[A'] }
\\
{\raisebox{-1.25mm}{\input{hb-src-f.pstex_t}}}
	\ar@3 [ur] ^-{ C[A[f]] }
	\ar@3 [dr] _-{ C[B[f]] }
	\ar@3 [rr] ^-{ C[sA[F]] } _-{ C[sB[F]] }
&& {\raisebox{-1.25mm}{\input{hb-src-g.pstex_t}}}
	\ar@3 [ur] ^-{ C[A[g]] }
	\ar@3 [dr] _-{ C[B[g]] }
	\ar@{} [rr] |-{\approx_{\Gamma_{\Sigma}}}
&& {\raisebox{-1.25mm}{\input{hb-h.pstex_t}}}
\\
& {\raisebox{-1.25mm}{\input{hb-tgt2-f.pstex_t}}}
	\ar@3 [rr] _-{ C[tB[F]] }
	\ar@{} [ur] |-{\approx_{\Gamma_{\Sigma}}}
&& {\raisebox{-1.25mm}{\input{hb-tgt2-g.pstex_t}}}
	\ar@3 [ur] _-{ C[B'] }
}
$$

\noindent One composes the $4$-cells of $\Sigma^\top(\Gamma_{\Sigma})$ of that diagram, to get that $(C[tA[F]]\star_2 C[A'], C[tB[F]] \star_2 C[B'])$ is a confluence that satisfies the required equivalence that concludes the first part of the proof:
$$
C[A[f]] \star_2 C[tA[F]]\star_2 C[A'] \: \approx_{\Gamma_{\Sigma}} \: C[B[f]] \star_2 C[tB[F]] \star_2 C[B'].
$$

\noindent The remainder of the proof is exactly the same as in~\ref{Subsection:BranchingHomotopyBases}. 
\end{proof}

\subsubsection{Theorem}
\label{SubSubSection:MainTheorem1}

\begin{em}
A finite, convergent and finitely indexed $3$-polygraph has finite derivation type.    
\end{em}

\subsection{The \pdf{3}-polygraph of permutations}
\label{Subsection:PolyPerm}

Here we see an example of a $3$-polygraph that is finite, convergent, right-indexed and, thus, with an infinite number of critical branchings, yet with finite derivation type thanks to finite indexation. Another proof for termination and the ideas we use here for proving confluence can be found in~\cite{Lafont03}.

\subsubsection{Definition}

The $3$-polygraph $\Sigma$ has one $0$-cell, one $1$-cell, one $2$-cell~$\figeps{tau}$, and the following two $3$-cells:
$$
\scalefigeps{-2.8mm}{tautau} \overset{\alpha}{\tfl} \scalefigeps{-2.8mm}{id2} 
\qquad\text{and}\qquad 
\scalefigeps{-4.3mm}{yb1} \overset{\beta}{\tfl} \scalefigeps{-4.3mm}{yb2}\:.
$$

\subsubsection{Termination}

We consider the following $\Sigma_2^*$-module $M_{X,\ast,\Zb}$ and derivation $d$ of $\Sigma_2^*$ into $M_{X,\ast,\Zb}$:
$$
X\left(\:\figeps{1-cell}\:\right) \:=\: \Nb,
\qquad
X\left(\figeps{tau}\right) (i,j) \:=\: (j+1,i),
$$
$$
d\left( \figeps{tau} \right) (i,j) \:=\: i.
$$

\noindent The $2$-functor $X$ and the derivation $d$ satisfy the conditions of Theorem~\ref{Theorem:TerminationDerivation}. Indeed, the following required (in)equalities hold:
$$
X\left(\scalefigeps{-2.8mm}{tautau}\right) (i,j) 
	\:=\: (i+1,j+1) 
	\:\geq\: (i,j)
	\:=\: X\left(\scalefigeps{-2.8mm}{id2}\right) (i,j),
$$
$$
X\left(\scalefigeps{-4.3mm}{yb1}\right) (i,j,k) 
	\:=\: (k+2,j+1,i)
	\:=\: X\left(\scalefigeps{-4.3mm}{yb2}\right) (i,j,k), 
$$
$$
d\left(\scalefigeps{-2.8mm}{tautau}\right) (i,j) 
	\:=\: i+j+1
	\:>\: 0
	\:=\: d\left(\scalefigeps{-2.8mm}{id2}\right) (i,j),
$$
$$
d\left(\scalefigeps{-4.3mm}{yb1}\right) (i,j,k) 
	\:=\: 2i+j+1
	\:>\: 2i+j
	\:=\: d\left(\scalefigeps{-4.3mm}{yb2}\right) (i,j,k).
$$

\subsubsection{Normal forms}

First, we note that, if $f$ is a $2$-cell of $\Sigma^*$ such that $d(f)(0,\dots,0)=0$ holds, then~$f$ is a normal form. Otherwise, there exists a context $C$ and a $2$-cell $g$ such that $f=C[g]$ holds and~$g$ is the source of one of the two $3$-cells of $\Sigma$. As a consequence, there exists a family $(i_1,\dots,i_n)$ of natural numbers, with $n=2$ or $n=3$, such that the following inequalities hold:
$$
d(f)(0,\dots,0) \:\geq\: d(g)(i_1,\dots,i_n) \:\geq\: 1.
$$

\noindent Now, let us define $N_0$ as the set of $2$-cells given by the following inductive construction:
$$
\figeps{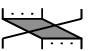} 
	\quad=\quad \figeps{tau} 
	\quad\text{or}\quad \scalefigeps{-2.8mm}{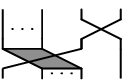} \:.
$$

\noindent We check that the relation
$$
X\left(\figeps{taun1}\right)(i_1,\dots,i_n,j) \:=\: (j+n,i_1,\dots,i_n).
$$

\noindent is satisfied. We proceed by structural induction, using the definition and the functoriality of $X$, to get
$$
X\left(\figeps{tau}\right)(i,j) \:=\: (j+1,i) 
$$

\noindent and
\begin{align*}
X\left(\scalefigeps{-2.8mm}{taun1-tau}\right)(i_1,\dots,i_n,i_{n+1},j) 
\:&=\: \left( X\left(\figeps{taun1}\right) \times \id_{\Nb}
\right) (i_1,\dots,i_n,j+1,i_{n+1})\\
\:&=\: (j+n+1,i_1,\dots,i_{n+1}).
\end{align*}

\noindent Then, we prove that the $2$-cells of $N_0$ are normal forms, still by structural induction. For the base case, we have, by definition of $d$:
$$
d\left(\figeps{tau}\right) (0,0) \:=\: 0.
$$

\noindent For the inductive case, we have, using the fact that $d$ is a derivation:
$$
d\left(\scalefigeps{-2.8mm}{taun1-tau}\right)(0,\dots,0) \:=\: d\left(\figeps{taun1}\right)(0,\dots,0) + d\left(\figeps{tau}\right)(0,0) \:=\: 0.
$$

\noindent Finally, let us denote by $N$ the set of $2$-cells of $\Sigma^*$ given by the following inductive graphical scheme: 
$$
\scalefigeps{-2mm}{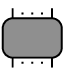} 
	\quad = \quad \ast 
	\quad \text{or} \quad \scalefigeps{-2mm}{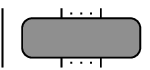}
	\quad \text{or} \quad \scalefigeps{-4mm}{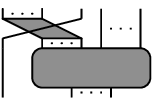} \:.
$$

\noindent We prove that the $2$-cells of $N$ are normal forms, by structural induction. We have $d(\ast) \:=\: 0$,  
$$
d\left(\scalefigeps{-2mm}{1-normal}\right)(i_1,\dots,i_n+1) 
\:=\: d\left(\:\figeps{1-cell}\:\right)(i_1) 
+ d\left(\scalefigeps{-2mm}{normal}\right)(i_2,\dots,i_n) 
\:=\: 0
$$

\noindent and, using the values of $X$ on $N_0$, 
\begin{align*}
& d\left(\scalefigeps{-4mm}{normal-taun1}\right)(i_1,\dots,i_m,j,k_1,\dots,k_n)  \\
\:=\:& d\left(\figeps{taun1}\right)(i_1,\dots,i_m,j) 
+ d\left(\scalefigeps{-2mm}{normal}\right)(i_1,\dots,i_m,k_1,\dots,k_n) 
\:=\: 0.
\end{align*}

\noindent Conversely, let us prove that every normal form of $\Sigma$ is in $N$. We proceed by induction on the pair $(m,n)$ of natural numbers, where $m$ is the size of the $2$-cells and $n$ is the size of their source.

The $2$-cells of $\Sigma^*$ with size $0$ are the $1_n$, where $n$ denotes the $1$-cell with size $n$. All of them are normal forms. Moreover, they belong to $N$: $1_0$ is $\ast$ and, for every natural number $n$, $1_{n+1}=1_1\star_0 1_n$. Moreover, the only $2$-cell of $\Sigma^*$ whose source has size $0$ is $1_0=\ast$, which is a normal form and belongs to~$N$.

Then, let us fix two non-zero natural numbers $m$ and $n$. We assume that, every normal form~$g$ of $\Sigma$ and such that $(\norm{g},\abs{s g})<(m,n)$ holds is in $N$, where we compare pairs of natural numbers with the product order.

Let us consider a normal form $f$ of $\Sigma$, with size $m$ and whose source has size $n$. Since $\norm{f}=m\geq 1$ and since $\figeps{tau}$ is the only $2$-cell of $\Sigma_{\text{Perm}}$, there exists a $2$-cell $g$ such that $f$ decomposes into:
$$
\raisebox{-2mm}{\begin{picture}(0,0)%
\includegraphics{nf-f.pstex}%
\end{picture}%
\setlength{\unitlength}{4144sp}%
\begingroup\makeatletter\ifx\SetFigFont\undefined%
\gdef\SetFigFont#1#2#3#4#5{%
  \reset@font\fontsize{#1}{#2pt}%
  \fontfamily{#3}\fontseries{#4}\fontshape{#5}%
  \selectfont}%
\fi\endgroup%
\begin{picture}(294,294)(-911,-1108)
\put(-764,-1006){\makebox(0,0)[b]{\smash{{\SetFigFont{8}{9.6}{\familydefault}{\mddefault}{\updefault}{\color[rgb]{0,0,0}$f$}%
}}}}
\end{picture}%
} \:=\: \raisebox{-4mm}{\begin{picture}(0,0)%
\includegraphics{nf-tau-g.pstex}%
\end{picture}%
\setlength{\unitlength}{4144sp}%
\begingroup\makeatletter\ifx\SetFigFont\undefined%
\gdef\SetFigFont#1#2#3#4#5{%
  \reset@font\fontsize{#1}{#2pt}%
  \fontfamily{#3}\fontseries{#4}\fontshape{#5}%
  \selectfont}%
\fi\endgroup%
\begin{picture}(654,429)(-866,-568)
\put(-539,-466){\makebox(0,0)[b]{\smash{{\SetFigFont{8}{9.6}{\familydefault}{\mddefault}{\updefault}{\color[rgb]{0,0,0}$g$}%
}}}}
\end{picture}%
} \:.
$$

\noindent Since $f$ is a normal form, then so does $g$. Moreover, $g$ has size $m-1$ and its source has size $n$. We apply the induction hypothesis to $g$: this $2$-cell is in $N$. Its source is $n\geq 1$, so that $g\neq\ast$; there remains two possibilities, by definition of $N$: 
$$
\raisebox{-2mm}{\begin{picture}(0,0)%
\includegraphics{nf-g.pstex}%
\end{picture}%
\setlength{\unitlength}{4144sp}%
\begingroup\makeatletter\ifx\SetFigFont\undefined%
\gdef\SetFigFont#1#2#3#4#5{%
  \reset@font\fontsize{#1}{#2pt}%
  \fontfamily{#3}\fontseries{#4}\fontshape{#5}%
  \selectfont}%
\fi\endgroup%
\begin{picture}(294,294)(-911,-1108)
\put(-764,-1006){\makebox(0,0)[b]{\smash{{\SetFigFont{8}{9.6}{\familydefault}{\mddefault}{\updefault}{\color[rgb]{0,0,0}$g$}%
}}}}
\end{picture}%
} \:=\: \raisebox{-2mm}{\begin{picture}(0,0)%
\includegraphics{nf-1-h.pstex}%
\end{picture}%
\setlength{\unitlength}{4144sp}%
\begingroup\makeatletter\ifx\SetFigFont\undefined%
\gdef\SetFigFont#1#2#3#4#5{%
  \reset@font\fontsize{#1}{#2pt}%
  \fontfamily{#3}\fontseries{#4}\fontshape{#5}%
  \selectfont}%
\fi\endgroup%
\begin{picture}(384,294)(-821,-883)
\put(-584,-781){\makebox(0,0)[b]{\smash{{\SetFigFont{8}{9.6}{\familydefault}{\mddefault}{\updefault}{\color[rgb]{0,0,0}$h$}%
}}}}
\end{picture}%
} 
\qquad\text{or}\qquad
\raisebox{-2mm}{\begin{picture}(0,0)%
\includegraphics{nf-g.pstex}%
\end{picture}%
\setlength{\unitlength}{4144sp}%
\begingroup\makeatletter\ifx\SetFigFont\undefined%
\gdef\SetFigFont#1#2#3#4#5{%
  \reset@font\fontsize{#1}{#2pt}%
  \fontfamily{#3}\fontseries{#4}\fontshape{#5}%
  \selectfont}%
\fi\endgroup%
\begin{picture}(294,294)(-911,-1108)
\put(-764,-1006){\makebox(0,0)[b]{\smash{{\SetFigFont{8}{9.6}{\familydefault}{\mddefault}{\updefault}{\color[rgb]{0,0,0}$g$}%
}}}}
\end{picture}%
} \:=\: \raisebox{-4mm}{\begin{picture}(0,0)%
\includegraphics{nf-taun1-h.pstex}%
\end{picture}%
\setlength{\unitlength}{4144sp}%
\begingroup\makeatletter\ifx\SetFigFont\undefined%
\gdef\SetFigFont#1#2#3#4#5{%
  \reset@font\fontsize{#1}{#2pt}%
  \fontfamily{#3}\fontseries{#4}\fontshape{#5}%
  \selectfont}%
\fi\endgroup%
\begin{picture}(609,429)(-821,-568)
\put(-449,-466){\makebox(0,0)[b]{\smash{{\SetFigFont{8}{9.6}{\familydefault}{\mddefault}{\updefault}{\color[rgb]{0,0,0}$h$}%
}}}}
\end{picture}%
} \:.
$$

\noindent In the first case, the $2$-cell $h$ is a normal form, has size $m-1$ and its source has size $n-1$. By induction hypothesis, we know that $h$ is in $N$. There are two subcases for the decomposition of $f$:
$$
\raisebox{-2mm}{\begin{picture}(0,0)%
\includegraphics{nf-f.pstex}%
\end{picture}%
\setlength{\unitlength}{4144sp}%
\begingroup\makeatletter\ifx\SetFigFont\undefined%
\gdef\SetFigFont#1#2#3#4#5{%
  \reset@font\fontsize{#1}{#2pt}%
  \fontfamily{#3}\fontseries{#4}\fontshape{#5}%
  \selectfont}%
\fi\endgroup%
\begin{picture}(294,294)(-911,-1108)
\put(-764,-1006){\makebox(0,0)[b]{\smash{{\SetFigFont{8}{9.6}{\familydefault}{\mddefault}{\updefault}{\color[rgb]{0,0,0}$f$}%
}}}}
\end{picture}%
} \:=\: \raisebox{-4mm}{\begin{picture}(0,0)%
\includegraphics{nf-tau-1-h.pstex}%
\end{picture}%
\setlength{\unitlength}{4144sp}%
\begingroup\makeatletter\ifx\SetFigFont\undefined%
\gdef\SetFigFont#1#2#3#4#5{%
  \reset@font\fontsize{#1}{#2pt}%
  \fontfamily{#3}\fontseries{#4}\fontshape{#5}%
  \selectfont}%
\fi\endgroup%
\begin{picture}(429,429)(-821,-568)
\put(-584,-466){\makebox(0,0)[b]{\smash{{\SetFigFont{8}{9.6}{\familydefault}{\mddefault}{\updefault}{\color[rgb]{0,0,0}$h$}%
}}}}
\end{picture}%
}
\qquad\text{or}\qquad
\raisebox{-2mm}{\begin{picture}(0,0)%
\includegraphics{nf-f.pstex}%
\end{picture}%
\setlength{\unitlength}{4144sp}%
\begingroup\makeatletter\ifx\SetFigFont\undefined%
\gdef\SetFigFont#1#2#3#4#5{%
  \reset@font\fontsize{#1}{#2pt}%
  \fontfamily{#3}\fontseries{#4}\fontshape{#5}%
  \selectfont}%
\fi\endgroup%
\begin{picture}(294,294)(-911,-1108)
\put(-764,-1006){\makebox(0,0)[b]{\smash{{\SetFigFont{8}{9.6}{\familydefault}{\mddefault}{\updefault}{\color[rgb]{0,0,0}$f$}%
}}}}
\end{picture}%
} \:=\: \raisebox{-4mm}{\begin{picture}(0,0)%
\includegraphics{nf-1-tau-h.pstex}%
\end{picture}%
\setlength{\unitlength}{4144sp}%
\begingroup\makeatletter\ifx\SetFigFont\undefined%
\gdef\SetFigFont#1#2#3#4#5{%
  \reset@font\fontsize{#1}{#2pt}%
  \fontfamily{#3}\fontseries{#4}\fontshape{#5}%
  \selectfont}%
\fi\endgroup%
\begin{picture}(744,429)(-821,-568)
\put(-404,-466){\makebox(0,0)[b]{\smash{{\SetFigFont{8}{9.6}{\familydefault}{\mddefault}{\updefault}{\color[rgb]{0,0,0}$h$}%
}}}}
\end{picture}%
} \:.
$$

\noindent The first decomposition is a proof that $f$ is in $N$, since $h$ is in $N$ and $\figeps{tau}$ is in $N_0$. The second decomposition tells us that $f=\figeps{1-cell}\star_0 f'$, where $f'$ is a normal form (otherwise $f$ would not), has size $m$ and its source has size $n-1$; we apply the induction hypothesis to get that $f'$ is in $N$; then we get that $f$ is in $N$.

Let us examine the second case: the $2$-cell $h$ is a normal form, has size at most $m-2$ and its source has size $n-1$; hence, by induction hypothesis, $h$ is in $N$. There are three subpossibilities:
$$
\raisebox{-2mm}{\begin{picture}(0,0)%
\includegraphics{nf-f.pstex}%
\end{picture}%
\setlength{\unitlength}{4144sp}%
\begingroup\makeatletter\ifx\SetFigFont\undefined%
\gdef\SetFigFont#1#2#3#4#5{%
  \reset@font\fontsize{#1}{#2pt}%
  \fontfamily{#3}\fontseries{#4}\fontshape{#5}%
  \selectfont}%
\fi\endgroup%
\begin{picture}(294,294)(-911,-1108)
\put(-764,-1006){\makebox(0,0)[b]{\smash{{\SetFigFont{8}{9.6}{\familydefault}{\mddefault}{\updefault}{\color[rgb]{0,0,0}$f$}%
}}}}
\end{picture}%
} \:=\: \raisebox{-9mm}{\begin{picture}(0,0)%
\includegraphics{nf-2a.pstex}%
\end{picture}%
\setlength{\unitlength}{4144sp}%
\begingroup\makeatletter\ifx\SetFigFont\undefined%
\gdef\SetFigFont#1#2#3#4#5{%
  \reset@font\fontsize{#1}{#2pt}%
  \fontfamily{#3}\fontseries{#4}\fontshape{#5}%
  \selectfont}%
\fi\endgroup%
\begin{picture}(1374,879)(124,-208)
\put(811,-106){\makebox(0,0)[b]{\smash{{\SetFigFont{8}{9.6}{\familydefault}{\mddefault}{\updefault}{\color[rgb]{0,0,0}$h$}%
}}}}
\end{picture}%
}
\qquad\text{or}\qquad
\raisebox{-2mm}{\begin{picture}(0,0)%
\includegraphics{nf-f.pstex}%
\end{picture}%
\setlength{\unitlength}{4144sp}%
\begingroup\makeatletter\ifx\SetFigFont\undefined%
\gdef\SetFigFont#1#2#3#4#5{%
  \reset@font\fontsize{#1}{#2pt}%
  \fontfamily{#3}\fontseries{#4}\fontshape{#5}%
  \selectfont}%
\fi\endgroup%
\begin{picture}(294,294)(-911,-1108)
\put(-764,-1006){\makebox(0,0)[b]{\smash{{\SetFigFont{8}{9.6}{\familydefault}{\mddefault}{\updefault}{\color[rgb]{0,0,0}$f$}%
}}}}
\end{picture}%
} \:=\: \raisebox{-6mm}{\begin{picture}(0,0)%
\includegraphics{nf-2b.pstex}%
\end{picture}%
\setlength{\unitlength}{4144sp}%
\begingroup\makeatletter\ifx\SetFigFont\undefined%
\gdef\SetFigFont#1#2#3#4#5{%
  \reset@font\fontsize{#1}{#2pt}%
  \fontfamily{#3}\fontseries{#4}\fontshape{#5}%
  \selectfont}%
\fi\endgroup%
\begin{picture}(834,564)(34,242)
\put(451,344){\makebox(0,0)[b]{\smash{{\SetFigFont{8}{9.6}{\familydefault}{\mddefault}{\updefault}{\color[rgb]{0,0,0}$h$}%
}}}}
\end{picture}%
}
\qquad\text{or}\qquad
\raisebox{-2mm}{\begin{picture}(0,0)%
\includegraphics{nf-f.pstex}%
\end{picture}%
\setlength{\unitlength}{4144sp}%
\begingroup\makeatletter\ifx\SetFigFont\undefined%
\gdef\SetFigFont#1#2#3#4#5{%
  \reset@font\fontsize{#1}{#2pt}%
  \fontfamily{#3}\fontseries{#4}\fontshape{#5}%
  \selectfont}%
\fi\endgroup%
\begin{picture}(294,294)(-911,-1108)
\put(-764,-1006){\makebox(0,0)[b]{\smash{{\SetFigFont{8}{9.6}{\familydefault}{\mddefault}{\updefault}{\color[rgb]{0,0,0}$f$}%
}}}}
\end{picture}%
} \:=\: \raisebox{-4mm}{\begin{picture}(0,0)%
\includegraphics{nf-2c.pstex}%
\end{picture}%
\setlength{\unitlength}{4144sp}%
\begingroup\makeatletter\ifx\SetFigFont\undefined%
\gdef\SetFigFont#1#2#3#4#5{%
  \reset@font\fontsize{#1}{#2pt}%
  \fontfamily{#3}\fontseries{#4}\fontshape{#5}%
  \selectfont}%
\fi\endgroup%
\begin{picture}(1104,429)(-146,62)
\put(406,164){\makebox(0,0)[b]{\smash{{\SetFigFont{8}{9.6}{\familydefault}{\mddefault}{\updefault}{\color[rgb]{0,0,0}$h$}%
}}}}
\end{picture}%
} \:.
$$

\noindent The first subcase is, in fact, impossible since $f$ would contain the source of a $3$-cell, which contradicts the assumption that $f$ is a normal form. The second case gives that $f$ is in $N$. In the third case, we have a decomposition of $f$ into $(f'\star_0 1_p)\star_1(1_1\star_0 f'')$ where $f'$ is in $N_0$ and $f''$ is a normal form (otherwise $f$ would not), has size at most $m-1$ and has source $n-1$: thus, we apply the induction hypothesis to get that $f''$ and, hence, $f$ are in $N$.

\subsubsection{Confluence}
\label{Subsubsection:ConfluencePerm}

The $3$-polygraph $\Sigma$ has three regular and one right-indexed critical branchings, with the following sources:
$$
\scalefigeps{-4mm}{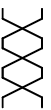} \:\:, 
\qquad\qquad
\scalefigeps{-5.5mm}{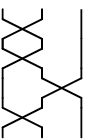} \:\:,
\qquad\qquad
\scalefigeps{-5.5mm}{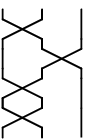} \:\:, 
\qquad\qquad
\raisebox{-8mm}{\begin{picture}(0,0)%
\includegraphics{paire-critique-yb.pstex}%
\end{picture}%
\setlength{\unitlength}{4144sp}%
\begingroup\makeatletter\ifx\SetFigFont\undefined%
\gdef\SetFigFont#1#2#3#4#5{%
  \reset@font\fontsize{#1}{#2pt}%
  \fontfamily{#3}\fontseries{#4}\fontshape{#5}%
  \selectfont}%
\fi\endgroup%
\begin{picture}(834,834)(79,-73)
\put(631,299){\makebox(0,0)[b]{\smash{{\SetFigFont{8}{9.6}{\familydefault}{\mddefault}{\updefault}{\color[rgb]{0,0,0}$k$}%
}}}}
\end{picture}%
} \:\:.
$$

\noindent From Theorem~\ref{Theorem:RightIndexedConfluence}, we know that, to get confluence of $\Sigma$, it is sufficient to prove that the three regular critical branchings are confluent and that each normal instance of the right-indexed one is. First, we check that the three regular critical branchings are confluent:
$$
\xymatrix{
{\figeps{tautautau}}
	\ar@3 @/^6ex/ [rr] ^-{\alpha} _-{\hole}="1"
	\ar@3 @/_6ex/ [rr] _-{\alpha} ^-{\hole}="2"
	\ar@4 "1";"2" |-*+[o]{\scriptstyle \alpha\alpha}
&& {\figeps{tau}}
}
$$
$$
\xymatrix{
{\figeps{tau-yb1}}
	\ar@3 [rrr] ^-{\alpha} _-{\hole}="1"
	\ar@3 [dr] _-{\beta}
&&& {\figeps{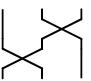}}
\\
& {\figeps{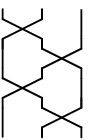}}
	\ar@3 [r] _-{\beta} ^-{\hole}="2"
	\ar@4"1";"2" |-*+[o]{\scriptstyle \alpha\beta}
& {\figeps{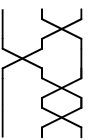}}
	\ar@3 [ur] _-{\alpha}
}
\qquad
\xymatrix{
& {\figeps{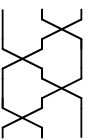}}
	\ar@3 [r] ^-{\beta} _-{\hole}="3"
& {\figeps{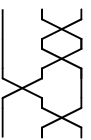}}
	\ar@3 [dr] _-{\alpha}
\\
{\figeps{yb1-tau}}
	\ar@3 [ur] ^-{\beta}
	\ar@3 [rrr] _-{\alpha} ^-{\hole}="4"
	\ar@4 "3";"4" |-*+[o]{\scriptstyle \beta\alpha}
&&& {\figeps{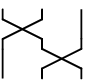}}
}
$$

\noindent From the inductive characterization of the set $N$ of normal forms we have given, we deduce that there are two normal instances of the right-indexed critical branching: for $k=\figeps{1-cell}$ and $k=\figeps{tau}$. We check that both are confluent. For $k=\figeps{1-cell}\;$, we have:
$$
\xymatrix{
& {\figeps{1X-X1-1X-1X-X1}}
	\ar@3 [r] ^-{\alpha} _-{\hole}="1"
& {\figeps{1X-X1-X1}}
	\ar@3 [dr] ^-{\alpha}
\\
{\figeps{yb-yb-1} }
	\ar@3 [ur] ^-{\beta}
	\ar@3 [dr] _-{\beta}
&&& {\figeps{1X}}
\\
& {\figeps{X1-1X-1X-X1-1X}}
	\ar@3 [r] _-{\alpha} ^-{\hole}="2"
& {\figeps{X1-X1-1X}}
	\ar@3 [ur] _-{\alpha} 
\ar@4 "1";"2" |-*+[o]{\scriptstyle \beta\beta\left(\:\raisebox{-1mm}{\smallfigeps{1-cell}}\:\right)}
}
$$

\noindent And, for $k=\figeps{tau}$, we have:
$$
\xymatrix{
& {\figeps{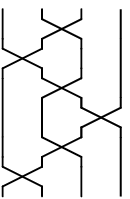}}
	\ar@3 [r] ^-{\beta} 
& {\figeps{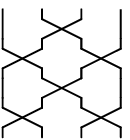}}
	\ar@3 [r] ^-{\beta} 
& {\figeps{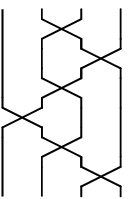}}
	\ar@3 [dr] ^-{\beta}
\\
{\figeps{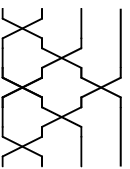}}
	\ar@3 [ur] ^-{\beta}
	\ar@3 [dr] _-{\beta}
&&&& {\figeps{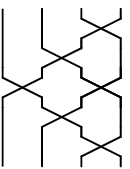}}
\\
& {\figeps{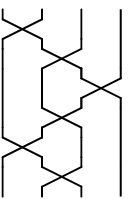}}
	\ar@3 [r] _-{\beta} 
& {\figeps{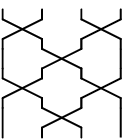}}
	\ar@3 [r] _-{\beta} 
& {\figeps{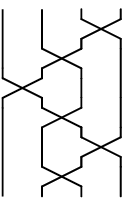}}
	\ar@3 [ur] _-{\beta} 
\ar@4 "1,3";"3,3" |-*+[o]{\scriptstyle \beta\beta\left(\raisebox{-1mm}{\smallfigeps{tau}}\right)}
}
$$

\subsubsection{Theorem}

\begin{em}
The $3$-polygraph $\Sigma$ has finite derivation type.
\end{em}

\begin{proof} 
The $3$-polygraph $\Sigma$ is finite, convergent and finitely indexed. Thus Theorem~\ref{Theorem:RightIndexedTDF} tells us that $\Sigma$ has finite derivation type. More precisely, the five $4$-cells $\alpha\alpha$, $\alpha\beta$, $\beta\alpha$, $\beta\beta\left(\:\figeps{1-cell}\:\right)$ and $\beta\beta\left(\;\figeps{tau}\;\right)$ form a homotopy basis of the track $3$-category $\Sigma^{\top}$.
\end{proof}

\subsection{The main counterexample}
\label{Subsection:main_counter_example}

We prove here that, without finite indexation, finiteness and convergence alone are not sufficient enough to ensure that a $3$-polygraph has finite derivation type.

Let us consider the $3$-polygraph $\Sigma$ with one $0$-cell, one $1$-cell, three
$2$-cells $\figeps{o}$, $\figeps{n}$ and $\figeps{u}$ and the following four $3$-cells: 
$$
\scalefigeps{-2mm}{on} \overset{\alpha}{\tfl} \scalefigeps{-2mm}{no} 
\:, \qquad
\scalefigeps{-2mm}{ou} \overset{\beta}{\tfl} \scalefigeps{-2mm}{uo}
\:, \qquad
\scalefigeps{-2mm}{nu} \overset{\gamma}{\tfl} \scalefigeps{-2mm}{id1} 
\:, \qquad
\scalefigeps{-2mm}{un} \overset{\delta}{\tfl} \scalefigeps{-2mm}{id1} 
\:.
$$

\noindent We define by induction on the natural number $k$ the $2$-cell $\figeps{o}^k$ as follows:
$$
\figeps{o} ^0 \:=\: \figeps{1-cell} 
	\qquad \text{and} \qquad
\figeps{o} ^{k+1} \:=\: \figeps{o}^k\star_1\figeps{o}.
$$

\subsubsection{Termination}

To prove that the $3$-polygraph $\Sigma$ terminates, we proceed in two steps. First, we consider the derivation $\norm{\cdot}_{\smallfigeps{n}}$, into the trivial module $M_{\ast,\ast,\Zb}$. It satisfies the equalities
$$
\norm{s\alpha}_{\smallfigeps{n}} \:=\: 1 \:=\: \norm{t\alpha}_{\smallfigeps{n}}
\qquad\text{and}\qquad
\norm{s\beta}_{\smallfigeps{n}} \:=\: 0 \:=\: \norm{t\beta}_{\smallfigeps{n}}
$$

\noindent and the strict inequalities
$$
\norm{s\gamma}_{\smallfigeps{n}} \:=\: 1 \:>\: 0 \:=\: \norm{t\gamma}_{\smallfigeps{n}}
\qquad\text{and}\qquad
\norm{s\delta}_{\smallfigeps{n}} \:=\: 1 \:>\: 0 \:=\: \norm{t\delta}_{\smallfigeps{n}}.
$$

\noindent As a consequence, one gets that, if the $3$-polygraph $\Sigma' = (\Sigma_2,\ens{\alpha,\beta})$ terminates, then so does the $3$-polygraph $\Sigma$. Indeed, otherwise, there would exist an infinite reduction sequence $(f_n)_{n\in\Nb}$ in $\Sigma$ and, thus, an infinite decreasing sequence $(\norm{f_n}_{\smallfigeps{n}})_{n\in\Nb}$ of natural numbers; moreover, this last sequence would be strictly decreasing at each step $n$ that is generated by either $\gamma$ or $\delta$. Thus, after some natural number~$p$, this sequence could be generated by $\alpha$ and $\beta$ only. This would yield an infinite reduction sequence $(f_n)_{n\geq p}$ in $\Sigma'$, which is impossible by hypothesis. Let us note that one could have used the derivation~$\norm{\cdot}_{\smallfigeps{u}}$ with the same results.

To prove that $\Sigma'$ terminates, we consider the derivation $d$ into the $\Sigma_2^*$-module $M_{X,Y,\Zb}$ given by: 
$$
X\left(\:\figeps{1-cell}\:\right) \:=\: \Nb, \qquad
X\left(\figeps{n}\right) \:=\: (0,0), \qquad
X\left(\figeps{o}\right)(i) \:=\: i+1,
$$
$$
Y\left(\:\figeps{1-cell}\:\right) \:=\: \Nb, \qquad
Y\left(\figeps{u}\right) \:=\: (0,0), \qquad
Y\left(\figeps{o}\right)(i) \:=\: i+1,
$$
$$
d\left(\figeps{n}\right)(i,j) \:=\: i, \qquad 
d\left(\figeps{u}\right)(i,j) \:=\: i, \qquad 
d\left(\figeps{o}\right)(i,j) \:=\: 0. 
$$

\noindent Since $d$ is a derivation, one gets:
\begin{align*}
d(\alpha)
\:=&\: 
	d\left(\scalefigeps{-2mm}{on}\right)
	- d\left(\scalefigeps{-2mm}{no}\right) \\
\:=&\: 
	d\left(\figeps{n}\right) \star_1 \left(\figeps{o}\:\figeps{1-cell}\:\right) 
	\:+\: \figeps{n} \star_1 \left( d\left(\figeps{o}\right) \star_0 \:\figeps{1-cell}\: \right) 
	\:-\: d\left(\figeps{n}\right) \star_1 \left(\:\figeps{1-cell}\:\figeps{o}\right) 
	\:-\: \figeps{n} \star_1 \left( \:\figeps{1-cell}\: \star_0 d\left(\figeps{o}\right) \right) \:.
\end{align*}

\noindent Thus, for every natural numbers $i$ and $j$, one gets:
\begin{align*}
d(\alpha) (i,j) 
\:=&\:
	d\left(\figeps{n}\right) (i+1,j) 
	+ d\left(\figeps{o}\right) (0,i)
	- d\left(\figeps{n}\right) (i,j+1) 
	- d\left(\figeps{o}\right) (0,j) \\
\:=&\: (i+1) \:+\: 0 \:-\: i \:-\: 0 \\
\:=&\: 1.
\end{align*}

\noindent Similarly, one gets $d(\beta)(i,j)=1$ for every natural numbers $i$ and $j$, yielding, thanks to Theorem~\ref{Theorem:TerminationDerivation}, the termination of $\Sigma'$ and, thus, of $\Sigma$.

\subsubsection{Normal forms}

Let $f$ be a $2$-cell of $\Sigma^*$, that cannot be reduced by the $3$-cells $\gamma$ and $\delta$ and which satisfies:
$$
d(f)(0,\dots,0) \:=\: 0.
$$

\noindent Then $f$ is a normal form. Indeed, otherwise, there exists a context $C$ such that $f=C[g]$, with either $g=s\alpha$ or $g=s\beta$. As a consequence, there exist two natural numbers $i$ and $j$ such that the following inequalities hold:
$$
d(f)(0,\dots,0) \geq d(g)(i,j) \geq 1.
$$

\noindent Now, we define $N$ as the set of $2$-cells given by the following inductive construction scheme:
$$
\begin{array}{l c l c l c l}
\scalefigeps{-2mm}{normal} 
& \qquad = \qquad 
	& \text{(a)} \quad \ast
& \qquad \text{or} \qquad 
	& \text{(b)} \quad \raisebox{-3mm}{\begin{picture}(0,0)%
\includegraphics{C-k-f0.pstex}%
\end{picture}%
\setlength{\unitlength}{4144sp}%
\begingroup\makeatletter\ifx\SetFigFont\undefined%
\gdef\SetFigFont#1#2#3#4#5{%
  \reset@font\fontsize{#1}{#2pt}%
  \fontfamily{#3}\fontseries{#4}\fontshape{#5}%
  \selectfont}%
\fi\endgroup%
\begin{picture}(577,378)(-908,290)
\put(-516,524){\makebox(0,0)[lb]{\smash{{\SetFigFont{8}{9.6}{\familydefault}{\mddefault}{\updefault}{\color[rgb]{0,0,0}$k$}%
}}}}
\end{picture}%
} \scalefigeps{-2mm}{normal} 
& \qquad \text{or} \qquad
	& \text{(c)} \quad \raisebox{-2mm}{\begin{picture}(0,0)%
\includegraphics{n-k-nn.pstex}%
\end{picture}%
\setlength{\unitlength}{4144sp}%
\begingroup\makeatletter\ifx\SetFigFont\undefined%
\gdef\SetFigFont#1#2#3#4#5{%
  \reset@font\fontsize{#1}{#2pt}%
  \fontfamily{#3}\fontseries{#4}\fontshape{#5}%
  \selectfont}%
\fi\endgroup%
\begin{picture}(669,291)(-686,-73)
\put(-202, 74){\makebox(0,0)[lb]{\smash{{\SetFigFont{8}{9.6}{\familydefault}{\mddefault}{\updefault}{\color[rgb]{0,0,0}$k$}%
}}}}
\end{picture}%
} \scalefigeps{-2mm}{normal} 
\\
&&& \text{or} 
	& \text{(d)} \quad \raisebox{-2mm}{\begin{picture}(0,0)%
\includegraphics{u-k-uu.pstex}%
\end{picture}%
\setlength{\unitlength}{4144sp}%
\begingroup\makeatletter\ifx\SetFigFont\undefined%
\gdef\SetFigFont#1#2#3#4#5{%
  \reset@font\fontsize{#1}{#2pt}%
  \fontfamily{#3}\fontseries{#4}\fontshape{#5}%
  \selectfont}%
\fi\endgroup%
\begin{picture}(669,318)(-686,290)
\put(-202,524){\makebox(0,0)[lb]{\smash{{\SetFigFont{8}{9.6}{\familydefault}{\mddefault}{\updefault}{\color[rgb]{0,0,0}$k$}%
}}}}
\end{picture}%
} \scalefigeps{-2mm}{normal} 
& \text{or} 
& \text{(e)} \quad \figeps{o}^k \: \scalefigeps{-2mm}{normal} \:.
\end{array}
$$

\noindent We use the special graphical representations $\figeps{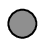}$, $\figeps{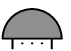}$ and $\figeps{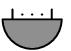}$ for $2$-cells of $N$ which have, respectively, degenerate source and target, degenerate source, degenerate target.

We start by checking that the $2$-cells of $N$ are normal forms. For that, one proceeds by structural induction, using the construction scheme, in order to prove two properties.

The first one is that each $2$-cell of $N$ is irreducible by the $3$-cells $\gamma$ and $\delta$: this is an observation that the given construction scheme does not allow any $2$-cell of $N$ to contain either $\figeps{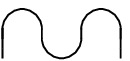}$ or $\figeps{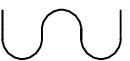}$.

The second property is that, for a $2$-cell $f$ of $N$, one
has $d(f)(0,\dots,0)=0$. For the base case, \ie, when $f$ is
built using construction rule (a), one has $d(\ast)=0$ since
$d$ is a derivation. Then, for the induction, there are four
cases, depending on the construction rule used to build $f$: 
\begin{align*}
\text{(b)} \quad
& d\left( \raisebox{-3mm}{\begin{picture}(0,0)%
\includegraphics{C-k-f0.pstex}%
\end{picture}%
\setlength{\unitlength}{4144sp}%
\begingroup\makeatletter\ifx\SetFigFont\undefined%
\gdef\SetFigFont#1#2#3#4#5{%
  \reset@font\fontsize{#1}{#2pt}%
  \fontfamily{#3}\fontseries{#4}\fontshape{#5}%
  \selectfont}%
\fi\endgroup%
\begin{picture}(577,378)(-908,290)
\put(-516,524){\makebox(0,0)[lb]{\smash{{\SetFigFont{8}{9.6}{\familydefault}{\mddefault}{\updefault}{\color[rgb]{0,0,0}$k$}%
}}}}
\end{picture}%
} \scalefigeps{-2mm}{normal}  \right) (0,\dots,0) \\
=\:& d\left(\figeps{n}\right)(0,k) \:+\:
		d\left(\figeps{u}\right)(0,k) \:+\:
		k\cdot d\left(\figeps{o}\right)(0,0) \:+\:
		d\left(\figeps{f0}\right)  \:+\:
		d\left(\scalefigeps{-2mm}{normal} \right)(0,\dots,0) \\
=\:& 0.
\end{align*}
\begin{align*}
\text{(c)} \quad
& d\left( \raisebox{-2mm}{\begin{picture}(0,0)%
\includegraphics{n-k-nn.pstex}%
\end{picture}%
\setlength{\unitlength}{4144sp}%
\begingroup\makeatletter\ifx\SetFigFont\undefined%
\gdef\SetFigFont#1#2#3#4#5{%
  \reset@font\fontsize{#1}{#2pt}%
  \fontfamily{#3}\fontseries{#4}\fontshape{#5}%
  \selectfont}%
\fi\endgroup%
\begin{picture}(669,291)(-686,-73)
\put(-202, 74){\makebox(0,0)[lb]{\smash{{\SetFigFont{8}{9.6}{\familydefault}{\mddefault}{\updefault}{\color[rgb]{0,0,0}$k$}%
}}}}
\end{picture}%
} \scalefigeps{-2mm}{normal}  \right) (0,\dots,0) \\
=\:& d\left(\figeps{n}\right)(0,k) \:+\:
		k\cdot d\left(\figeps{o}\right)(0,0) \:+\:
		d\left(\figeps{nn}\right)(0,\dots,0)  \:+\:
		d\left(\scalefigeps{-2mm}{normal} \right)(0,\dots,0) \\
=\:& 0.
\end{align*}
\begin{align*}
\text{(d)} \quad
& d\left( \raisebox{-2mm}{\begin{picture}(0,0)%
\includegraphics{u-k-uu.pstex}%
\end{picture}%
\setlength{\unitlength}{4144sp}%
\begingroup\makeatletter\ifx\SetFigFont\undefined%
\gdef\SetFigFont#1#2#3#4#5{%
  \reset@font\fontsize{#1}{#2pt}%
  \fontfamily{#3}\fontseries{#4}\fontshape{#5}%
  \selectfont}%
\fi\endgroup%
\begin{picture}(669,318)(-686,290)
\put(-202,524){\makebox(0,0)[lb]{\smash{{\SetFigFont{8}{9.6}{\familydefault}{\mddefault}{\updefault}{\color[rgb]{0,0,0}$k$}%
}}}}
\end{picture}%
} \scalefigeps{-2mm}{normal}  \right) (0,\dots,0) \\
=\:& d\left(\figeps{u}\right)(0,k) \:+\:
		k\cdot d\left(\figeps{o}\right)(0,0) \:+\:
		d\left(\figeps{uu}\right)(0,\dots,0)  \:+\:
		d\left(\scalefigeps{-2mm}{normal} \right)(0,\dots,0) \\
=\:& 0.
\end{align*}
\begin{align*}
\text{(e)} \quad
& d\left( \figeps{o}^k \scalefigeps{-2mm}{normal}  \right) (0,\dots,0) \\
=\:& k\cdot d\left(\figeps{o}\right)(0,0) \:+\:
		d\left(\scalefigeps{-2mm}{normal} \right)(0,\dots,0) \\
=\:& 0.
\end{align*}

\noindent Now, let us prove that every $2$-cell of $\Sigma^*$ that is a normal form is contained in the set $N$. We proceed by induction on the triple $(m,n,p)$ of natural numbers, where $m$ is the size of the $2$-cells, $n$ the size of their source, $p$ the size of their target.

The only $2$-cells of $\Sigma^*$ with size $0$ are the $1_n$, where $n$ denotes the $1$-cell with size $n$. All of them are normal forms and belong to $N$. Indeed, each $1_n$ can be formed, from $\ast$, by $n$ subsequent applications of the construction rule (e) with $k=0$.

The $2$-cells of $\Sigma^*$ with size $1$ are the $1_p\star_0\phi\star_0 1_q$, where $\phi$ is one of $\figeps{n}$, $\figeps{u}$ and~$\figeps{o}$. Such a $2$-cell is always a normal form and belongs to $N$. Indeed, we have seen that $1_q$ is in $N$. Then we get that $\phi\star_0 1_q$ is in~$N$, by case analysis based on $\phi$.
\begin{itemize}
\item If $\phi$ is $\figeps{n}$, construction rule (c) with $\scalefigeps{-2mm}{normal}=1_q$, $\figeps{nn}=\ast$ and $k=0$.
\item If $\phi=\figeps{u}$, construction rule (d) with $\scalefigeps{-2mm}{normal}=1_q$, $\figeps{uu}=\ast$ and $k=0$.
\item If $\phi=\figeps{o}$, construction rule (e) with $\scalefigeps{-2mm}{normal}=1_q$ and $k=1$.
\end{itemize}

\noindent Finally, $1_p\star_0 \phi\star_0 1_q$ is in $N$, built using construction rule (e), applied $p$ times in sequence with $k=0$ and starting from $\scalefigeps{-2mm}{normal}=\phi\star_0 1_q$.

Now, let us fix a non-zero natural number $m$ and two natural numbers $n$ and $p$. We assume that we have proved the result for each normal form $g$ with size at most $m-1$ or with size $m$ and such that the inequality $(\abs{sg},\abs{tg})<(n,p)$ holds. 

Let us consider a normal form $f$ such that $\norm{f}=m$, $\abs{sf}=n$ and $\abs{tf}=p$ hold. Since~$f$ has size at least $1$, there exists a $2$-cell $g$ such that $f$ decomposes in one of the three following ways:
$$
\raisebox{-2mm}{\begin{picture}(0,0)%
\includegraphics{nf-f.pstex}%
\end{picture}%
\setlength{\unitlength}{4144sp}%
\begingroup\makeatletter\ifx\SetFigFont\undefined%
\gdef\SetFigFont#1#2#3#4#5{%
  \reset@font\fontsize{#1}{#2pt}%
  \fontfamily{#3}\fontseries{#4}\fontshape{#5}%
  \selectfont}%
\fi\endgroup%
\begin{picture}(294,294)(-911,-1108)
\put(-764,-1006){\makebox(0,0)[b]{\smash{{\SetFigFont{8}{9.6}{\familydefault}{\mddefault}{\updefault}{\color[rgb]{0,0,0}$f$}%
}}}}
\end{picture}%
} 
	\qquad=\qquad 
\raisebox{-3mm}{\begin{picture}(0,0)%
\includegraphics{nf-n-g.pstex}%
\end{picture}%
\setlength{\unitlength}{4144sp}%
\begingroup\makeatletter\ifx\SetFigFont\undefined%
\gdef\SetFigFont#1#2#3#4#5{%
  \reset@font\fontsize{#1}{#2pt}%
  \fontfamily{#3}\fontseries{#4}\fontshape{#5}%
  \selectfont}%
\fi\endgroup%
\begin{picture}(654,429)(79,332)
\put(406,457){\makebox(0,0)[b]{\smash{{\SetFigFont{8}{9.6}{\familydefault}{\mddefault}{\updefault}{\color[rgb]{0,0,0}$g$}%
}}}}
\end{picture}%
}
	\qquad\text{or}\qquad
\raisebox{-3mm}{\begin{picture}(0,0)%
\includegraphics{nf-u-g.pstex}%
\end{picture}%
\setlength{\unitlength}{4144sp}%
\begingroup\makeatletter\ifx\SetFigFont\undefined%
\gdef\SetFigFont#1#2#3#4#5{%
  \reset@font\fontsize{#1}{#2pt}%
  \fontfamily{#3}\fontseries{#4}\fontshape{#5}%
  \selectfont}%
\fi\endgroup%
\begin{picture}(654,429)(-866,-568)
\put(-539,-444){\makebox(0,0)[b]{\smash{{\SetFigFont{8}{9.6}{\familydefault}{\mddefault}{\updefault}{\color[rgb]{0,0,0}$g$}%
}}}}
\end{picture}%
}
	\qquad\text{or}\qquad
\raisebox{-3mm}{\begin{picture}(0,0)%
\includegraphics{nf-o-g.pstex}%
\end{picture}%
\setlength{\unitlength}{4144sp}%
\begingroup\makeatletter\ifx\SetFigFont\undefined%
\gdef\SetFigFont#1#2#3#4#5{%
  \reset@font\fontsize{#1}{#2pt}%
  \fontfamily{#3}\fontseries{#4}\fontshape{#5}%
  \selectfont}%
\fi\endgroup%
\begin{picture}(474,429)(-776,-568)
\put(-539,-444){\makebox(0,0)[b]{\smash{{\SetFigFont{8}{9.6}{\familydefault}{\mddefault}{\updefault}{\color[rgb]{0,0,0}$g$}%
}}}}
\end{picture}%
}\:.
$$

\noindent One denotes by $\phi$ the generating $2$-cell corresponding to each of those decompositions: $\figeps{n}$, $\figeps{u}$ and~$\figeps{o}$, respectively. Since $f$ is a normal form, so does $g$ and $g$ has size $m-1$: we apply the induction hypothesis to it, so that we know that $g$ is in $N$. Thus, $g$ decomposes into one of the five following ways, corresponding to the five construction rules of~$N$:
$$
\begin{array}{l c l c l c l}
\raisebox{-2mm}{\begin{picture}(0,0)%
\includegraphics{nf-g.pstex}%
\end{picture}%
\setlength{\unitlength}{4144sp}%
\begingroup\makeatletter\ifx\SetFigFont\undefined%
\gdef\SetFigFont#1#2#3#4#5{%
  \reset@font\fontsize{#1}{#2pt}%
  \fontfamily{#3}\fontseries{#4}\fontshape{#5}%
  \selectfont}%
\fi\endgroup%
\begin{picture}(294,294)(-911,-1108)
\put(-764,-1006){\makebox(0,0)[b]{\smash{{\SetFigFont{8}{9.6}{\familydefault}{\mddefault}{\updefault}{\color[rgb]{0,0,0}$g$}%
}}}}
\end{picture}%
} 
& \qquad = \qquad 
	& \text{(i)} \quad \ast
& \qquad \text{or} \qquad 
	& \text{(ii)} \quad \raisebox{-3mm}{\begin{picture}(0,0)%
\includegraphics{C-k-f0.pstex}%
\end{picture}%
\setlength{\unitlength}{4144sp}%
\begingroup\makeatletter\ifx\SetFigFont\undefined%
\gdef\SetFigFont#1#2#3#4#5{%
  \reset@font\fontsize{#1}{#2pt}%
  \fontfamily{#3}\fontseries{#4}\fontshape{#5}%
  \selectfont}%
\fi\endgroup%
\begin{picture}(577,378)(-908,290)
\put(-516,524){\makebox(0,0)[lb]{\smash{{\SetFigFont{8}{9.6}{\familydefault}{\mddefault}{\updefault}{\color[rgb]{0,0,0}$k$}%
}}}}
\end{picture}%
} \raisebox{-2mm}{\begin{picture}(0,0)%
\includegraphics{nf-h.pstex}%
\end{picture}%
\setlength{\unitlength}{4144sp}%
\begingroup\makeatletter\ifx\SetFigFont\undefined%
\gdef\SetFigFont#1#2#3#4#5{%
  \reset@font\fontsize{#1}{#2pt}%
  \fontfamily{#3}\fontseries{#4}\fontshape{#5}%
  \selectfont}%
\fi\endgroup%
\begin{picture}(294,294)(-911,-1108)
\put(-764,-1006){\makebox(0,0)[b]{\smash{{\SetFigFont{8}{9.6}{\familydefault}{\mddefault}{\updefault}{\color[rgb]{0,0,0}$h$}%
}}}}
\end{picture}%
} 
& \qquad \text{or} \qquad
	& \text{(iii)} \quad \raisebox{-2mm}{\begin{picture}(0,0)%
\includegraphics{n-k-nn.pstex}%
\end{picture}%
\setlength{\unitlength}{4144sp}%
\begingroup\makeatletter\ifx\SetFigFont\undefined%
\gdef\SetFigFont#1#2#3#4#5{%
  \reset@font\fontsize{#1}{#2pt}%
  \fontfamily{#3}\fontseries{#4}\fontshape{#5}%
  \selectfont}%
\fi\endgroup%
\begin{picture}(669,291)(-686,-73)
\put(-202, 74){\makebox(0,0)[lb]{\smash{{\SetFigFont{8}{9.6}{\familydefault}{\mddefault}{\updefault}{\color[rgb]{0,0,0}$k$}%
}}}}
\end{picture}%
} \raisebox{-2mm}{\begin{picture}(0,0)%
\includegraphics{nf-h.pstex}%
\end{picture}%
\setlength{\unitlength}{4144sp}%
\begingroup\makeatletter\ifx\SetFigFont\undefined%
\gdef\SetFigFont#1#2#3#4#5{%
  \reset@font\fontsize{#1}{#2pt}%
  \fontfamily{#3}\fontseries{#4}\fontshape{#5}%
  \selectfont}%
\fi\endgroup%
\begin{picture}(294,294)(-911,-1108)
\put(-764,-1006){\makebox(0,0)[b]{\smash{{\SetFigFont{8}{9.6}{\familydefault}{\mddefault}{\updefault}{\color[rgb]{0,0,0}$h$}%
}}}}
\end{picture}%
} 
\\
&&& \text{or} 
	& \text{(iv)} \quad \raisebox{-2mm}{\begin{picture}(0,0)%
\includegraphics{u-k-uu.pstex}%
\end{picture}%
\setlength{\unitlength}{4144sp}%
\begingroup\makeatletter\ifx\SetFigFont\undefined%
\gdef\SetFigFont#1#2#3#4#5{%
  \reset@font\fontsize{#1}{#2pt}%
  \fontfamily{#3}\fontseries{#4}\fontshape{#5}%
  \selectfont}%
\fi\endgroup%
\begin{picture}(669,318)(-686,290)
\put(-202,524){\makebox(0,0)[lb]{\smash{{\SetFigFont{8}{9.6}{\familydefault}{\mddefault}{\updefault}{\color[rgb]{0,0,0}$k$}%
}}}}
\end{picture}%
} \raisebox{-2mm}{\begin{picture}(0,0)%
\includegraphics{nf-h.pstex}%
\end{picture}%
\setlength{\unitlength}{4144sp}%
\begingroup\makeatletter\ifx\SetFigFont\undefined%
\gdef\SetFigFont#1#2#3#4#5{%
  \reset@font\fontsize{#1}{#2pt}%
  \fontfamily{#3}\fontseries{#4}\fontshape{#5}%
  \selectfont}%
\fi\endgroup%
\begin{picture}(294,294)(-911,-1108)
\put(-764,-1006){\makebox(0,0)[b]{\smash{{\SetFigFont{8}{9.6}{\familydefault}{\mddefault}{\updefault}{\color[rgb]{0,0,0}$h$}%
}}}}
\end{picture}%
} 
& \text{or} 
& \text{(v)} \quad \figeps{o}^k \raisebox{-2mm}{\begin{picture}(0,0)%
\includegraphics{nf-h.pstex}%
\end{picture}%
\setlength{\unitlength}{4144sp}%
\begingroup\makeatletter\ifx\SetFigFont\undefined%
\gdef\SetFigFont#1#2#3#4#5{%
  \reset@font\fontsize{#1}{#2pt}%
  \fontfamily{#3}\fontseries{#4}\fontshape{#5}%
  \selectfont}%
\fi\endgroup%
\begin{picture}(294,294)(-911,-1108)
\put(-764,-1006){\makebox(0,0)[b]{\smash{{\SetFigFont{8}{9.6}{\familydefault}{\mddefault}{\updefault}{\color[rgb]{0,0,0}$h$}%
}}}}
\end{picture}%
} \:.
\end{array}
$$

\noindent We study all the possible decompositions of $f$, depending on the one of $g$ and on $\phi$. In case (i), \ie, when $g=\ast$, we have $\phi=\figeps{u}$, since this is the only possibility to have $t\phi$ degenerate. We have already seen that~$\figeps{u}$ is in $N$. In case (ii), one has the following possibilities, depending on $\phi$:
$$
\raisebox{-2mm}{\begin{picture}(0,0)%
\includegraphics{nf-f.pstex}%
\end{picture}%
\setlength{\unitlength}{4144sp}%
\begingroup\makeatletter\ifx\SetFigFont\undefined%
\gdef\SetFigFont#1#2#3#4#5{%
  \reset@font\fontsize{#1}{#2pt}%
  \fontfamily{#3}\fontseries{#4}\fontshape{#5}%
  \selectfont}%
\fi\endgroup%
\begin{picture}(294,294)(-911,-1108)
\put(-764,-1006){\makebox(0,0)[b]{\smash{{\SetFigFont{8}{9.6}{\familydefault}{\mddefault}{\updefault}{\color[rgb]{0,0,0}$f$}%
}}}}
\end{picture}%
} 
	\qquad=\qquad 
\raisebox{-3mm}{\begin{picture}(0,0)%
\includegraphics{C-k-f0.pstex}%
\end{picture}%
\setlength{\unitlength}{4144sp}%
\begingroup\makeatletter\ifx\SetFigFont\undefined%
\gdef\SetFigFont#1#2#3#4#5{%
  \reset@font\fontsize{#1}{#2pt}%
  \fontfamily{#3}\fontseries{#4}\fontshape{#5}%
  \selectfont}%
\fi\endgroup%
\begin{picture}(577,378)(-908,290)
\put(-516,524){\makebox(0,0)[lb]{\smash{{\SetFigFont{8}{9.6}{\familydefault}{\mddefault}{\updefault}{\color[rgb]{0,0,0}$k$}%
}}}}
\end{picture}%
} \raisebox{-3mm}{\begin{picture}(0,0)%
\includegraphics{nf-n-h.pstex}%
\end{picture}%
\setlength{\unitlength}{4144sp}%
\begingroup\makeatletter\ifx\SetFigFont\undefined%
\gdef\SetFigFont#1#2#3#4#5{%
  \reset@font\fontsize{#1}{#2pt}%
  \fontfamily{#3}\fontseries{#4}\fontshape{#5}%
  \selectfont}%
\fi\endgroup%
\begin{picture}(654,429)(79,332)
\put(406,434){\makebox(0,0)[b]{\smash{{\SetFigFont{8}{9.6}{\familydefault}{\mddefault}{\updefault}{\color[rgb]{0,0,0}$h$}%
}}}}
\end{picture}%
}
	\qquad\text{or}\qquad
\raisebox{-3mm}{\begin{picture}(0,0)%
\includegraphics{C-k-f0.pstex}%
\end{picture}%
\setlength{\unitlength}{4144sp}%
\begingroup\makeatletter\ifx\SetFigFont\undefined%
\gdef\SetFigFont#1#2#3#4#5{%
  \reset@font\fontsize{#1}{#2pt}%
  \fontfamily{#3}\fontseries{#4}\fontshape{#5}%
  \selectfont}%
\fi\endgroup%
\begin{picture}(577,378)(-908,290)
\put(-516,524){\makebox(0,0)[lb]{\smash{{\SetFigFont{8}{9.6}{\familydefault}{\mddefault}{\updefault}{\color[rgb]{0,0,0}$k$}%
}}}}
\end{picture}%
} \raisebox{-3mm}{\begin{picture}(0,0)%
\includegraphics{nf-u-h.pstex}%
\end{picture}%
\setlength{\unitlength}{4144sp}%
\begingroup\makeatletter\ifx\SetFigFont\undefined%
\gdef\SetFigFont#1#2#3#4#5{%
  \reset@font\fontsize{#1}{#2pt}%
  \fontfamily{#3}\fontseries{#4}\fontshape{#5}%
  \selectfont}%
\fi\endgroup%
\begin{picture}(654,429)(-866,-568)
\put(-539,-466){\makebox(0,0)[b]{\smash{{\SetFigFont{8}{9.6}{\familydefault}{\mddefault}{\updefault}{\color[rgb]{0,0,0}$h$}%
}}}}
\end{picture}%
}
	\qquad\text{or}\qquad
\raisebox{-3mm}{\begin{picture}(0,0)%
\includegraphics{C-k-f0.pstex}%
\end{picture}%
\setlength{\unitlength}{4144sp}%
\begingroup\makeatletter\ifx\SetFigFont\undefined%
\gdef\SetFigFont#1#2#3#4#5{%
  \reset@font\fontsize{#1}{#2pt}%
  \fontfamily{#3}\fontseries{#4}\fontshape{#5}%
  \selectfont}%
\fi\endgroup%
\begin{picture}(577,378)(-908,290)
\put(-516,524){\makebox(0,0)[lb]{\smash{{\SetFigFont{8}{9.6}{\familydefault}{\mddefault}{\updefault}{\color[rgb]{0,0,0}$k$}%
}}}}
\end{picture}%
} \raisebox{-3mm}{\begin{picture}(0,0)%
\includegraphics{nf-o-h.pstex}%
\end{picture}%
\setlength{\unitlength}{4144sp}%
\begingroup\makeatletter\ifx\SetFigFont\undefined%
\gdef\SetFigFont#1#2#3#4#5{%
  \reset@font\fontsize{#1}{#2pt}%
  \fontfamily{#3}\fontseries{#4}\fontshape{#5}%
  \selectfont}%
\fi\endgroup%
\begin{picture}(474,429)(-776,-568)
\put(-539,-466){\makebox(0,0)[b]{\smash{{\SetFigFont{8}{9.6}{\familydefault}{\mddefault}{\updefault}{\color[rgb]{0,0,0}$h$}%
}}}}
\end{picture}%
}\:.
$$

\noindent The following $2$-cells must be normal forms, since $f$ is, and they have size at most $m-2$:
$$
\raisebox{-3mm}{\begin{picture}(0,0)%
\includegraphics{nf-n-h.pstex}%
\end{picture}%
\setlength{\unitlength}{4144sp}%
\begingroup\makeatletter\ifx\SetFigFont\undefined%
\gdef\SetFigFont#1#2#3#4#5{%
  \reset@font\fontsize{#1}{#2pt}%
  \fontfamily{#3}\fontseries{#4}\fontshape{#5}%
  \selectfont}%
\fi\endgroup%
\begin{picture}(654,429)(79,332)
\put(406,434){\makebox(0,0)[b]{\smash{{\SetFigFont{8}{9.6}{\familydefault}{\mddefault}{\updefault}{\color[rgb]{0,0,0}$h$}%
}}}}
\end{picture}%
} \:,
\qquad
\raisebox{-3mm}{\begin{picture}(0,0)%
\includegraphics{nf-u-h.pstex}%
\end{picture}%
\setlength{\unitlength}{4144sp}%
\begingroup\makeatletter\ifx\SetFigFont\undefined%
\gdef\SetFigFont#1#2#3#4#5{%
  \reset@font\fontsize{#1}{#2pt}%
  \fontfamily{#3}\fontseries{#4}\fontshape{#5}%
  \selectfont}%
\fi\endgroup%
\begin{picture}(654,429)(-866,-568)
\put(-539,-466){\makebox(0,0)[b]{\smash{{\SetFigFont{8}{9.6}{\familydefault}{\mddefault}{\updefault}{\color[rgb]{0,0,0}$h$}%
}}}}
\end{picture}%
} \:,
\qquad
\raisebox{-3mm}{\begin{picture}(0,0)%
\includegraphics{nf-o-h.pstex}%
\end{picture}%
\setlength{\unitlength}{4144sp}%
\begingroup\makeatletter\ifx\SetFigFont\undefined%
\gdef\SetFigFont#1#2#3#4#5{%
  \reset@font\fontsize{#1}{#2pt}%
  \fontfamily{#3}\fontseries{#4}\fontshape{#5}%
  \selectfont}%
\fi\endgroup%
\begin{picture}(474,429)(-776,-568)
\put(-539,-466){\makebox(0,0)[b]{\smash{{\SetFigFont{8}{9.6}{\familydefault}{\mddefault}{\updefault}{\color[rgb]{0,0,0}$h$}%
}}}}
\end{picture}%
} \:.
$$

\noindent We apply the induction hypothesis to each one, concluding that they all belong to $N$. Thus $f$ is in $N$, built by construction rule (b). Case (iii) is similar to case (ii), with the $2$-cell $\raisebox{-3mm}{\begin{picture}(0,0)%
\includegraphics{C-k-f0.pstex}%
\end{picture}%
\setlength{\unitlength}{4144sp}%
\begingroup\makeatletter\ifx\SetFigFont\undefined%
\gdef\SetFigFont#1#2#3#4#5{%
  \reset@font\fontsize{#1}{#2pt}%
  \fontfamily{#3}\fontseries{#4}\fontshape{#5}%
  \selectfont}%
\fi\endgroup%
\begin{picture}(577,378)(-908,290)
\put(-516,524){\makebox(0,0)[lb]{\smash{{\SetFigFont{8}{9.6}{\familydefault}{\mddefault}{\updefault}{\color[rgb]{0,0,0}$k$}%
}}}}
\end{picture}%
}$ replaced by $\raisebox{-2mm}{\begin{picture}(0,0)%
\includegraphics{n-k-nn.pstex}%
\end{picture}%
\setlength{\unitlength}{4144sp}%
\begingroup\makeatletter\ifx\SetFigFont\undefined%
\gdef\SetFigFont#1#2#3#4#5{%
  \reset@font\fontsize{#1}{#2pt}%
  \fontfamily{#3}\fontseries{#4}\fontshape{#5}%
  \selectfont}%
\fi\endgroup%
\begin{picture}(669,291)(-686,-73)
\put(-202, 74){\makebox(0,0)[lb]{\smash{{\SetFigFont{8}{9.6}{\familydefault}{\mddefault}{\updefault}{\color[rgb]{0,0,0}$k$}%
}}}}
\end{picture}%
}$. In case (iv), the reasoning depends on $\phi$:

\begin{itemize}

\item When $\phi=\figeps{n}$, one has the following possibilities, depending where $\phi$ connects to~$g$:
$$
\begin{array}{c c c c c c c}
\raisebox{-2mm}{\begin{picture}(0,0)%
\includegraphics{nf-f.pstex}%
\end{picture}%
\setlength{\unitlength}{4144sp}%
\begingroup\makeatletter\ifx\SetFigFont\undefined%
\gdef\SetFigFont#1#2#3#4#5{%
  \reset@font\fontsize{#1}{#2pt}%
  \fontfamily{#3}\fontseries{#4}\fontshape{#5}%
  \selectfont}%
\fi\endgroup%
\begin{picture}(294,294)(-911,-1108)
\put(-764,-1006){\makebox(0,0)[b]{\smash{{\SetFigFont{8}{9.6}{\familydefault}{\mddefault}{\updefault}{\color[rgb]{0,0,0}$f$}%
}}}}
\end{picture}%
} 
	& \quad = \quad 
& \raisebox{-3mm}{\begin{picture}(0,0)%
\includegraphics{n-u-k-uu-1.pstex}%
\end{picture}%
\setlength{\unitlength}{4144sp}%
\begingroup\makeatletter\ifx\SetFigFont\undefined%
\gdef\SetFigFont#1#2#3#4#5{%
  \reset@font\fontsize{#1}{#2pt}%
  \fontfamily{#3}\fontseries{#4}\fontshape{#5}%
  \selectfont}%
\fi\endgroup%
\begin{picture}(759,381)(-731,290)
\put(-156,524){\makebox(0,0)[lb]{\smash{{\SetFigFont{8}{9.6}{\familydefault}{\mddefault}{\updefault}{\color[rgb]{0,0,0}$k$}%
}}}}
\end{picture}%
} \raisebox{-2mm}{\begin{picture}(0,0)%
\includegraphics{nf-h.pstex}%
\end{picture}%
\setlength{\unitlength}{4144sp}%
\begingroup\makeatletter\ifx\SetFigFont\undefined%
\gdef\SetFigFont#1#2#3#4#5{%
  \reset@font\fontsize{#1}{#2pt}%
  \fontfamily{#3}\fontseries{#4}\fontshape{#5}%
  \selectfont}%
\fi\endgroup%
\begin{picture}(294,294)(-911,-1108)
\put(-764,-1006){\makebox(0,0)[b]{\smash{{\SetFigFont{8}{9.6}{\familydefault}{\mddefault}{\updefault}{\color[rgb]{0,0,0}$h$}%
}}}}
\end{picture}%
} 
	& \quad \text{or} \quad
& \raisebox{-3mm}{\begin{picture}(0,0)%
\includegraphics{n-u-k-uu-2.pstex}%
\end{picture}%
\setlength{\unitlength}{4144sp}%
\begingroup\makeatletter\ifx\SetFigFont\undefined%
\gdef\SetFigFont#1#2#3#4#5{%
  \reset@font\fontsize{#1}{#2pt}%
  \fontfamily{#3}\fontseries{#4}\fontshape{#5}%
  \selectfont}%
\fi\endgroup%
\begin{picture}(961,381)(-821,290)
\put(-44,569){\makebox(0,0)[lb]{\smash{{\SetFigFont{8}{9.6}{\familydefault}{\mddefault}{\updefault}{\color[rgb]{0,0,0}$k$}%
}}}}
\end{picture}%
} \raisebox{-2mm}{\begin{picture}(0,0)%
\includegraphics{nf-h.pstex}%
\end{picture}%
\setlength{\unitlength}{4144sp}%
\begingroup\makeatletter\ifx\SetFigFont\undefined%
\gdef\SetFigFont#1#2#3#4#5{%
  \reset@font\fontsize{#1}{#2pt}%
  \fontfamily{#3}\fontseries{#4}\fontshape{#5}%
  \selectfont}%
\fi\endgroup%
\begin{picture}(294,294)(-911,-1108)
\put(-764,-1006){\makebox(0,0)[b]{\smash{{\SetFigFont{8}{9.6}{\familydefault}{\mddefault}{\updefault}{\color[rgb]{0,0,0}$h$}%
}}}}
\end{picture}%
} 
	& \quad \text{or} \quad 
& \raisebox{-3mm}{\begin{picture}(0,0)%
\includegraphics{n-u-k-uu-3.pstex}%
\end{picture}%
\setlength{\unitlength}{4144sp}%
\begingroup\makeatletter\ifx\SetFigFont\undefined%
\gdef\SetFigFont#1#2#3#4#5{%
  \reset@font\fontsize{#1}{#2pt}%
  \fontfamily{#3}\fontseries{#4}\fontshape{#5}%
  \selectfont}%
\fi\endgroup%
\begin{picture}(759,381)(-731,290)
\put(-156,524){\makebox(0,0)[lb]{\smash{{\SetFigFont{8}{9.6}{\familydefault}{\mddefault}{\updefault}{\color[rgb]{0,0,0}$k$}%
}}}}
\end{picture}%
} \raisebox{-2mm}{\begin{picture}(0,0)%
\includegraphics{nf-h.pstex}%
\end{picture}%
\setlength{\unitlength}{4144sp}%
\begingroup\makeatletter\ifx\SetFigFont\undefined%
\gdef\SetFigFont#1#2#3#4#5{%
  \reset@font\fontsize{#1}{#2pt}%
  \fontfamily{#3}\fontseries{#4}\fontshape{#5}%
  \selectfont}%
\fi\endgroup%
\begin{picture}(294,294)(-911,-1108)
\put(-764,-1006){\makebox(0,0)[b]{\smash{{\SetFigFont{8}{9.6}{\familydefault}{\mddefault}{\updefault}{\color[rgb]{0,0,0}$h$}%
}}}}
\end{picture}%
} 
\\
	&&& \text{or}
& \raisebox{-3mm}{\input{n-u-k-uu-h.pstex_t}}
	&\text{or}
& \raisebox{-2mm}{\begin{picture}(0,0)%
\includegraphics{u-k-uu.pstex}%
\end{picture}%
\setlength{\unitlength}{4144sp}%
\begingroup\makeatletter\ifx\SetFigFont\undefined%
\gdef\SetFigFont#1#2#3#4#5{%
  \reset@font\fontsize{#1}{#2pt}%
  \fontfamily{#3}\fontseries{#4}\fontshape{#5}%
  \selectfont}%
\fi\endgroup%
\begin{picture}(669,318)(-686,290)
\put(-202,524){\makebox(0,0)[lb]{\smash{{\SetFigFont{8}{9.6}{\familydefault}{\mddefault}{\updefault}{\color[rgb]{0,0,0}$k$}%
}}}}
\end{picture}%
} \raisebox{-3mm}{\begin{picture}(0,0)%
\includegraphics{nf-n-h.pstex}%
\end{picture}%
\setlength{\unitlength}{4144sp}%
\begingroup\makeatletter\ifx\SetFigFont\undefined%
\gdef\SetFigFont#1#2#3#4#5{%
  \reset@font\fontsize{#1}{#2pt}%
  \fontfamily{#3}\fontseries{#4}\fontshape{#5}%
  \selectfont}%
\fi\endgroup%
\begin{picture}(654,429)(79,332)
\put(406,434){\makebox(0,0)[b]{\smash{{\SetFigFont{8}{9.6}{\familydefault}{\mddefault}{\updefault}{\color[rgb]{0,0,0}$h$}%
}}}}
\end{picture}%
} \:.
\end{array}
$$

\noindent The first and third case cannot occur. Indeed, one proves, by structural induction, that a normal form with source of size at least $1$ and with degenerate target has the following shape:
$$
\raisebox{-3mm}{\begin{picture}(0,0)%
\includegraphics{C-k-f0.pstex}%
\end{picture}%
\setlength{\unitlength}{4144sp}%
\begingroup\makeatletter\ifx\SetFigFont\undefined%
\gdef\SetFigFont#1#2#3#4#5{%
  \reset@font\fontsize{#1}{#2pt}%
  \fontfamily{#3}\fontseries{#4}\fontshape{#5}%
  \selectfont}%
\fi\endgroup%
\begin{picture}(577,378)(-908,290)
\put(-516,524){\makebox(0,0)[lb]{\smash{{\SetFigFont{8}{9.6}{\familydefault}{\mddefault}{\updefault}{\color[rgb]{0,0,0}$k$}%
}}}}
\end{picture}%
} \cdots\quad \raisebox{-3mm}{\begin{picture}(0,0)%
\includegraphics{C-k-f0.pstex}%
\end{picture}%
\setlength{\unitlength}{4144sp}%
\begingroup\makeatletter\ifx\SetFigFont\undefined%
\gdef\SetFigFont#1#2#3#4#5{%
  \reset@font\fontsize{#1}{#2pt}%
  \fontfamily{#3}\fontseries{#4}\fontshape{#5}%
  \selectfont}%
\fi\endgroup%
\begin{picture}(577,378)(-908,290)
\put(-516,524){\makebox(0,0)[lb]{\smash{{\SetFigFont{8}{9.6}{\familydefault}{\mddefault}{\updefault}{\color[rgb]{0,0,0}$k$}%
}}}}
\end{picture}%
} \raisebox{-2mm}{\begin{picture}(0,0)%
\includegraphics{u-k-uu.pstex}%
\end{picture}%
\setlength{\unitlength}{4144sp}%
\begingroup\makeatletter\ifx\SetFigFont\undefined%
\gdef\SetFigFont#1#2#3#4#5{%
  \reset@font\fontsize{#1}{#2pt}%
  \fontfamily{#3}\fontseries{#4}\fontshape{#5}%
  \selectfont}%
\fi\endgroup%
\begin{picture}(669,318)(-686,290)
\put(-202,524){\makebox(0,0)[lb]{\smash{{\SetFigFont{8}{9.6}{\familydefault}{\mddefault}{\updefault}{\color[rgb]{0,0,0}$k$}%
}}}}
\end{picture}%
} \:.
$$

\noindent As a consequence, such a decomposition of $f$ would contain either $\figeps{nu}$ or $\figeps{un}$, preventing it from being a normal form.

For the second case, one applies the induction hypothesis to the $2$-cell $\scalefigeps{-2mm}{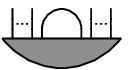}$: indeed, it is a $2$-cell with size at most $m-1$ that must be a normal form, otherwise $f$ would not. Thus, $f$ is built from $2$-cells of $N$ following construction rule (d) and, as such, is in $N$.

The fourth decomposition contains either $\figeps{on}$ or $\figeps{un}$, respectively when $k\geq 1$ and $k=0$. Thus it is not possible that $f$ decomposes this way, since it is a normal form.

For the fifth decomposition, one applies the induction hypothesis to $\raisebox{-3mm}{\begin{picture}(0,0)%
\includegraphics{nf-n-h.pstex}%
\end{picture}%
\setlength{\unitlength}{4144sp}%
\begingroup\makeatletter\ifx\SetFigFont\undefined%
\gdef\SetFigFont#1#2#3#4#5{%
  \reset@font\fontsize{#1}{#2pt}%
  \fontfamily{#3}\fontseries{#4}\fontshape{#5}%
  \selectfont}%
\fi\endgroup%
\begin{picture}(654,429)(79,332)
\put(406,434){\makebox(0,0)[b]{\smash{{\SetFigFont{8}{9.6}{\familydefault}{\mddefault}{\updefault}{\color[rgb]{0,0,0}$h$}%
}}}}
\end{picture}%
}$, which is a $2$-cell that must be a normal form, with size at most $m-1$. 

\item When $\phi=\figeps{u}$, one has the following possible decompositions of $f$:
$$
\raisebox{-2mm}{\begin{picture}(0,0)%
\includegraphics{nf-f.pstex}%
\end{picture}%
\setlength{\unitlength}{4144sp}%
\begingroup\makeatletter\ifx\SetFigFont\undefined%
\gdef\SetFigFont#1#2#3#4#5{%
  \reset@font\fontsize{#1}{#2pt}%
  \fontfamily{#3}\fontseries{#4}\fontshape{#5}%
  \selectfont}%
\fi\endgroup%
\begin{picture}(294,294)(-911,-1108)
\put(-764,-1006){\makebox(0,0)[b]{\smash{{\SetFigFont{8}{9.6}{\familydefault}{\mddefault}{\updefault}{\color[rgb]{0,0,0}$f$}%
}}}}
\end{picture}%
} 
	\quad = \quad 
\figeps{u} \quad \raisebox{-3mm}{\begin{picture}(0,0)%
\includegraphics{u-k-uu.pstex}%
\end{picture}%
\setlength{\unitlength}{4144sp}%
\begingroup\makeatletter\ifx\SetFigFont\undefined%
\gdef\SetFigFont#1#2#3#4#5{%
  \reset@font\fontsize{#1}{#2pt}%
  \fontfamily{#3}\fontseries{#4}\fontshape{#5}%
  \selectfont}%
\fi\endgroup%
\begin{picture}(669,318)(-686,290)
\put(-202,524){\makebox(0,0)[lb]{\smash{{\SetFigFont{8}{9.6}{\familydefault}{\mddefault}{\updefault}{\color[rgb]{0,0,0}$k$}%
}}}}
\end{picture}%
} \raisebox{-2mm}{\begin{picture}(0,0)%
\includegraphics{nf-h.pstex}%
\end{picture}%
\setlength{\unitlength}{4144sp}%
\begingroup\makeatletter\ifx\SetFigFont\undefined%
\gdef\SetFigFont#1#2#3#4#5{%
  \reset@font\fontsize{#1}{#2pt}%
  \fontfamily{#3}\fontseries{#4}\fontshape{#5}%
  \selectfont}%
\fi\endgroup%
\begin{picture}(294,294)(-911,-1108)
\put(-764,-1006){\makebox(0,0)[b]{\smash{{\SetFigFont{8}{9.6}{\familydefault}{\mddefault}{\updefault}{\color[rgb]{0,0,0}$h$}%
}}}}
\end{picture}%
} 
	\quad \text{or} \quad
\raisebox{-3mm}{\begin{picture}(0,0)%
\includegraphics{u-u-k-uu.pstex}%
\end{picture}%
\setlength{\unitlength}{4144sp}%
\begingroup\makeatletter\ifx\SetFigFont\undefined%
\gdef\SetFigFont#1#2#3#4#5{%
  \reset@font\fontsize{#1}{#2pt}%
  \fontfamily{#3}\fontseries{#4}\fontshape{#5}%
  \selectfont}%
\fi\endgroup%
\begin{picture}(961,381)(-821,290)
\put(-44,569){\makebox(0,0)[lb]{\smash{{\SetFigFont{8}{9.6}{\familydefault}{\mddefault}{\updefault}{\color[rgb]{0,0,0}$k$}%
}}}}
\end{picture}%
} \raisebox{-2mm}{\begin{picture}(0,0)%
\includegraphics{nf-h.pstex}%
\end{picture}%
\setlength{\unitlength}{4144sp}%
\begingroup\makeatletter\ifx\SetFigFont\undefined%
\gdef\SetFigFont#1#2#3#4#5{%
  \reset@font\fontsize{#1}{#2pt}%
  \fontfamily{#3}\fontseries{#4}\fontshape{#5}%
  \selectfont}%
\fi\endgroup%
\begin{picture}(294,294)(-911,-1108)
\put(-764,-1006){\makebox(0,0)[b]{\smash{{\SetFigFont{8}{9.6}{\familydefault}{\mddefault}{\updefault}{\color[rgb]{0,0,0}$h$}%
}}}}
\end{picture}%
} 
	\quad \text{or} \quad 
\raisebox{-3mm}{\begin{picture}(0,0)%
\includegraphics{u-k-uu.pstex}%
\end{picture}%
\setlength{\unitlength}{4144sp}%
\begingroup\makeatletter\ifx\SetFigFont\undefined%
\gdef\SetFigFont#1#2#3#4#5{%
  \reset@font\fontsize{#1}{#2pt}%
  \fontfamily{#3}\fontseries{#4}\fontshape{#5}%
  \selectfont}%
\fi\endgroup%
\begin{picture}(669,318)(-686,290)
\put(-202,524){\makebox(0,0)[lb]{\smash{{\SetFigFont{8}{9.6}{\familydefault}{\mddefault}{\updefault}{\color[rgb]{0,0,0}$k$}%
}}}}
\end{picture}%
} \raisebox{-3mm}{\begin{picture}(0,0)%
\includegraphics{nf-u-h.pstex}%
\end{picture}%
\setlength{\unitlength}{4144sp}%
\begingroup\makeatletter\ifx\SetFigFont\undefined%
\gdef\SetFigFont#1#2#3#4#5{%
  \reset@font\fontsize{#1}{#2pt}%
  \fontfamily{#3}\fontseries{#4}\fontshape{#5}%
  \selectfont}%
\fi\endgroup%
\begin{picture}(654,429)(-866,-568)
\put(-539,-466){\makebox(0,0)[b]{\smash{{\SetFigFont{8}{9.6}{\familydefault}{\mddefault}{\updefault}{\color[rgb]{0,0,0}$h$}%
}}}}
\end{picture}%
} \:. 
$$

\noindent The first case shows that $f$ is in $N$: indeed, it is built with construction rule (d), applied with $\figeps{uu}=\ast$, $k=0$ and $\scalefigeps{-2mm}{normal}=\raisebox{-3mm}{\begin{picture}(0,0)%
\includegraphics{u-k-uu.pstex}%
\end{picture}%
\setlength{\unitlength}{4144sp}%
\begingroup\makeatletter\ifx\SetFigFont\undefined%
\gdef\SetFigFont#1#2#3#4#5{%
  \reset@font\fontsize{#1}{#2pt}%
  \fontfamily{#3}\fontseries{#4}\fontshape{#5}%
  \selectfont}%
\fi\endgroup%
\begin{picture}(669,318)(-686,290)
\put(-202,524){\makebox(0,0)[lb]{\smash{{\SetFigFont{8}{9.6}{\familydefault}{\mddefault}{\updefault}{\color[rgb]{0,0,0}$k$}%
}}}}
\end{picture}%
} \raisebox{-2mm}{\begin{picture}(0,0)%
\includegraphics{nf-h.pstex}%
\end{picture}%
\setlength{\unitlength}{4144sp}%
\begingroup\makeatletter\ifx\SetFigFont\undefined%
\gdef\SetFigFont#1#2#3#4#5{%
  \reset@font\fontsize{#1}{#2pt}%
  \fontfamily{#3}\fontseries{#4}\fontshape{#5}%
  \selectfont}%
\fi\endgroup%
\begin{picture}(294,294)(-911,-1108)
\put(-764,-1006){\makebox(0,0)[b]{\smash{{\SetFigFont{8}{9.6}{\familydefault}{\mddefault}{\updefault}{\color[rgb]{0,0,0}$h$}%
}}}}
\end{picture}%
}$, which is $g$ and, as such, belongs to $N$.

In the second case, we apply the induction hypothesis to $\scalefigeps{-2mm}{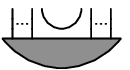}$: it is a normal form of size at most $m-1$. Thus $f$ is built with construction rule (d).

In the third case, one applies the induction hypothesis to $\raisebox{-3mm}{\begin{picture}(0,0)%
\includegraphics{nf-u-h.pstex}%
\end{picture}%
\setlength{\unitlength}{4144sp}%
\begingroup\makeatletter\ifx\SetFigFont\undefined%
\gdef\SetFigFont#1#2#3#4#5{%
  \reset@font\fontsize{#1}{#2pt}%
  \fontfamily{#3}\fontseries{#4}\fontshape{#5}%
  \selectfont}%
\fi\endgroup%
\begin{picture}(654,429)(-866,-568)
\put(-539,-466){\makebox(0,0)[b]{\smash{{\SetFigFont{8}{9.6}{\familydefault}{\mddefault}{\updefault}{\color[rgb]{0,0,0}$h$}%
}}}}
\end{picture}%
}$: it is a normal form of size at most $m-1$. We conclude that $f$ is built with construction rule (d).

\item When $\phi=\figeps{o}$, the possible decompositions of $f$ are:
$$
\begin{array}{c c c c c}
\raisebox{-2mm}{\begin{picture}(0,0)%
\includegraphics{nf-f.pstex}%
\end{picture}%
\setlength{\unitlength}{4144sp}%
\begingroup\makeatletter\ifx\SetFigFont\undefined%
\gdef\SetFigFont#1#2#3#4#5{%
  \reset@font\fontsize{#1}{#2pt}%
  \fontfamily{#3}\fontseries{#4}\fontshape{#5}%
  \selectfont}%
\fi\endgroup%
\begin{picture}(294,294)(-911,-1108)
\put(-764,-1006){\makebox(0,0)[b]{\smash{{\SetFigFont{8}{9.6}{\familydefault}{\mddefault}{\updefault}{\color[rgb]{0,0,0}$f$}%
}}}}
\end{picture}%
} 
	& \quad = \quad 
& \raisebox{-3mm}{\begin{picture}(0,0)%
\includegraphics{o-u-k-uu-1.pstex}%
\end{picture}%
\setlength{\unitlength}{4144sp}%
\begingroup\makeatletter\ifx\SetFigFont\undefined%
\gdef\SetFigFont#1#2#3#4#5{%
  \reset@font\fontsize{#1}{#2pt}%
  \fontfamily{#3}\fontseries{#4}\fontshape{#5}%
  \selectfont}%
\fi\endgroup%
\begin{picture}(699,318)(-716,290)
\put(-202,524){\makebox(0,0)[lb]{\smash{{\SetFigFont{8}{9.6}{\familydefault}{\mddefault}{\updefault}{\color[rgb]{0,0,0}$k$}%
}}}}
\end{picture}%
} \raisebox{-2mm}{\begin{picture}(0,0)%
\includegraphics{nf-h.pstex}%
\end{picture}%
\setlength{\unitlength}{4144sp}%
\begingroup\makeatletter\ifx\SetFigFont\undefined%
\gdef\SetFigFont#1#2#3#4#5{%
  \reset@font\fontsize{#1}{#2pt}%
  \fontfamily{#3}\fontseries{#4}\fontshape{#5}%
  \selectfont}%
\fi\endgroup%
\begin{picture}(294,294)(-911,-1108)
\put(-764,-1006){\makebox(0,0)[b]{\smash{{\SetFigFont{8}{9.6}{\familydefault}{\mddefault}{\updefault}{\color[rgb]{0,0,0}$h$}%
}}}}
\end{picture}%
} 
	& \quad \text{or} \quad
& \raisebox{-4mm}{\begin{picture}(0,0)%
\includegraphics{o-u-k-uu-2.pstex}%
\end{picture}%
\setlength{\unitlength}{4144sp}%
\begingroup\makeatletter\ifx\SetFigFont\undefined%
\gdef\SetFigFont#1#2#3#4#5{%
  \reset@font\fontsize{#1}{#2pt}%
  \fontfamily{#3}\fontseries{#4}\fontshape{#5}%
  \selectfont}%
\fi\endgroup%
\begin{picture}(848,426)(-776,290)
\put(-112,614){\makebox(0,0)[lb]{\smash{{\SetFigFont{8}{9.6}{\familydefault}{\mddefault}{\updefault}{\color[rgb]{0,0,0}$k$}%
}}}}
\end{picture}%
} \raisebox{-2mm}{\begin{picture}(0,0)%
\includegraphics{nf-h.pstex}%
\end{picture}%
\setlength{\unitlength}{4144sp}%
\begingroup\makeatletter\ifx\SetFigFont\undefined%
\gdef\SetFigFont#1#2#3#4#5{%
  \reset@font\fontsize{#1}{#2pt}%
  \fontfamily{#3}\fontseries{#4}\fontshape{#5}%
  \selectfont}%
\fi\endgroup%
\begin{picture}(294,294)(-911,-1108)
\put(-764,-1006){\makebox(0,0)[b]{\smash{{\SetFigFont{8}{9.6}{\familydefault}{\mddefault}{\updefault}{\color[rgb]{0,0,0}$h$}%
}}}}
\end{picture}%
} 
\\
& \text{or}
& \raisebox{-3mm}{\begin{picture}(0,0)%
\includegraphics{o-u-k-uu-3.pstex}%
\end{picture}%
\setlength{\unitlength}{4144sp}%
\begingroup\makeatletter\ifx\SetFigFont\undefined%
\gdef\SetFigFont#1#2#3#4#5{%
  \reset@font\fontsize{#1}{#2pt}%
  \fontfamily{#3}\fontseries{#4}\fontshape{#5}%
  \selectfont}%
\fi\endgroup%
\begin{picture}(807,318)(-686,290)
\put(-202,524){\makebox(0,0)[lb]{\smash{{\SetFigFont{8}{9.6}{\familydefault}{\mddefault}{\updefault}{\color[rgb]{0,0,0}$k+1$}%
}}}}
\end{picture}%
} \raisebox{-2mm}{\begin{picture}(0,0)%
\includegraphics{nf-h.pstex}%
\end{picture}%
\setlength{\unitlength}{4144sp}%
\begingroup\makeatletter\ifx\SetFigFont\undefined%
\gdef\SetFigFont#1#2#3#4#5{%
  \reset@font\fontsize{#1}{#2pt}%
  \fontfamily{#3}\fontseries{#4}\fontshape{#5}%
  \selectfont}%
\fi\endgroup%
\begin{picture}(294,294)(-911,-1108)
\put(-764,-1006){\makebox(0,0)[b]{\smash{{\SetFigFont{8}{9.6}{\familydefault}{\mddefault}{\updefault}{\color[rgb]{0,0,0}$h$}%
}}}}
\end{picture}%
} 
	& \quad \text{or} \quad
& \raisebox{-2mm}{\begin{picture}(0,0)%
\includegraphics{u-k-uu.pstex}%
\end{picture}%
\setlength{\unitlength}{4144sp}%
\begingroup\makeatletter\ifx\SetFigFont\undefined%
\gdef\SetFigFont#1#2#3#4#5{%
  \reset@font\fontsize{#1}{#2pt}%
  \fontfamily{#3}\fontseries{#4}\fontshape{#5}%
  \selectfont}%
\fi\endgroup%
\begin{picture}(669,318)(-686,290)
\put(-202,524){\makebox(0,0)[lb]{\smash{{\SetFigFont{8}{9.6}{\familydefault}{\mddefault}{\updefault}{\color[rgb]{0,0,0}$k$}%
}}}}
\end{picture}%
} \raisebox{-3mm}{\begin{picture}(0,0)%
\includegraphics{nf-o-h.pstex}%
\end{picture}%
\setlength{\unitlength}{4144sp}%
\begingroup\makeatletter\ifx\SetFigFont\undefined%
\gdef\SetFigFont#1#2#3#4#5{%
  \reset@font\fontsize{#1}{#2pt}%
  \fontfamily{#3}\fontseries{#4}\fontshape{#5}%
  \selectfont}%
\fi\endgroup%
\begin{picture}(474,429)(-776,-568)
\put(-539,-466){\makebox(0,0)[b]{\smash{{\SetFigFont{8}{9.6}{\familydefault}{\mddefault}{\updefault}{\color[rgb]{0,0,0}$h$}%
}}}}
\end{picture}%
} \:.
\end{array}
$$

\noindent The first case cannot occur: otherwise, $f$ would contain $\figeps{ou}$ and, thus, it would not be a normal form.

In the second case, we apply the induction hypothesis to $\scalefigeps{-2mm}{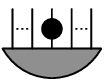}$: this is a normal form with size at most $m-1$. This proves that $f$ is in $N$, built following construction rule (d).

In the third case, $f$ is in $N$, built following construction rule (d).

In the fourth case, we apply the induction hypothesis to $\raisebox{-3mm}{\begin{picture}(0,0)%
\includegraphics{nf-o-h.pstex}%
\end{picture}%
\setlength{\unitlength}{4144sp}%
\begingroup\makeatletter\ifx\SetFigFont\undefined%
\gdef\SetFigFont#1#2#3#4#5{%
  \reset@font\fontsize{#1}{#2pt}%
  \fontfamily{#3}\fontseries{#4}\fontshape{#5}%
  \selectfont}%
\fi\endgroup%
\begin{picture}(474,429)(-776,-568)
\put(-539,-466){\makebox(0,0)[b]{\smash{{\SetFigFont{8}{9.6}{\familydefault}{\mddefault}{\updefault}{\color[rgb]{0,0,0}$h$}%
}}}}
\end{picture}%
}$: this is a normal form with size at most $m-1$. Thus $f$ is in $N$, built with construction rule (d).
\end{itemize}

\noindent The final case~(v) also depends on the values of $\phi$: 
\begin{itemize}
\item When $\phi=\figeps{n}$, we have the following possible decompositions of $f$:
$$
\raisebox{-2mm}{\begin{picture}(0,0)%
\includegraphics{nf-f.pstex}%
\end{picture}%
\setlength{\unitlength}{4144sp}%
\begingroup\makeatletter\ifx\SetFigFont\undefined%
\gdef\SetFigFont#1#2#3#4#5{%
  \reset@font\fontsize{#1}{#2pt}%
  \fontfamily{#3}\fontseries{#4}\fontshape{#5}%
  \selectfont}%
\fi\endgroup%
\begin{picture}(294,294)(-911,-1108)
\put(-764,-1006){\makebox(0,0)[b]{\smash{{\SetFigFont{8}{9.6}{\familydefault}{\mddefault}{\updefault}{\color[rgb]{0,0,0}$f$}%
}}}}
\end{picture}%
} 
	\quad = \quad 
\raisebox{-3mm}{\input{n-k-h.pstex_t}} 
	\quad \text{or} \quad
\figeps{o}^k \raisebox{-3mm}{\begin{picture}(0,0)%
\includegraphics{nf-n-h.pstex}%
\end{picture}%
\setlength{\unitlength}{4144sp}%
\begingroup\makeatletter\ifx\SetFigFont\undefined%
\gdef\SetFigFont#1#2#3#4#5{%
  \reset@font\fontsize{#1}{#2pt}%
  \fontfamily{#3}\fontseries{#4}\fontshape{#5}%
  \selectfont}%
\fi\endgroup%
\begin{picture}(654,429)(79,332)
\put(406,434){\makebox(0,0)[b]{\smash{{\SetFigFont{8}{9.6}{\familydefault}{\mddefault}{\updefault}{\color[rgb]{0,0,0}$h$}%
}}}}
\end{picture}%
} \:.
$$

\noindent In the first case, one must have $k=0$: otherwise, $f$ would contain $\figeps{on}$ which is not a normal form. Thus the $2$-cell $h$ is a normal form of size $m-1$: we apply the induction hypothesis to get that $h$ is in $N$. Then, by structural induction on $h$, one shows that it has one of the following two shapes:
$$
\raisebox{-2mm}{\begin{picture}(0,0)%
\includegraphics{nf-h.pstex}%
\end{picture}%
\setlength{\unitlength}{4144sp}%
\begingroup\makeatletter\ifx\SetFigFont\undefined%
\gdef\SetFigFont#1#2#3#4#5{%
  \reset@font\fontsize{#1}{#2pt}%
  \fontfamily{#3}\fontseries{#4}\fontshape{#5}%
  \selectfont}%
\fi\endgroup%
\begin{picture}(294,294)(-911,-1108)
\put(-764,-1006){\makebox(0,0)[b]{\smash{{\SetFigFont{8}{9.6}{\familydefault}{\mddefault}{\updefault}{\color[rgb]{0,0,0}$h$}%
}}}}
\end{picture}%
} 
	\quad = \quad 
\scalefigeps{-2mm}{nn}\: \raisebox{-3mm}{\begin{picture}(0,0)%
\includegraphics{u-k-uu.pstex}%
\end{picture}%
\setlength{\unitlength}{4144sp}%
\begingroup\makeatletter\ifx\SetFigFont\undefined%
\gdef\SetFigFont#1#2#3#4#5{%
  \reset@font\fontsize{#1}{#2pt}%
  \fontfamily{#3}\fontseries{#4}\fontshape{#5}%
  \selectfont}%
\fi\endgroup%
\begin{picture}(669,318)(-686,290)
\put(-202,524){\makebox(0,0)[lb]{\smash{{\SetFigFont{8}{9.6}{\familydefault}{\mddefault}{\updefault}{\color[rgb]{0,0,0}$k$}%
}}}}
\end{picture}%
} \scalefigeps{-2mm}{normal}
	\quad\text{or}\quad
\scalefigeps{-2mm}{nn} \: \figeps{o}^k \scalefigeps{-2mm}{normal} \:.
$$

\noindent The first decomposition is impossible since, otherwise, $f$ would contain $\figeps{nu}$ and, thus, it would not be a normal form. The second decomposition gives that $f$ is in $N$, built from case (c).

In the second case, the $2$-cell $\raisebox{-3mm}{\begin{picture}(0,0)%
\includegraphics{nf-n-h.pstex}%
\end{picture}%
\setlength{\unitlength}{4144sp}%
\begingroup\makeatletter\ifx\SetFigFont\undefined%
\gdef\SetFigFont#1#2#3#4#5{%
  \reset@font\fontsize{#1}{#2pt}%
  \fontfamily{#3}\fontseries{#4}\fontshape{#5}%
  \selectfont}%
\fi\endgroup%
\begin{picture}(654,429)(79,332)
\put(406,434){\makebox(0,0)[b]{\smash{{\SetFigFont{8}{9.6}{\familydefault}{\mddefault}{\updefault}{\color[rgb]{0,0,0}$h$}%
}}}}
\end{picture}%
}$ is a normal form. Moreover, if $k\geq 1$, it has size at most $m-1$, and, if $k=0$, it has size $m$, while its source and target have sizes $n-1$ and $p-1$, respectively. Thus, in either situation, we can apply the induction hypothesis to conclude that this $2$-cell is in $N$. As a consequence, $f$ is in $N$, built with construction rule (e).

\item When $\phi=\figeps{u}$, we have the following possible decompositions of $f$:
$$
\raisebox{-2mm}{\begin{picture}(0,0)%
\includegraphics{nf-f.pstex}%
\end{picture}%
\setlength{\unitlength}{4144sp}%
\begingroup\makeatletter\ifx\SetFigFont\undefined%
\gdef\SetFigFont#1#2#3#4#5{%
  \reset@font\fontsize{#1}{#2pt}%
  \fontfamily{#3}\fontseries{#4}\fontshape{#5}%
  \selectfont}%
\fi\endgroup%
\begin{picture}(294,294)(-911,-1108)
\put(-764,-1006){\makebox(0,0)[b]{\smash{{\SetFigFont{8}{9.6}{\familydefault}{\mddefault}{\updefault}{\color[rgb]{0,0,0}$f$}%
}}}}
\end{picture}%
} 
	\quad = \quad 
\figeps{u} \: \figeps{o}^k \raisebox{-2mm}{\begin{picture}(0,0)%
\includegraphics{nf-h.pstex}%
\end{picture}%
\setlength{\unitlength}{4144sp}%
\begingroup\makeatletter\ifx\SetFigFont\undefined%
\gdef\SetFigFont#1#2#3#4#5{%
  \reset@font\fontsize{#1}{#2pt}%
  \fontfamily{#3}\fontseries{#4}\fontshape{#5}%
  \selectfont}%
\fi\endgroup%
\begin{picture}(294,294)(-911,-1108)
\put(-764,-1006){\makebox(0,0)[b]{\smash{{\SetFigFont{8}{9.6}{\familydefault}{\mddefault}{\updefault}{\color[rgb]{0,0,0}$h$}%
}}}}
\end{picture}%
} 
	\quad \text{or} \quad
\figeps{o}^k \raisebox{-3mm}{\begin{picture}(0,0)%
\includegraphics{nf-u-h.pstex}%
\end{picture}%
\setlength{\unitlength}{4144sp}%
\begingroup\makeatletter\ifx\SetFigFont\undefined%
\gdef\SetFigFont#1#2#3#4#5{%
  \reset@font\fontsize{#1}{#2pt}%
  \fontfamily{#3}\fontseries{#4}\fontshape{#5}%
  \selectfont}%
\fi\endgroup%
\begin{picture}(654,429)(-866,-568)
\put(-539,-466){\makebox(0,0)[b]{\smash{{\SetFigFont{8}{9.6}{\familydefault}{\mddefault}{\updefault}{\color[rgb]{0,0,0}$h$}%
}}}}
\end{picture}%
} \:.
$$

\noindent In the first case, $f$ is in $N$, built from $h$ in two subsequent steps, with construction rules (e), then~(d).

In the second case, one can apply the induction hypothesis to $\raisebox{-3mm}{\begin{picture}(0,0)%
\includegraphics{nf-u-h.pstex}%
\end{picture}%
\setlength{\unitlength}{4144sp}%
\begingroup\makeatletter\ifx\SetFigFont\undefined%
\gdef\SetFigFont#1#2#3#4#5{%
  \reset@font\fontsize{#1}{#2pt}%
  \fontfamily{#3}\fontseries{#4}\fontshape{#5}%
  \selectfont}%
\fi\endgroup%
\begin{picture}(654,429)(-866,-568)
\put(-539,-466){\makebox(0,0)[b]{\smash{{\SetFigFont{8}{9.6}{\familydefault}{\mddefault}{\updefault}{\color[rgb]{0,0,0}$h$}%
}}}}
\end{picture}%
}$. Indeed, it is a  normal form, with either size at most $m-1$, when $k\geq 1$, or with size $m$ and source and target of sizes $n-1$ and $p-1$, respectively. Thus this $2$-cell is in $N$, and so does $f$, which is built following construction rule (e).

\item When $\phi=\figeps{o}$, we have the following possible decompositions of $f$:
$$
\raisebox{-2mm}{\begin{picture}(0,0)%
\includegraphics{nf-f.pstex}%
\end{picture}%
\setlength{\unitlength}{4144sp}%
\begingroup\makeatletter\ifx\SetFigFont\undefined%
\gdef\SetFigFont#1#2#3#4#5{%
  \reset@font\fontsize{#1}{#2pt}%
  \fontfamily{#3}\fontseries{#4}\fontshape{#5}%
  \selectfont}%
\fi\endgroup%
\begin{picture}(294,294)(-911,-1108)
\put(-764,-1006){\makebox(0,0)[b]{\smash{{\SetFigFont{8}{9.6}{\familydefault}{\mddefault}{\updefault}{\color[rgb]{0,0,0}$f$}%
}}}}
\end{picture}%
} 
	\quad = \quad 
\figeps{o}^{k+1} \raisebox{-2mm}{\begin{picture}(0,0)%
\includegraphics{nf-h.pstex}%
\end{picture}%
\setlength{\unitlength}{4144sp}%
\begingroup\makeatletter\ifx\SetFigFont\undefined%
\gdef\SetFigFont#1#2#3#4#5{%
  \reset@font\fontsize{#1}{#2pt}%
  \fontfamily{#3}\fontseries{#4}\fontshape{#5}%
  \selectfont}%
\fi\endgroup%
\begin{picture}(294,294)(-911,-1108)
\put(-764,-1006){\makebox(0,0)[b]{\smash{{\SetFigFont{8}{9.6}{\familydefault}{\mddefault}{\updefault}{\color[rgb]{0,0,0}$h$}%
}}}}
\end{picture}%
} 
	\quad \text{or} \quad
\figeps{o}^k \raisebox{-3mm}{\begin{picture}(0,0)%
\includegraphics{nf-o-h.pstex}%
\end{picture}%
\setlength{\unitlength}{4144sp}%
\begingroup\makeatletter\ifx\SetFigFont\undefined%
\gdef\SetFigFont#1#2#3#4#5{%
  \reset@font\fontsize{#1}{#2pt}%
  \fontfamily{#3}\fontseries{#4}\fontshape{#5}%
  \selectfont}%
\fi\endgroup%
\begin{picture}(474,429)(-776,-568)
\put(-539,-466){\makebox(0,0)[b]{\smash{{\SetFigFont{8}{9.6}{\familydefault}{\mddefault}{\updefault}{\color[rgb]{0,0,0}$h$}%
}}}}
\end{picture}%
} \:.
$$

\noindent In the first case, $f$ is built from $h$ by application of construction rule (e) and, as such, is in $N$. 

In the second case, one applies the induction hypothesis to $\raisebox{-3mm}{\begin{picture}(0,0)%
\includegraphics{nf-o-h.pstex}%
\end{picture}%
\setlength{\unitlength}{4144sp}%
\begingroup\makeatletter\ifx\SetFigFont\undefined%
\gdef\SetFigFont#1#2#3#4#5{%
  \reset@font\fontsize{#1}{#2pt}%
  \fontfamily{#3}\fontseries{#4}\fontshape{#5}%
  \selectfont}%
\fi\endgroup%
\begin{picture}(474,429)(-776,-568)
\put(-539,-466){\makebox(0,0)[b]{\smash{{\SetFigFont{8}{9.6}{\familydefault}{\mddefault}{\updefault}{\color[rgb]{0,0,0}$h$}%
}}}}
\end{picture}%
}$, which is a  normal form, with either size at most $m-1$, when $k\geq 1$, or with size $m$ and source and target of sizes $n-1$ and $p-1$, respectively. As a consequence, this $2$-cell is in $N$, proving that $f$ is built following construction rule (e) and, thus, it is in $N$.
\end{itemize}

\noindent To conclude, we have proved that the normal forms of $\Sigma^*$ are exactly the $2$-cells of $N$. In particular, we denote by~$N_0$ the set of normal forms with degenerate source and target. From the inductive scheme defining~$N$, we deduce that the following two construction rules characterize~$N_0$:
$$
\figeps{f0} 
	\quad = \quad 
\ast 
	\quad\text{or}\quad
\raisebox{-3mm}{\begin{picture}(0,0)%
\includegraphics{C-k-f0.pstex}%
\end{picture}%
\setlength{\unitlength}{4144sp}%
\begingroup\makeatletter\ifx\SetFigFont\undefined%
\gdef\SetFigFont#1#2#3#4#5{%
  \reset@font\fontsize{#1}{#2pt}%
  \fontfamily{#3}\fontseries{#4}\fontshape{#5}%
  \selectfont}%
\fi\endgroup%
\begin{picture}(577,378)(-908,290)
\put(-516,524){\makebox(0,0)[lb]{\smash{{\SetFigFont{8}{9.6}{\familydefault}{\mddefault}{\updefault}{\color[rgb]{0,0,0}$k$}%
}}}}
\end{picture}%
} \figeps{f0} \:.
$$

\subsubsection{Confluence}

Let us examine the critical branchings of $\Sigma$. The $3$-polygraph $\Sigma$ has four regular critical branchings, whose sources are:
$$
\scalefigeps{-2mm}{nun}  \:, 
\qquad
\scalefigeps{-2mm}{unu} \:,
\qquad  
\scalefigeps{-2mm}{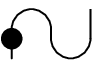} \:,
\qquad
\scalefigeps{-2mm}{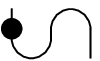} \:.
$$

\noindent It also has one right-indexed critical branching, generated by the $3$-cells $\alpha$ and $\beta$, with source:
$$
\raisebox{-3mm}{\begin{picture}(0,0)%
\includegraphics{onouf.pstex}%
\end{picture}%
\setlength{\unitlength}{4144sp}%
\begingroup\makeatletter\ifx\SetFigFont\undefined%
\gdef\SetFigFont#1#2#3#4#5{%
  \reset@font\fontsize{#1}{#2pt}%
  \fontfamily{#3}\fontseries{#4}\fontshape{#5}%
  \selectfont}%
\fi\endgroup%
\begin{picture}(560,474)(-907,242)
\put(-539,434){\makebox(0,0)[b]{\smash{{\SetFigFont{8}{9.6}{\familydefault}{\mddefault}{\updefault}{\color[rgb]{0,0,0}$k$}%
}}}}
\end{picture}%
} \:.
$$

\noindent Thus $\Sigma$ is a terminating and right-indexed $3$-polygraph. By application of Theorem~\ref{Theorem:RightIndexedConfluence}, we get confluence of $\Sigma$ by proving that its four regular critical branchings and all normal instances of its right-indexed critical branchings are confluent.

For the regular ones, we have the following confluence diagrams:
$$
\xymatrix{
{\scalefigeps{-2mm}{nun}}
	\ar@3 @/^6ex/ [rr] ^{\gamma} _{\hole}="1"
	\ar@3 @/_6ex/ [rr] _{\delta} ^{\hole}="2"
&& {\figeps{n}}
\ar@4 "1";"2" |-*+[o]{\scriptstyle \gamma\delta}
}
\qquad\qquad
\xymatrix{
{\scalefigeps{-2mm}{unu}}
	\ar@3 @/^6ex/ [rr] ^{\delta} _{\hole}="1"
	\ar@3 @/_6ex/ [rr] _{\gamma} ^{\hole}="2"
&& {\figeps{u}}
\ar@4 "1";"2" |-*+[o]{\scriptstyle \delta\gamma}
}
$$
$$
\xymatrix{
& {\scalefigeps{-2mm}{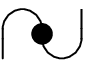}}
	\ar@3 [r] ^{\beta} _-{\hole}="source"
& {\scalefigeps{-2mm}{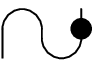}}
	\ar@3 [dr] ^{\gamma} \\
{\scalefigeps{-2mm}{onu}}
	\ar@3 [ur] ^{\alpha}
	\ar@3 [rrr] _{\gamma} ^-{\hole}="target"
&&& {\figeps{o}}
\ar@4 "source" ; "target" |-*+[o]{\scriptstyle \alpha\gamma}
}
\qquad\qquad
\xymatrix{
& {\scalefigeps{-2mm}{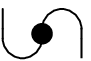}}
	\ar@3 [r] ^{\alpha} _-{\hole}="source"
& {\scalefigeps{-2mm}{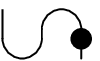}}
	\ar@3 [dr] ^{\delta} \\
{\scalefigeps{-2mm}{oun}}
	\ar@3 [ur] ^{\beta}
	\ar@3 [rrr] _{\delta} ^-{\hole}="target"
&&& {\figeps{o}}
\ar@4 "source" ; "target" |-*+[o]{\scriptstyle \beta\delta}
}
$$

\noindent From the characterization of normal forms of $\Sigma$, the normal instances of the right-indexed critical branching $\alpha\beta\left(\figeps{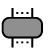}\right)$ are the instances corresponding to the following $2$-cells where, in the latter, $\figeps{f0}$ and $n$ respectively range over $N_0$ and $\Nb$:
$$
\figeps{f} \: = \: \scalefigeps{-3mm}{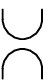} \:,
\qquad
\figeps{f} \: = \: \figeps{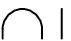} \:,
\qquad
\figeps{f} \: = \: \figeps{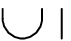} \:,
\qquad
\figeps{f} \: = \: \figeps{f0} \: \figeps{o} ^n \:.
$$

\noindent Now we check that, for each one of these
$2$-cells, the corresponding critical branching
$\alpha\beta\left(\figeps{f}\right)$ is confluent. Let us
note that, for the first three cases, there are several
possible confluence diagrams, because they also contain
regular critical branchings of $\Sigma$. 

\noindent For $\figeps{f} \:=\: \scalefigeps{-3mm}{u1n} \:$,
we choose the following one: 
$$
\xymatrix{
& {\scalefigeps{-3mm}{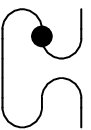}}
	\ar@3 [r] ^{\delta}
& {\figeps{nou}}
	\ar@3 [r] ^{\beta}
& {\figeps{nuo}}
	\ar@3 [dr] ^{\gamma}
\\
{\scalefigeps{-3mm}{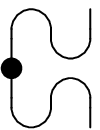}}
	\ar@3 [ur] ^{\alpha} 
	\ar@3 [dr] _{\beta} 
&&&& {\figeps{o}}
\\
& {\scalefigeps{-3mm}{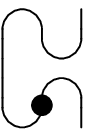}}
	\ar@3 [r] _{\gamma}
& {\figeps{uon}}
	\ar@3 [r] _{\alpha}
& {\figeps{uno}}
	\ar@3 [ur] _{\delta}
\ar@4 "1,3";"3,3" |-*+[o]{\scriptstyle \alpha\beta\left(\raisebox{-1.5mm}{\smallfigeps{u1n}}\right)}
}
$$

\noindent For $\figeps{f} \: = \: \figeps{n0i} \:$:
$$
\xymatrix{
&& {\scalefigeps{-2mm}{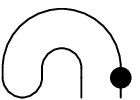}}
	\ar@3 [drr] ^{\delta} 
\\
{\scalefigeps{-2mm}{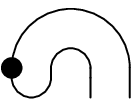}}
	\ar@3 [urr] ^{\alpha}
	\ar@3 [dr] _{\beta}
&&&& {\figeps{no}}
\\
& {\scalefigeps{-2mm}{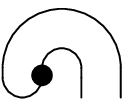}}
	\ar@3 [r] _{\alpha} 
& {\scalefigeps{-2mm}{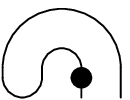}}
	\ar@3 [r] _{\delta}
& {\figeps{on}}
	\ar@3 [ur] _{\alpha}
\ar@4 "1,3";"3,3" |-*+[o]{\scriptstyle \alpha\beta\left(\smallfigeps{n0i}\right)} 
}
$$

\noindent For $\figeps{f} \: = \: \figeps{u0i} \:$:
$$
\xymatrix{
& {\scalefigeps{-2mm}{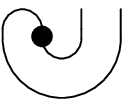}}
	\ar@3 [r] ^{\beta} 
& {\scalefigeps{-2mm}{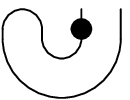}}
	\ar@3 [r] ^{\gamma}
& {\figeps{ou}}
	\ar@3 [dr] ^{\beta}
\\
{\scalefigeps{-2mm}{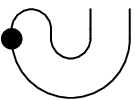}}
	\ar@3 [ur] ^{\alpha}
	\ar@3 [drr] _{\beta}
&&&& {\figeps{uo}}
\\
&& {\scalefigeps{-2mm}{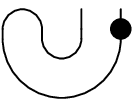}}
	\ar@3 [urr] _{\gamma} 
\ar@4 "1,3";"3,3" |-*+[o]{\scriptstyle \alpha\beta\left(\smallfigeps{u0i}\right)}
}
$$

\noindent Finally, for $\figeps{f} \: = \: \figeps{f0} \: \figeps{o} ^n \:$:
$$
\xymatrix{
{\raisebox{-1.25mm}{\begin{picture}(0,0)%
\includegraphics{n-of0o-u.pstex}%
\end{picture}%
\setlength{\unitlength}{4144sp}%
\begingroup\makeatletter\ifx\SetFigFont\undefined%
\gdef\SetFigFont#1#2#3#4#5{%
  \reset@font\fontsize{#1}{#2pt}%
  \fontfamily{#3}\fontseries{#4}\fontshape{#5}%
  \selectfont}%
\fi\endgroup%
\begin{picture}(665,558)(-952,200)
\put(-472,501){\makebox(0,0)[lb]{\smash{{\SetFigFont{8}{9.6}{\familydefault}{\mddefault}{\updefault}{\color[rgb]{0,0,0}$n$}%
}}}}
\end{picture}%
}}
	\ar @3 @/^6ex/ [rr] ^-{\alpha} _-{\hole}="1"
	\ar @3 @/_6ex/ [rr] _-{\beta} ^-{\hole}="2"
&& {\raisebox{-1.25mm}{\begin{picture}(0,0)%
\includegraphics{n-1f0oo-u.pstex}%
\end{picture}%
\setlength{\unitlength}{4144sp}%
\begingroup\makeatletter\ifx\SetFigFont\undefined%
\gdef\SetFigFont#1#2#3#4#5{%
  \reset@font\fontsize{#1}{#2pt}%
  \fontfamily{#3}\fontseries{#4}\fontshape{#5}%
  \selectfont}%
\fi\endgroup%
\begin{picture}(763,558)(-911,200)
\put(-472,501){\makebox(0,0)[lb]{\smash{{\SetFigFont{8}{9.6}{\familydefault}{\mddefault}{\updefault}{\color[rgb]{0,0,0}$n+1$}%
}}}}
\end{picture}%
}}
\ar@4 "1";"2" |-{\scriptstyle
  \alpha\beta\left(\figeps{f0}\:\figeps{o}^n\right)} 
}
$$

\subsubsection{Homotopy basis}

The $3$-polygraph $\Sigma$ is convergent and right-indexed. Thus, Theorem~\ref{Theorem:RightIndexedTDF} tells us that the following $4$-cells form a homotopy basis of $\Sigma^\top$:
$$
\gamma\delta, \: 
\delta\gamma, \:
\alpha\gamma, \: 
\beta\delta, \: 
\alpha\beta\left(\scalefigeps{-3mm}{u1n}\right), \:
\alpha\beta\left(\figeps{n0i}\right), \:
\alpha\beta\left(\figeps{u0i}\right),
$$

\noindent plus, for every $\figeps{f0}$ in $N_0$ and $n$ in $\Nb$, the $4$-cell
$$
\alpha\beta\left(\figeps{f0}\:\figeps{o}^n\right).
$$

\noindent In fact, the $4$-cells $\alpha\beta\left(\scalefigeps{-3mm}{u1n}\right)$, $\alpha\beta\left(\figeps{n0i}\right)$ and $\alpha\beta\left(\figeps{u0i}\right)$ are superfluous. Indeed, the $3$-spheres forming their boundaries are also the boundaries of $4$-cells of $\Sigma^\top(\alpha\gamma,\beta\delta)$, as diagrammatically proved thereafter.

\noindent For $\alpha\beta\left(\scalefigeps{-3mm}{u1n}\right)$:
$$
\xymatrix{
& {\figeps{nou}}
	\ar@3 [rr] ^{\beta}
	\ar@{} [dr] |{=}
&& {\figeps{nuo}}
	\ar@3 @/^6ex/ [ddrrr] ^{\gamma} 
	\ar@{} "2,5" ; "1,4" |{=}
\\
{\scalefigeps{-3mm}{nou-un}}
	\ar@3 [ur] ^{\delta}
	\ar@3 [rr] |*+[o]{\scriptstyle \beta} _(0.33){\hole}="1"
&& {\scalefigeps{-3mm}{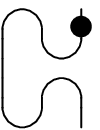}}
	\ar@3 [ur] |*+[o]{\scriptstyle \delta}
	\ar@3 [rr] |*+[o]{\scriptstyle \gamma}
&& {\figeps{oun}}
	\ar@3 [drr] |*+[o]{\scriptstyle \delta}
\\
{\scalefigeps{-3mm}{nu-o-un}}
	\ar@3 [u] ^{\alpha} 
	\ar@3 [d] _{\beta}	
	\ar@3  [urrrr] |*+[o]{\scriptstyle \gamma} ^(0.33){\hole}="2"
	\ar@3  [drrrr] |*+[o]{\scriptstyle \delta} _(0.33){\hole}="4"
	\ar@{} [rrrrrr] |(0.66){=}
&&&&&& {\figeps{o}}
\\
{\scalefigeps{-3mm}{nu-uon}}
	\ar@3 [rr] |*+[o]{\scriptstyle\alpha} ^(0.33){\hole}="3"
	\ar@3 [dr] _{\gamma}
&& {\scalefigeps{-3mm}{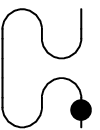}}
	\ar@3 [rr] |*+[o]{\scriptstyle \delta}
	\ar@3 [dr] |*+[o]{\scriptstyle \gamma}
	\ar@{} [dl] |{=}
&& {\figeps{onu}}
	\ar@3 [urr] |*+[o]{\scriptstyle \gamma}
\\
& {\figeps{uon}}
	\ar@3 [rr] _{\alpha}
&& {\figeps{uno}}
	\ar@3 @/_6ex/ [uurrr] _{\delta} 
	\ar@{} "4,5" ; "5,4" |{=}
\ar@4 "1";"2" |-*+[o]{\scriptstyle \alpha\gamma}
\ar@4 "3";"4" |-*+[o]{\scriptstyle \beta\delta}
}
$$

\noindent For $\alpha\beta\left(\figeps{n0i}\right)$:
$$
\xymatrix{
& {\scalefigeps{-2mm}{no-un}}
	\ar@3 [r] ^-{\delta} _-{}="3"
& {\figeps{no}}
\\
{\scalefigeps{-2mm}{n-oun}}
	\ar@3 [ur] ^{\alpha}
	\ar@3 [dr] _{\beta}
	\ar@3 [rrr] |-*+[o]{\scriptstyle\delta} _-{\hole}="2" ^-{}="4"
&&& {\figeps{on}}
	\ar@3 [ul] _{\alpha}
\\
& {\scalefigeps{-2mm}{n-uon}}
	\ar@3 [r] _-{\alpha} ^-{\hole}="1"
& {\scalefigeps{-2mm}{n-uno}}
	\ar@3 [ur] _{\delta}
\ar@4 "1"; "2" |-*+[o]{\scriptstyle\beta\delta}
\ar@{} "3"; "4" |-{=}
}
$$

\noindent And, for $\alpha\beta\left(\figeps{u0i}\right)\:$:
$$
\xymatrix{
& {\scalefigeps{-2mm}{nou-u}}
	\ar@3 [r] ^-{\beta} _-{\hole}="1"
& {\scalefigeps{-2mm}{nuo-u}}
	\ar@3 [dr] ^{\gamma}
\\
{\scalefigeps{-2mm}{onu-u}}
	\ar@3 [ur] ^{\alpha}
	\ar@3 [dr] _{\beta}
	\ar@3 [rrr] |-*+[o]{\scriptstyle\gamma} ^-{\hole}="2" _{}="3"
&&& {\figeps{ou}}
	\ar@3 [dl] ^{\beta}
\\
& {\scalefigeps{-2mm}{nu-uo}}
	\ar@3 [r] _-{\gamma} ^-{}="4"
& {\figeps{uo}}
\ar@4 "1" ; "2" |-*+[o]{\scriptstyle \alpha\gamma}
\ar@{} "3" ; "4" |-{=}
}
$$

\noindent We denote by $\Gamma_0$ the family made of the $4$-cells $\gamma\delta$, $\delta\gamma$, $\alpha\gamma$ and $\beta\delta$. Then, for every natural number $n$, one defines:
$$
\Gamma_{n+1} \:=\: \Gamma_n \:\amalg\: \left\{ \alpha\beta\left(\figeps{f0}\:\figeps{o}^n\right) \:,\: \figeps{f0} \in N_0 \right\} \:.
$$

\noindent Thus, the following set of $4$-cells is a homotopy basis of $\Sigma$:
$$
\Gamma \:=\: \bigcup_{n\in\Nb} \Gamma_n.
$$

\noindent For every natural number $n$, we denote by $\xi_n$ the $4$-cell $\alpha\beta\left(\figeps{o}^n\right)$ of $\Gamma_{n+1}$, hence of $\Gamma$.

\subsubsection{Lemma} \label{Lemma:ProofCounterExample0}

\begin{em}
Let $n$ be a natural number. There is no $4$-cell of $\Sigma^\top(\Gamma_0)$ with the same boundary as~$\xi_n$, \ie:
$$
s\xi_n \not\approx_{\Gamma_0} t\xi_n.
$$
\end{em}

\begin{proof}
Let us assume, on the contrary, that there exists a $4$-cell $\Phi$ in $\Sigma^\top(\Gamma_0)$ such that both $s\Phi=s\xi_n$ and $t\Phi=t\xi_n$ hold. We consider the derivation $d$ of $\Sigma^{\top}$ into the trivial module that takes the following values on the generating $3$-cells:  
$$
d(\alpha) = 1, \qquad d(\beta) = -1,  \qquad d(\gamma)=0, \qquad d(\delta)=0.
$$

\noindent Then, we check that, for any $4$-cell $\Psi$ of $\Sigma^\top(\Gamma_0)$, we have $d(s\Psi)=d(t\Psi)$. Since $d$ is a derivation, it is sufficient to check this equality on the generating $4$-cells of $\Gamma_0$: 
\begin{itemize}
\item $d(s\gamma\delta) = d(\gamma) = 0 \:$ and 
$\: d(t\gamma\delta) = d(\delta) = 0$.
\item $d(s\delta\gamma) = d(\delta) = 0 \:$ and 
$\: d(t\delta\gamma) = d(\gamma) = 0$.
\item $d(s\alpha\gamma) = d(\alpha) + d(\beta) + d(\gamma) =
  0\: $ and 
$\: d(t\alpha\gamma)=d(\gamma)=0$.
\item $d(s\beta\delta) = d(\beta) + d(\alpha) + d(\delta) =
  0\: $ and $\: d(t\beta\delta) = d(\delta) = 0$. 
\end{itemize}

\noindent Thus, since $\Phi$ is in $\Sigma^\top(\Gamma_0)$, one must have $d(s\Phi)=d(t\Phi)$. However, one has:
$$
d(s\Phi) \:=\: d(\alpha) \:=\: 1
\qquad\text{and}\qquad
d(t\Phi) \:=\: d(\beta) \:=\: -1.
$$

\noindent This proves that such a $4$-cell $\Phi$ cannot exist in $\Sigma^{\top}(\Gamma_0)$.
\end{proof}

\subsubsection{Lemma} \label{Lemma:ProofCounterExample1}

\begin{em}
Let $n$ be a natural number. There is no $4$-cell of $\Sigma^\top(\Gamma_n)$ with the same boundary as~$\xi_n$, \ie:
$$
s\xi_n \not\approx_{\Gamma_n} t\xi_n.
$$
\end{em}

\begin{proof}
On the contrary, let us assume that $\Phi$ is a $4$-cell of $\Sigma^\top(\Gamma_n)$ such that both $s\Phi=s\xi_n$ and $t\Phi=t\xi_n$ hold. As a direct consequence, we have:
$$
s_2\Phi \:=\: s_2\xi_n \:=\: \raisebox{-3mm}{\begin{picture}(0,0)%
\includegraphics{n-oon-u.pstex}%
\end{picture}%
\setlength{\unitlength}{4144sp}%
\begingroup\makeatletter\ifx\SetFigFont\undefined%
\gdef\SetFigFont#1#2#3#4#5{%
  \reset@font\fontsize{#1}{#2pt}%
  \fontfamily{#3}\fontseries{#4}\fontshape{#5}%
  \selectfont}%
\fi\endgroup%
\begin{picture}(485,380)(38,290)
\put(339,501){\makebox(0,0)[lb]{\smash{{\SetFigFont{8}{9.6}{\familydefault}{\mddefault}{\updefault}{\color[rgb]{0,0,0}$n$}%
}}}}
\end{picture}%
}.
$$

\noindent Hence, the normal form of $s_2\Phi$ is $\raisebox{-3mm}{\begin{picture}(0,0)%
\includegraphics{n-1on-u.pstex}%
\end{picture}%
\setlength{\unitlength}{4144sp}%
\begingroup\makeatletter\ifx\SetFigFont\undefined%
\gdef\SetFigFont#1#2#3#4#5{%
  \reset@font\fontsize{#1}{#2pt}%
  \fontfamily{#3}\fontseries{#4}\fontshape{#5}%
  \selectfont}%
\fi\endgroup%
\begin{picture}(584,380)(78,290)
\put(339,501){\makebox(0,0)[lb]{\smash{{\SetFigFont{8}{9.6}{\familydefault}{\mddefault}{\updefault}{\color[rgb]{0,0,0}$n+1$}%
}}}}
\end{picture}%
}$. Now, let us prove that~$\Phi$ cannot contain any occurrence of a generating $4$-cell $\alpha\beta\left(\figeps{f0}\:\figeps{o}^k\right)$ or its inverse, with $k<n$. If that was the case, there would exist $4$-cells $\Psi_1$, $\Psi_2$ in $\Sigma^\top(\Gamma_n)$, a context $C$ of $\Sigma^\top$, an $\epsilon$ in $\ens{-1,1}$, a $2$-cell~$\figeps{f0}$ and a $k$ in $\ens{0,\dots,n-1}$ such that the $4$-cell $\Phi$ decomposes this way: 
$$
\Phi \:=\: \Psi_1 \star_3 
C \left[ \alpha\beta\left(\figeps{f0}\:\figeps{o}^k\right)^{\epsilon} \right] 
\star_3 \Psi_2.
$$

\noindent As a consequence, we would have:
$$
s_2\Phi 
	\:=\: s_2 \left( C \left[ \alpha\beta\left(\figeps{f0}\:\figeps{o}^k\right)^{\epsilon} \right] \right)
	\:\approx_{\Sigma_3}\: (s_2 C) \left[ s_2 \alpha\beta\left(\figeps{f0}\:\figeps{o}^k\right) \right]
	\:=\: \raisebox{-8mm}{\input{c-n-ook-u-f0.pstex_t}} \:.
$$

\noindent Since $\Sigma$ is convergent, this implies that $s_2\Phi$ and the rightmost $2$-cell have the same normal form. One denotes by $D$ the context of $\Sigma_2^*$ such that $D[\ast]$ is the normal form of $(s_2 C)[\ast]$. Then, the following equality holds:
$$
\raisebox{-3mm}{\begin{picture}(0,0)%
\includegraphics{n-1on-u.pstex}%
\end{picture}%
\setlength{\unitlength}{4144sp}%
\begingroup\makeatletter\ifx\SetFigFont\undefined%
\gdef\SetFigFont#1#2#3#4#5{%
  \reset@font\fontsize{#1}{#2pt}%
  \fontfamily{#3}\fontseries{#4}\fontshape{#5}%
  \selectfont}%
\fi\endgroup%
\begin{picture}(584,380)(78,290)
\put(339,501){\makebox(0,0)[lb]{\smash{{\SetFigFont{8}{9.6}{\familydefault}{\mddefault}{\updefault}{\color[rgb]{0,0,0}$n+1$}%
}}}}
\end{picture}%
} \:=\: \raisebox{-8mm}{\input{c-n-1ok-u-f0.pstex_t}}.
$$

\noindent Let us prove that this is not possible. For that, we define the derivation $d$ of $\Sigma_2^*$ into the module $M_{X,\ast,G}$ given thereafter:
\begin{itemize}
\item The abelian group $G$ is freely generated by the set $\Nb$ of natural numbers. The natural number $n$, seen as a generator of $G$, is denoted by $a_n$. 

\item The $2$-functor $X:\Sigma^*_2\fl\Set$ is generated by the values:
$$
X\left(\:\figeps{1-cell}\:\right) \:=\: \Nb, \qquad 
X\left(\figeps{n}\right) \:=\: (0,0), \qquad 
X\left(\figeps{o}\right) (i) \:=\: i+1.
$$

\item The derivation $d$ is given by:
$$
d\left(\figeps{n}\right) \:=\: 0, \qquad
d\left(\figeps{u}\right) (i,j) \:=\: a_j, \qquad
d\left(\figeps{o}\right) (i) \:=\: 0.
$$
\end{itemize}

\noindent Then, on the one hand, we have:
$$
d\left( \raisebox{-3mm}{\begin{picture}(0,0)%
\includegraphics{n-1on-u.pstex}%
\end{picture}%
\setlength{\unitlength}{4144sp}%
\begingroup\makeatletter\ifx\SetFigFont\undefined%
\gdef\SetFigFont#1#2#3#4#5{%
  \reset@font\fontsize{#1}{#2pt}%
  \fontfamily{#3}\fontseries{#4}\fontshape{#5}%
  \selectfont}%
\fi\endgroup%
\begin{picture}(584,380)(78,290)
\put(339,501){\makebox(0,0)[lb]{\smash{{\SetFigFont{8}{9.6}{\familydefault}{\mddefault}{\updefault}{\color[rgb]{0,0,0}$n+1$}%
}}}}
\end{picture}%
} \right) \:=\: a_{n+1}.
$$

\noindent And, on the other hand, we use the fact that $d$ is a derivation to compute:
$$
d\left( \raisebox{-8mm}{\input{c-n-1ok-u-f0.pstex_t}} \right) 
	\:=\: d\left( \raisebox{-4mm}{\input{d-ast.pstex_t}} \right) 
		+ d\left( \raisebox{-3mm}{\begin{picture}(0,0)%
\includegraphics{n-1ok-u.pstex}%
\end{picture}%
\setlength{\unitlength}{4144sp}%
\begingroup\makeatletter\ifx\SetFigFont\undefined%
\gdef\SetFigFont#1#2#3#4#5{%
  \reset@font\fontsize{#1}{#2pt}%
  \fontfamily{#3}\fontseries{#4}\fontshape{#5}%
  \selectfont}%
\fi\endgroup%
\begin{picture}(584,380)(78,290)
\put(339,501){\makebox(0,0)[lb]{\smash{{\SetFigFont{8}{9.6}{\familydefault}{\mddefault}{\updefault}{\color[rgb]{0,0,0}$k+1$}%
}}}}
\end{picture}%
}\right) 
		+ d\left(\figeps{f0}\right) 
	\:=\: d(f) + a_{k+1} \:,
$$

\noindent where $f$ denotes $D\left[\figeps{f0}\right]$. Thus, we have $a_{n+1}=a_{k+1}+d(f)$, with $k<n$ and some $f$ in $\Sigma^*_2$. This is impossible because~$G$ is freely generated and $d$ sends any $2$-cell of $\Sigma^*_2$ to an element of $G$ written using the $a_i$'s with positive coefficients.

We conclude that the $4$-cell $\Phi$ is built from the $4$-cells of $\Gamma_0$ and their inverses only, \emph{i.e.} $\Phi$ is a $4$-cell of $\Sigma^\top(\Gamma_0)$. However, this would contradict Lemma~\ref{Lemma:ProofCounterExample0}.
\end{proof}

\subsubsection{Theorem}

\begin{em}
The $3$-polygraph $\Sigma$ does not have finite derivation type.
\end{em}

\begin{proof}
On the contrary, let us assume that $\Sigma$ does have finite derivation type. Then, by application of Proposition~\ref{PropExtractFiniteHomotopyBase}, there exists a finite subfamily $\Gamma'$ of $\Gamma$ which is a homotopy basis of $\Sigma^\top$.

Since $\Gamma'$ is finite, there exists some natural number $n$ such that $\Gamma'$ is contained in $\Gamma_n$. In particular, the $4$-cell $\xi_n$ is not in $\Gamma'$. However, since $\Gamma'$ is a homotopy basis and since $\Gamma'$ is contained in $\Gamma_n$, we have:
$$
s\xi_n \approx_{\Gamma_n} t\xi_n.
$$

\noindent We have seen in~\ref{Lemma:ProofCounterExample1} that this is not possible, thus contradicting the fact that one can extract a finite homotopy basis from $\Gamma$. As a consequence, the $3$-polygraph $\Sigma$ does not have finite derivation type.
\end{proof}

\subsubsection{A variant of the counterexample}

In the previous $3$-polygraph, one can think that the problem comes from the complicated normal forms, especially from the fact that one can find normal forms $\figeps{f0}$ of $N_0$ everywhere in a given $2$-cell. Here we give another example, similar to the first one but with more simple normal forms. It is a bit more contrived, which led us to prefer the other one for the main exposition.

Let $\Xi$ be the $3$-polygraph with the following generating cells:
\begin{itemize}
\item Two $0$-cells, denoted by $\xi$ and $\eta$ and, in the diagrammatic representations,  respectively pictured by a white background and by a gray one. 
\item Two $1$-cells $\onecell{p}{\xi}{\eta}$ and $\onecell{q}{\eta}{\xi}$. By abuse, both are pictured by a wire, leaving the backgrounds discriminate them.
\item Four $2$-cells $\figeps{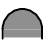}$, $\figeps{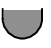}$, $\figeps{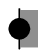}$, and $\figeps{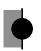}$.
\item Two $3$-cells $\:\figeps{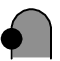} \: \overset{\alpha'}{\tfl} \: \figeps{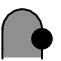}\:$ and $\figeps{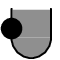} \:  \overset{\beta'}{\tfl} \:\figeps{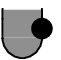}$.
\end{itemize}

\noindent Following the same reasoning steps as in the previous example, one proves that the finite $3$-polygraph~$\Xi$ is convergent. But it lacks finite indexation and finite derivation type. Indeed, the following family of $4$-cells, indexed by the natural number $n$, form a minimal homotopy basis of $\Xi^\top$:
$$
\xymatrix{
{ \raisebox{-2mm}{\begin{picture}(0,0)%
\includegraphics{n-oon-u-b.pstex}%
\end{picture}%
\setlength{\unitlength}{4144sp}%
\begingroup\makeatletter\ifx\SetFigFont\undefined%
\gdef\SetFigFont#1#2#3#4#5{%
  \reset@font\fontsize{#1}{#2pt}%
  \fontfamily{#3}\fontseries{#4}\fontshape{#5}%
  \selectfont}%
\fi\endgroup%
\begin{picture}(440,379)(38,290)
\put(294,524){\makebox(0,0)[lb]{\smash{{\SetFigFont{8}{9.6}{\familydefault}{\mddefault}{\updefault}{\color[rgb]{0,0,0}$n$}%
}}}}
\end{picture}%
} \!\!\!\!} 
	\ar @3 @/^7ex/ [rrr] ^-{\alpha'} _-{\hole}="1"
	\ar @3 @/_7ex/ [rrr] _-{\beta'} ^-{\hole}="2"
&&& { \raisebox{-2mm}{\begin{picture}(0,0)%
\includegraphics{n-1on-b.pstex}%
\end{picture}%
\setlength{\unitlength}{4144sp}%
\begingroup\makeatletter\ifx\SetFigFont\undefined%
\gdef\SetFigFont#1#2#3#4#5{%
  \reset@font\fontsize{#1}{#2pt}%
  \fontfamily{#3}\fontseries{#4}\fontshape{#5}%
  \selectfont}%
\fi\endgroup%
\begin{picture}(539,379)(78,290)
\put(294,524){\makebox(0,0)[lb]{\smash{{\SetFigFont{8}{9.6}{\familydefault}{\mddefault}{\updefault}{\color[rgb]{0,0,0}$n+1$}%
}}}}
\end{picture}%
} } 
\ar @4 "1";"2" |-*+[o]{\scriptstyle \alpha'\beta'\left(\figeps{o-c}^{n}\right)}
}
$$

\begin{small}
\renewcommand{\refname}{\textbf{References}}
\bibliographystyle{amsplain}
\bibliography{../bibliographie}
\end{small}

\end{document}